\documentclass[a4paper,12pt]{amsart}

\usepackage[X2,T2A]{fontenc}
\usepackage[utf8]{inputenc}
\usepackage[russian,english]{babel}

\usepackage{geometry}
\usepackage[parfill]{parskip}
\usepackage{graphicx}
\usepackage{amssymb}
\usepackage{epstopdf}
\usepackage[toc]{appendix}
\usepackage{csquotes}
\usepackage{hyperref}
\usepackage{footmisc}
\usepackage{amsmath} 
\usepackage{amsthm}
\usepackage{tikz}
\usetikzlibrary{calc}

\DeclareTextSymbolDefault{\cyryat}{X2}

\title{Lobachevsky's Views of Geometry}
\author{Andrei Rodin}

\begin{document}
\selectlanguage{english}                                          
\begin{abstract}
Lobachevsky's discovery of Non-Euclidean geometry was a byproduct of his more ambitious plan of reforming the traditional Euclid-style foundations of geometry. This project aligned with d'Alembert’s and Condillac's practically-oriented visions of mathematics and  aimed at closing the perceived gap between the pure mathematics, on the one hand, and the natural sciences and engineering, on the other hand. Following d'Alembert Lobachevsky deliberately avoided the axiomatic development of his geometrical theory and combined instead some preliminary intuitive synthetic constructions with their rigorous analytic treatment. In addition to Non-Euclidean geometry, Lobachevsky’s project brought about an early form of Dimension theory.

The paper comprises a description of the local context of Lobachevsky's works, including some pointers to Russian geometry textbooks of the beginning of the 19th century. We summarise Lobachevsky's philosophical views which motivated his mathematical achievements and show how these motivations were misunderstood and misinterpreted in the beginning of the 20th century by Henri Poincar\'e and Ernst Cassirer. Finally we argue that Lobachevsky’s epistemic views remain pertinent in the context of the ongoing discussion about the reunion of mathematics and physics. 

The \emph{Appendix} comprises English translations of eight  Lobachevsky's documents of various length and character where he expresses his general epistemic views on mathematics and, in particular, geometry.

\end{abstract}    

\maketitle
\newpage

\tableofcontents

\section{Introduction}\label{intro}
In his popular 1961 paper titled \emph{The Revolution of Mathematics} prominent American mathematician Marshall Stone describes what he sees as a turning point of the history of his discipline:

\begin{quote}
While several important changes have taken place since 1900 in our conception of mathematics or in our points of view concerning it, the one which truly involves a revolution in ideas is the discovery that mathematics is entirely independent of the physical world. \cite[p.716]{Stone:1961}
\end{quote}

Stone continues explaining the character of modern mathematics by  pointing to its historical origin:

\begin{quote}
The discovery from which our current view of mathematics as a totally abstract, strictly logical, and entirely independent discipline has emerged was the discovery of non-euclidean plane geometry by the two Bolyais, Gauss, and Lobachevski early in the nineteenth century.  \cite[p.718]{Stone:1961}
\end{quote}

In the early 1960s, the view according to which the discovery of non-euclidean geometries was a pivotal moment of history when the development of mathematics first turned towards its contemporary state of ``entirely independent'' and ``strictly logical'' discipline has been also defended by Jean Dieudonn\'e \cite{Dieudonne:1961} and other proponents of then-new Structuralist Bourbaki-style \emph{architecture of mathematics} \cite{Bourbaki:1950}. In spite of certain loud disagreeing voices including Richard Feynman's \cite{Feynman:1965}, in the second half of the 20th century the Structuralist view of Stone and Dieudonn\'e became a mainstream. A critically-minded historian could, of course, immediately suspect that Stone's and  Dieudonn\'e's historical narratives were heavily teleological, i.e., intendedly designed to demonstrate how the ongoing project of rebuilding the foundations of mathematics pursued by its proponents was deeply rooted in the history of this discipline. But as far as such a historical narrative does not contradict facts it can be anyway accepted \textemdash\  at least as a legitimate interpretation of the history.  It turns out, however, that the real history of the 19th century geometry is drastically different and more surprising than the received historical narrative might suggest. The present paper describes general views on mathematics (and on geometry, in particular)  held by of one the pioneers of the new geometry, namely, by Nikolai Ivanovitch Lobachevsky. Just like Stone and Dieudonn\'e Lobachevsky held strong views about his contemporary mathematics and aimed at rebuilding its conceptual foundations. But Lobachevsky's project aimed, in fact, at the opposite direction.

In the nutshell, Lobachevsky's project of rebuilding the foundations of geometry aims at closing the gap between this mathematical discipline and natural sciences. In his view, traditional Euclid-style geometrical concepts like that of point or a line are too abstract and too remote from the experience and the practice of measurements, and for this reason need to be either replaced or given a wholly new form, which would help to better connect geometry to the natural sciences and the engineering. It was this radically naturalistic and practically-oriented conception of geometry as an art of measuring that eventually led Lobachevsky to his great mathematical discovery. Saying that Lobachevsky somehow preconceived the view on mathematics latter held by Stone and Dieudonn\'e is not a matter of interpretation because these explicitly contradict the available evidence, as we shall shortly see. 


It is important to reconstruct Lobachevsky's views on geometry accurately not only for better understanding the history but also for learning how to think about non-Euclidean geometries in a more application-friendly manner. Richard Feynman \cite{Feynman:1965}, Vladimir Arnold \cite{Arnold:1999}, Frieman Dyson \cite{Dyson:1972}, Vladimir Voevodsky \cite{Rodin:2021} and all other people who believe, notwithstanding the Structuralist critique, that the connection to natural science is essential for mathematics, might find in Lobachevsky an inspiring example of how a naturalistic stance in mathematics can be combined with a genuine mathematical inventiveness and intellectual freedom. 

Lobachevsky's views on mathematics were already discussed in some earlier literature, mostly published in Russian during the Soviet period (1922-1991). The first scholar who systematically studied Lobachevsky's general views on mathematics and science was, to the best of our knowledge, Nicolai Nicolaevitch Parfentiev (1877-1943) who was professor of mathematics in Kazan University and along with his supervisor Alexandre Vasilievitch Vasiliev (1853-1929)  worked on establishing Lobachevsky's legacy.\footnote{Vasiliev's \emph{opus magnus} was his biography Lobachevsky; it was accomplished and printed in 1927. In circumstances that remain so far only poorly known and understood the book was not openly sold but instead kept in stock for two years and in 1929 finally destroyed. By a lucky chance, a copy of proofs of Vasiliev's book survived and allowed for new editions of this important source in 1992 \cite{Vasiliev:1992} and 2022 \cite{Vasiliev:2022}. } In the late 1920s Parfentiev published a series of papers on what he called \emph{Lobachevsky's Philosophy of Nature}, see  \cite{Parfentiev:1930} and further references therein. In 1950-1951 Sofya Alexandrovna Yanovskaya  (1896-1966) published in two instalments a long paper titled \emph{Lobachevsky's Worldview}  \cite{Yanovskaya:1950a}, \cite{Yanovskaya:1951}. This paper contains a lot of valuable material including a critique of historical narratives concerning the discovery of non-Euclidean geometry found in the contemporary historical literature. It should be noted, however, that this Yanovskaya's work is written and published during the period when Stalin's ideological propaganda and censorship were particularly strong.  As a result Yanovskaya's paper often more resembles an ideological pamphlet than a piece of scholarship. Understandably, Soviet scholars of later generations typically tried to avoid the politically sensitive topics concerning Lobachevsky's philosophical leanings, and focused instead on his mathematics and on biographical details. It continued until 2007 when Vasily Yakovlevitch Perminov (1938-2024) published a paper where he made an interesting attempt to place Lobachevsky's philosophical views in the history of philosophy and in today's philosophy of mathematics \cite{Perminov:2007}. Among the Lobachevsky scholarship outside Russia I would like to mention a 1975 paper by Norman Daniels which presents implications of Lobachevsky's views on the role of mathematics in physics \cite{Daniels:1975}.

The rest of this work is organised as follows. In Section \ref{context} we describe the academic and intellectual context in which Lobachevsky pursued his research. Section \ref{timeline} provides a timeline of Lobachevsky's mathematical achievements with a natural emphasis on his work in the foundations of geometry. In the same Section precise references to the discussed Lobachevsky's publications are provided. In the following Section \ref{sources} we describe Lobachevsky's writings that have been used in this work for reconstructing the author's philosophical views on mathematics. It includes fragments of Lobachevsky's mathematical publication as well as his pedagogical notes and some other selected texts. All of these texts have been earlier published by Russian Lobachevsky's scholars, so the present work involves no new archival findings. The following Section \ref{summary} contains an analytic  summary of Lobachevsky's views based on our reading of the aforementioned sources, and in the next Section \ref{hilb} we once again compare Lobachevsky's understanding of his own mathematical work with its later axiomatic interpretation. In the concluding Section \ref{future} we argue that Lobachevsky's philosophical views are pertinent for today's mathematical research. The Appendix \ref{ap} to the present work comprises full English translations of  height Lobachevsky's documents listed in the \emph{Sources} Section  \ref{sources}. These translations from Russian are made by the author of the present work.

\section{The Context: Euclid, d'Alembert, Condillac and Teaching Geometry in the Russian Empire} \label{context}

Lobachevsky's achievements are usually analysed by historians of mathematics against the background of the earlier and the contemporary geometrical research pursued in Europe \cite{Houzel:1992}. In this Section we describe a more specific context of Lobachevsky's work that primarily concerns teaching mathematics rather than mathematical research. 

Nikolai Lobachevsky was born December 1, 1792\footnote{which was November 20, 1792 by the Julian Calendar used in Russia before 1918} in Nizhny Novgorod. In the age of 10 he was brought by his mother to Kazan and after an examination admitted (along with his brothers Alexei and Alexander) to the local Gymnasium on a State funding. In 1807 Nikolai Lobachevsky graduated from the Gymnasium and entered the Imperial Kazan University, which was formally founded only three years earlier as an extension of the same Gymnasium. In 1811 Lobachevsky graduated from Kazan University  with a Master degree and was immediately offered a teaching position in the same university. In 1814 Lobachevsky's mentor Johann Christian Martin Bartels (1769-1836) was appointed as the Dean of the newly founded Department of Mathematics and Physics, and promoted Lobachevsky to the rank of Adjunct Professor. The rest of Lobachevsky's career was made in the same university. In 1819 Lobachevsky replaced Bartels as the Dean of the Department, and in 1827 he was appointed as the Rector of the University; this latter position Lobachevsky held until 1846.\footnote{These and further biographical details of Lobachevsky are found in \cite{Kagan:1948},\cite{Kagan:1959}. It should be stressed that the first period of Lobachevsky's life before moving to Kazan still remains poorly known. Even Lobachevsky's date and place of birth were firmly established only in the 1950s \cite{Andronov:1956}, \cite{Polotovsky:2007}. There is a compelling evidence that Nikolai Lobachevsky's  legal father Ivan Lobachevsky who gave him the name was neither his biological farther nor he made any contribution to his education. The biological father of Nikolai Lobachevsky's  (and of his two brothers) \textendash\ and one who financially and morally supported the family \textendash\ was likely Sergei Stepanovitch Shebarshin (\selectlanguage{Russian} Сергей Степанович Шебаршин \selectlanguage{English}), an officer and land surveyor, who was a colleague of Ivan Lobachevsky \cite{Gudkov:1992}. These findings help to explain why during his life time Nikolai Lobachevsky systematically concealed  his family background and even his birthdate, so this information, rather surprisingly, is not found in the official university documentation. It is tempting to speculate that the art of land surveying might contribute to the formation of Lobachevsky's view of geometry in his childhood through the occupation of his father  (in fact, of both of his fathers). Some further biographical details of Nikolai Lobachevsky are found below in the main text in Section \ref{timeline}. }

As a university professor based outside the two Russian capitals Lobachevsky was hardly expected to do an original research. He was expected, primarily, to teach mathematics and produce usable teaching materials including textbooks. By the beginning of the 19th century mathematics and sciences in Russian Empire were already well institutionalised \cite{Yushkevich:1968}. Research was supposed to be done in the \emph{Academy of Sciences} founded in Saint-Petersburg by tsar Peter the First in 1724 \footnote{ \emph{Russian Academy of Sciences} is today's name of this institution. Since its foundation in 1724 it was many times reorganised and renamed, see \url{https://www.ras.ru/about/history.aspx} for an overview.}, while teaching mathematics and sciences was left to universities, military schools and gymnasia. The Academy also supervised teaching at the university level and was responsible for reviewing textbooks produced by university professors like Lobachevsky; eventually the Academy produced some teaching materials of the university level itself (which is the case of the influential textbooks in Differential and Integral Calculi by Leonhard Euler \cite{Euler:1755}, \cite{Euler:1768-1770}). The supervision of teaching mathematics and other subjects in gymnasia was the responsibility of local universities. The Academy also  supervised research, which was eventually done by some university  professors.\footnote{This is why in 1832 Lobachevsky's paper was reviewed by  academician Ostrogradsky, see \ref{sfoundations}. The separation of teaching and research at the institutional level is still present today's Russia in spite of recurrent administrative efforts aiming at their integration.} 

Lobachevsky's main mathematical discovery was a byproduct of his project aimed at producing a novel geometry textbook of the university level. In order to understand this project in its integrity it is essential to take into account the situation concerning teaching mathematics in Europe and in Russia during the decades preceding the beginning of Lobachevsky's career, i.e., the early 1810s.  

In 1757 Jean Le Rond d'Alembert in his and Diderot's \emph{Encyclopaedia} included the article  \emph{Geometry} \cite[vol. 7 (1757), p. 629A-639B]{Diderot&DAlembert:1751-65} which was a manifesto aiming at radical revision of the traditional way of teaching geometry based on Euclid's \emph{El\'ements} and its modernised versions. In particular, d'Alembert argues that geometry textbooks should begin with treating solids rather than points, lines and other abstract entities; he very dismissively talks about Euclid and the very idea of axiomatics \footnote{For a review of d'Alembert's ideas about mathematics and the recent secondary literature, see \cite{Lamande:2019}}. In the second half of the 18th century d'Alembert's program of reforming mathematical education and the re-orienting mathematics research towards practical applications was shared by many influential and powerful people in Europe who realised that mathematics was a crucial factor in gaining a military might and strengthening the imperial power. These developments led to the publication of a significant number of new practically-oriented mathematics textbooks. Here are  some examples.   

Duke de Choiseul (1719-1785), the Minister of Navy of French monarch of the time Louis XV, clearly saw that Navy officers and Navy engineers needed serious mathematical competence. In 1763 he appointed as the \emph{examiner of naval guards} mathematician \'Etienne B\'ezout (1730-1783) who in 1764-1767 produced   \emph{Course of mathematics for Guardians of the Flag and the Navy} in four parts; the second part containing the Elementary geometry was first published in 1765 \cite{Bezout:1765}. This practically-oriented general course by B\'ezout  was extremely successful as a textbook for military and engineering schools. It was many times re-edited (including adopted editions holding different titles) and translated in English, Portuguese and Russian. 

Another very popular geometry textbook of the new generation was composed by Adrien-Marie Legendre (1752-1833) and first published in 1794 \cite{Legendre:1794}. Legendre began his professional career in the Paris Military School (L'\'Ecole militaire de Paris) being appointed to a position of mathematics professor in this school by the recommendation of d'Alembert \cite[p. 392]{Ball:1908}. Interestingly, Legendre made his most remembered contribution to mathematics education not as a part of his professional obligations but as a private initiative during the two years of  the post-revolutionary  \emph{Terror} (1793-1794) when he was hiding in Paris in a private house. Legendre's \emph{Elements} were widely used for teaching Elementary geometry in Europe, Russia and elsewhere during the most of the following century. Unlike B\'ezout's \emph{Course} Legendre's \emph{Elements} was an all-purpose textbook suitable for general education. It preserved some traditional features of Euclid's \emph{Elements} but also involved significant modifications of Euclid's theoretical structure in the accordance with  d'Alembert programme. In particular, it emphasised the role of measurement \footnote{The opening sentence of Legendre's textbook \cite{Legendre:1794} reads ``Geometry is a science that studies the measurement of extension'' (my translation from French). As we shall see, Lobachevsky shares the same conception of geometry, which remains quite stable throughout his career.} and applied throughout the textbook the elementary symbolic algebra, which enabled computations. 

A more advanced course of Elementary geometry was published in 1803 by Sylvestre-François Lacroix \cite{Lacroix:1803}.  Lacroix follows  d'Alembert's program more closely than B\'ezout and Legendre: in his presentation of Elementary geometry he mentions no axiom, he begins his course by introducing solids and he more systematically than the two other authors applies analytic computational methods. In addition to this  general purpose geometry textbook  Lacroix authored more specific textbooks of Descriptive geometry and Trigonometry. 

The 18th century Imperial Russia\footnote{Russia officially declared its imperial status by Peter the First in 1721 as a result of its victory in the Great Northern War (1700-1721).} did not stay apart of these developments. Tsar Peter the First (reign 1682-1725) early in his reign realised the importance of mathematics eduction for Russia's military building, and he personally supported the creation and publication in 1703 of the first Russian Elementary Arithmetic textbook composed by Leonty Magnitski\footnote{``Magnitski'' is a nickname given to Leonty by Peter in 1700, which meant to tell that the man was attractive like a magnet.  Leonty Magnitski (1669-1739) was born as  Leonty Telyashin into a peasant family and made an outstanding career as a private teacher for powerful families in Moscow, which allowed him circa 1698 to meet Peter in person, see \cite{BRE:2017}.} on the basis of a wide range of European mathematical literature existing to the date. This textbook, among other things, first introduced in Russia the positional Hindu-Arabic numerals but also treated more advanced topics including logarithms. As a part of the same project in 1701 Peter created in Moscow a \emph{School of Mathematics and Navy} where he appointed Magnitsky and other teaching staff about a half of which were invited foreign nationals \footnote{For more details on Magnitski's \emph{Arithmetic} and the \emph{School of Mathematics and Navy} see \cite{Yushkevich:1968}}. 

The first original textbook on Elementary geometry was written and published in Russian in 1765 by Nikolai Gavrilovich Kurganov  (1726-1796), who earlier graduated from the Moscow School of Mathematics and Navy \cite{Kurganov:1765}\footnote{Earlier in 1757 the same author published a textbook on \emph{Universal Arithmetic}, i.e. Algebra \cite{Kurganov:1757}. Another relevant Russian mathematical textbook published during the same historical period is \emph{Abridged Mathematics} by Stepan Yakovlevich Rumovsky \cite{Rumovsky:1760} who was Leonhard Euler's student and later himself an academician.
For more details see  \cite{Yushkevich:1968}}. Like B\'ezout's geometry textbook  \cite{Bezout:1765}, which  appeared in Paris the same year, Kurganov's  \emph{General Geometry} was written for teaching future navy officers. Kurganov's geometry textbook is equally d'Alembertian in its content, its composition, and its purpose. 

It should be added that all the aforementioned French textbooks were at some point translated and published (typically more than once) in Russian, see \cite{Bezout:1794}, \cite{Bezout:1798}, \cite{Legendre:1819}, \cite{Lacroix:1822}, \cite{Lacroix:1835}. The role of linguistic barrier should not be overestimated in this case: in the 18th and 19th centuries Russian mathematics students as well as their professors were typically fluent both in German and in French. The significant number of mathematics textbooks published in Russian in the 18th and 19th centuries  was a result of policy aiming at the development of Russian mathematical terminology and promote using the Russian language in mathematics and science education.\footnote{A somewhat extreme example of the implementation of this policy is Russian translation of Euclid's \emph{Elements} published in 1784 by two navy officers Vassily Nikitin (1737-1809) and Prokhor Suvorov (1750-1815) \cite{Euclid:1784} for teaching purposes. The translators made a systematic attempt to replace Cyrillic transliterations of Latin, French and German mathematical terms, which were earlier conventionally used in Russian mathematical literature, by proper Russian words of Slavonic origin. Only few of Nikitin's and Suvorov's terminological inventions survived by the present date but some of them (indcluding the standard Russian equivalent ``\selectlanguage{Russian}точка \selectlanguage{English}'' for English ``point'') did so.}

The above overview demonstrates a specific difference between the situation in mathematical education in the beginning of the 19th century in France, on the one hand, and in Russia, on the other hand. At this point of its history Russian Empire managed to establish a network of military and engineering schools where mathematics was taught in the modern d'Alembertian style very much like in the contemporary France. But the historical background in the two cases was not the same. In France the new wave of d'Alembertian practically-oriented mathematics textbooks emerged against an earlier well-established tradition of learning and teaching geometry after Euclid's \emph{Elements} and its later adaptations. No comparable tradition existed in Russia before Magnitski's \emph{Arithmetic} was published in 1703. So while in France the new d'Alembertian approach in teaching geometry had to deal with the resistance of the earlier teaching tradition, no similar resistance was present in Russia where mathematics education before Magnitski  barely existed, and the d'Alembertian practically-oriented approach was promoted by the highest authorities to begin with. 

This does not mean, however, that in the 18th century Russia the controversy between the traditional Euclide-style approach and the ``revolutionary'' d'Alembertian approach in teaching geometry remained wholly unknown. Contemporary adopted editions of Euclid's \emph{Elements} published in European languages were read in Russia along with the more recent d'Alembert-style mathematics textbooks. The first recorded translation of such a source into Russian was published in 1739 \cite{Euclid:1739}; this ``first Russian Euclid'' was a translation from Latin of \emph{Euclid's Elements of Plane and Solid Geometry} by Andr\'e Tacquet (1612-1660) first published in 1654 \cite{Tacquet:1654}, and then many times republished. The translation into Russian was made by Ivan Satarov, a military physician serving in the Admiralty located in Saint-Petersburg \cite{Rybnikov:1941}, \cite{Lokot:2018}. The title page of the book holds the name of its editor Henry Farquharson (1674-1739), a Scottish mathematician hired by Peter back in 1699 and later associated with Magnitski; it is likely that Farquharson was responsible for the choice of Tacquet's textbook for this translation.

By 1820 there were already available four different Russian translations of Euclid's  \emph{Elements}. This list included (i) the aforementioned translation by Satarov of 1739 \cite{Euclid:1739},  (ii) translation by Nikolai Kurganov of 1769 \cite{Euclid:1769} made from French edition \cite{Euclide:1762}, (iii) translation by Prokhor Suvorov and Vassily Nikitin of 1784 \cite{Euclid:1784} from Gregory edition of 1703  \cite{Euclid:1703}, and (iv) translation by Foma (Thoma) Ivanovitch Petrushevsky of 1819 \cite{Euclid:1819} equally made from the Gregory edition. Among these people Foma Petrushevsky hold a strong maverick opinion according to which Euclid \emph{Elements} had to be taught in schools in its original non-modernised form. But as an ordinary  school teacher and later low-level educational administrator he had very little influence on educational matters which were controlled in Russia by the State via the \emph{Academy}. And the \emph{Academy} at this point of history was dominated by progressive mathematicians who did not support Petrushevsky's conservative view on mathematics education.  So  Petrushevsky's ideas remained very unpopular during his lifetime. 

Earlier in 1798, however, academician Semyon Emelyanovich Guriev (1766-1813), published a monograph titled \emph{Essay on Improvement of Element of Geometry}, where he tried to combine the rigour of traditional Euclid-style geometric reasoning with practical advantages of d'Alembertian analytic approach. The main point of Guriev's critique of recent d'Alembert-style geometry textbooks was a poor foundation of Differential and Integral Calculus, which he considered as a proper part of geometry. Since we are talking here about the state of Calculus before Cauchy and Weierstrass who in the second half of the 19th century provided Calculus with a satisfactory foundation \cite{Grabiner:1981} Guriev's critique from today's vintage viewpoint appears to be quite justified. In 1811 Guriev published his \emph{Foundation of Geometry} which was supposed to satisfy his methodological requirements. This geometry textbook remained d'Alembertian in its structure \textendash in particular, it started with solids in the three-dimensional space rather than lines and points \textendash but at the same time it made a special emphasis on foundations of differential and integral Calculus and for that matter relied on the traditional Euclid-style geometrical reasoning. Thus by the beginning of 19th century the controversy between the ``old'' and the ``new'' style of teaching geometry was realised and reflected upon in Russia in spite of the fact that the ``old'' Euclidean tradition was not deeply rooted in the history Russian mathematics education.   

Another influential predecessor of Lobachevsky who needs to be mentioned here is Timofey Fedorovich Osipovsky (1766-1832).  Osipovski started his teaching career in a primary school but in1803 became a university professor in Kharkiv University and later the Rector of this university \cite{Bakhmutskaya:1952}, \cite{Prudnikov:1952}. \footnote{Noticeably, in 1803 Osipovski was also offered a position in the \emph{Academy} but turned it down and opted for the university.} This spectacular promotion was due to the set of two mathematical textbooks composed and published by Osipovsky in 1801-1802, which immediately became very popular in Russia.  In 1801 appeared Osipovsky's volume on geometry \cite{Osipovski:1801} (which the author decided to call ``volume 2'' of his work), and in 1802 appeared his volume on Arithmetic  \cite{Osipovski:1802}; later in 1823 Osipvosky published the third volume of the series on the Theory of Analytic Functions, which included Differential and Integral Calculi and elements of Differential Equations and Calculus of Variations \cite{Osipovski:1823}.\footnote{For more information on mathematics education in the 18th century and early 19th century Russia see the introductory article to Lobachevsky's \emph{Geometry} by V.F. Kagan in \cite[vol.2, pp.9-29]{Lobachevsky:1946-51}. }

In additional to his work as a mathematical educator Osipovsky developed serious philosophical interests. He delivered public philosophical lectures in Kharkov University containing a systematic critique of Kant's views on space, time, and  causality  \cite{Osipovski:1807},  \cite{Osipovski:1813}\footnote{These lectures were delivered in 1807 and 1813 correspondingly}. Arguing against Immanuel Kant, Osipovsky supported a realistic and empirical understanding of these categories, which drew on the philosophical works of Denis Diderot, Jean Le Rond d'Alembert, \'Etienne Bonnot de Condillac and perhaps some other contemporary French authors. In 1805 in Moscow appeared  Osipovsky's translation of Condillac's \emph{Logic} \cite{Condillac:1780}, \cite{Condillac:1805}.  Following his French teachers Osipovsky developed his philosophical interests in a close connection to his mathematical interests. In particular, in his geometrical textbook of 1801 Osipovsky closely follows the plan suggested by d'Alembert in the \emph{Encyclopaedia} \cite[vol. 7 (1757), p. 629A-639B]{Diderot&DAlembert:1751-65}; his critique of Kant heavily draws on physical and mathematical ideas and examples.    

Although explicit references to other authors found in Lobachevsky's writings are very limited both in their number and their scope, there is no doubt that Lobachevsky was well familiar with all the aforementioned literature and developments. In the controversy about the ``old'' (i.e., Euclid-style) and the ``new'' (d'Alembert-style) mathematics textbooks Lobachevsky definitely takes the progressive d'Alembertian side. Following Osipovsky and d'Alembert,  Lobachevsky avoids to organise his mathematical writings axiomatically. He ignores not only the distinction between an axiom and a theorem but even a more basic distinction between a statement of a theorem and its proof.\footnote{About the influence of d'Alembert on Lobachevsky see \cite{Kagan:1948} and \cite{Houzel:1992}.} 

Lobachevsky's and Osipovsky's mathematical writings share not only their general d'Alembertian approach but also more specific features. The controversy between Analysis and Synthesis, which plays an important role in Lobachevsky, as we shall shortly see in \ref{sum5}, is also discussed in the \emph{Introduction} to Ospovsky's geometry textbook   \cite[1-2]{Osipovski:1801}; the same can be said about Osipovsky's treatment of the three spatial dimensions, which in a different form reappears in Lobachevsky \ref{sum4}. Lobachevsky's ideas concerning human senses as the primary source of human knowledge and the symbolic language as a tool of reasoning \ref{sum2}, \ref{sum6} are equally found in Osipovsky's philosophical writings. It appears that these ideas are borrowed by the two Russian mathematicians from the same source, namely, from Condillac's philosophical works, primarily from his \emph{Treatise on the Sensations} \cite{Condillac:1754}. Osipovsky's mathematical textbooks as well as his philosophical articles could motivate Lobachevsky to write his own geometry textbook in the d'Alembertian spirit. The impact of Osipovsky's works on Lobachevsky is definitely worth a further study.  

It should be stressed that the views on mathematics, science, knowledge, and other general issues shared by Osipovsky and Lobachevsky were hardly very original. At the times of the \emph{Industrial Revolution} similar views were shared by the majority of Lobachevsky's contemporaries in Russia, Germany and elsewhere, who were impressed and inspired by the growing impact of mathematics and mathematically-laden physics on engineering. What distinguished Lobachevsky among his many contemporaries was his mathematical genius, which allowed him to find so non-trivial application of those common philosophical ideas in mathematics. 

Not all of Lobachevsky's mathematical projects motivated by these ideas were successful, however. For example, his attempts to develop Dimension theory that could explain away geometrical abstractions (like the concepts of surface, line and point) in terms of concrete measurements \ref{sum4} had a very limited success. Seen retrospectively, they are very interesting but in the real history they remained almost unknown to the community did not have any impact on the later development of the Dimension theory \cite{Johnson:2025}. 

The idea that Euclid's Fifth Postulate (the Parallel Axiom) may fail was put forward and further developed by Lobachevsky in during the first half of the 19th century. Today we commonty qualify it as a great mathematical discovery. But Dieudonn\'e, Stone and other Mathematical Structuralists, who in the 1950s and 1960s developed their own peculiar understanding  of this idea, wrongly assumed that Lobachevsky thought about it in similar terms. Lobachevsky developed this idea in the in the same empirically- and practically-oriented conceptual framework centred on the concept of measurement, in which he attempted to develop his theory of spatial dimension. Lobachevsky convinced himself that at larger astronomical scales the Parallel Axiom cannot be justified by the everyday spatial intuitions and measurements of for this reason should be abandoned. What makes Lobachevsky a true creator \footnote{Calling Lobachevsky ``a creator'' rather than ``the creator'' of Hyperbolic geometry we take a distance from the continuing priority dispute. J\'anos Bolyai made the same discovery about the same time as Lobachevsky; remarkably he used a similar analytic technique. Friedrich Gauss, according to his own testimony, understood the core of Hyperbolic geometry back in 1790s. What allows us to avoid the priority dispute and this broader context is the fact that by all evidence Lobachevsky made his discovery independently of those authors. So Bolyai does not appear in our account, and Gauss appears only as a reader and eventual correspondent of Lobachevsky after Lobachevsky's main discoveries were already made \ref{new}.} of the new geometrical theory is the combination of these two factors: (i) his conceptualisation of Hyperbolic space as a possible model of physical space at large astronomical distances and (ii) Lobachevsky's invention of a new analytic calculus, which we call today the Hyperbolic Trigonometry. The axiomatic thinking played no role in Lobachevsky's mathematical work.

In what follows we focus on the conceptual part of Lobachevsky's mathematical thinking leaving aside many technical details.\footnote{For technical aspects of Lobachevsky's theory of Hyperbolic geometry see \cite{Braver:2011}.}

\section{The Timeline} \label{timeline}
In this Section we list, describe and briefly discuss Lobachevsky's geometrical works in the chronological order of their accomplishment. This list includes all published items as well as some unpublished items that have been lost. A similar list can be found in \cite[p. 290-297]{Lobachevsky:2010} but our list is more complete and contains our commentaries on various aspects of Lobachevsky's work.  In the end of this list we included some of Lobachevsky's mathematical publications on different subjects. The full list of all Lobachevsky's mathematical publications which appeared during his lifetime can be found in \cite{Houel:1870}. Translations and re-editions of Lobachevsky's geometrical works existing to the date are not fully covered in the present work but some of them are mentioned in the following commentaries. In order to identify Lobachevsky's writings, we, by default, use his \emph{Complete Collected Works} published back in 1946-51 \cite{Lobachevsky:1946-51} where all extant Lobachevsky's geometrical works have been re-edited and provided with extensive commentaries. 

\subsection{ \emph{Foundation of Geometry} of 1819 and  \emph{Geometry} of 1823}\label{fg1819-23}

In 1819 Lobachevsky submitted to Kazan University, with a hope of eventual publication, a manuscript titled \emph{Foundation of Geometry} (\selectlanguage{Russian}  \emph{Основание геометрии} \selectlanguage{English}). The manuscript has not been published and is now lost.  \cite[vol.2, p.125]{Lobachevsky:1946-51}

In 1823, responding to a call by the university administration that encouraged professors to publish their original textbooks \cite{Crowe:1995} Lobachevsky submitted another manuscript titled \emph{Geometry} for publication using government funds.\footnote{In 1822 Lobachevsky was appointed full professor in the Kazan University.} Most likely the lost 1819 manuscript was an earlier version of the 1823 manuscript. 

Lobachevsky's 1823 manuscript was blind-reviewed by Nicolaus Fuss (1755-1826), internationally renown mathematician and member of Russian \emph{Academy}. Fuss was Leonhard Euler's secretary and his family member; he authored two geometry textbooks which appeared in Russian during the same year \cite{Fuss:1823}, \cite{Fuss:1823a}.\footnote{Details of how Lobachevsky's  \emph{Geometry} was sent to Fuss for review, and the full text of the review, are found in \cite[vol.2, p.126-127]{Lobachevsky:1946-51}.}.

Fuss' review was negative. Acknowledging that Lobachevsky's manuscript contained some original ideas, Fuss stressed the fact that it lacked a logical order and could not be possibly used as a textbook. At least one Fuss' critical remark had a clear political connotation. It concerned Lobachevsky's choice (who followed Legendre and some other contemporary French authors) of \emph{centesimal} grades for measuring angles ($2\pi$ radians = 400 centesimal grades) instead of the traditional grades ($2\pi$ radians = 360 grades). Fuss disapproved on this Lobachevsky's choice claiming (wrongly) that the centesimal grades were invented during French Revolution, which he described in his review as a ``madness of the nation''. Although Lobachevsky was given a chance to resubmit a revised version of his manuscript he never used this opportunity.\footnote{See \ref{centesimal} below. For the history of \emph{centesimal} grades see \cite[ch.10]{Treese:2018}} 

After Lobachevsky's death in 1856 the manuscript of the \emph{Geometry} was lost. But in 1898 it was found in the university archives by Lobachevsky's first biographer Alexander Vasilievitch Vassiliev (1853-1929), and in 1909 published  as \cite{Lobachevsky:1909}. The standard modern edition of this Lobachevsky's work is  \cite[vol.2, p.48-123']{Lobachevsky:1946-51}. 

Lobachevsky's \emph{Geometry} is a tentative textbook in Elementary geometry, which can be described as \emph{radically} d'Alembertian. In the short \emph{Introduction} the author describes geometry as a ``part of pure mathematics that  prescribes ways of measuring space''  and describes spatial \emph{extension} as a ``property of bodies, which consists of the fact that if they are spread then they [eventually] touch each other''. Following d'Alembert \cite[vol. 7 (1757), p. 629A-639B]{Diderot&DAlembert:1751-65} Lobachevsky divides geometry into the \emph{Longimetry} (which treats measurements of distances), the \emph{Planimetry} (which treats measurements of surfaces), and the \emph{Stereometry} (which treats measurements of volumes).  Beware that the latter two terms are used by Lobachevsky in their old traditional sense but not in their today's sense, which is broader \cite[vol.2, p.43]{Lobachevsky:1946-51}, \cite[vol.4, p. 33-34]{Tropfke:1921}.

In \emph{Geometry} Lobachevsky makes no attempt to specify axioms (or first principles of some other sort) or formally separate definitions, statements of theorems, and their proofs. Instead Lobachevsky presents his thoughts in a form of narrative where all these elements are mixed. These features are characteristic for Lobachevsky's writing style which did not significantly change during his career. 

\emph{Geometry} has no content that could be identified as pertaining to Non-Euclidean geometry. Nevertheless this text can be hardly properly understood unless one keeps it in mind that Lobachevsky doubts the validity of Euclid's Axiom of Parallels and wants to develop as much of geometry as he can without using this Axiom (or any of its equivalents). Lobachevsky postpones using this Axiom (without calling it by this name, of course) until a very late stage of theoretical development of his theory, and remarks, referring to the Axiom of Parallels, that ``[n]o rigorous proof of this truth is known to the date; those alleged proofs which have been given [so far] qualify as explanations but not as mathematical proofs.'' \cite[vol.2, p.70]{Lobachevsky:1946-51} 
\footnote{This is not an unprecedented move since already Euclid similarly avoids using Postulate Five in the first 27 propositions of his \emph{Elements}. But Lobachevsky in his \emph{Geometry} goes much further than Euclid in this direction, which makes his order of presentation of geometrical contents very unusual.}. 

In spite of Fuss's prejudgement against the  centesimal grades his negative report on Lobachevsky's  \emph{Geometry} of 1823 appears today sound. This manuscript clearly featured an original approach in the foundations of geometry but it certainly could not be used as a geometry textbook. As V.F. Kagan rightly remarks in his commentary, ``this is an overview of geometry written for those who already know the Elementary geometry but not a textbook''  \cite[vol.2, p.40]{Lobachevsky:1946-51}.

\subsection{ \emph{Brief exposition of the principles of Geometry with the rigorous proof of the theorem of parallels} of 1826}\label{succinte}
\label{brief1826}
 in February 1826  Lobachevsky made another attempt to publish his studies in the foundations of geometry and submitted to the Division of Physics and Mathematics of Kazan University a manuscript titled  \emph{Brief exposition of the principles of Geometry with the rigorous proof of the theorem of parallels}.  The manuscript was written in French, and its original is now lost. \footnote{The original French title of the manuscript was \emph{Exposition succincte des principes de la G\'eom\'etrie avec une d\'emonstration rigoureuse du th\'eor\`eme des parall\`eles}}. The manuscript was accompanied by a note \textemdash\ later published in \cite[p.222]{Modzalevsky:1948} and  \cite[vol.1, p.412]{Lobachevsky:1946-51} \textemdash\ in which Lobachevsky asked his colleagues to evaluate his manuscript and, in the case of a positive verdict, publish it in the University's press. The note also includes the names of three faculties (Simonov, Kupfer and Brashman) who were entitled to give formal reports on this Lobachevsky's manuscript. But these reports never materialised. Most likely the reviewers could not understand Lobachevsky's manuscript and decided not to judge it rather than provide negative reports on the work of their well-connected colleague who very soon after the event (in July 1827) was appointed as the Rector of Kazan University. 
 
 The title of the manuscript suggests that it could contain another attempted proof of Euclid's Fifth Postulate. But in a footnote to his \emph{On foundations of Geometry} of 1829 (see \ref{sfoundations} below) Lobachevsky claimed that the contents of the initial fragment of this latter work was borrowed from his unpublished manuscript of 1826, which the author presented at the meeting of the Division of Physics and Mathematics of Kazan University on February 12, 1826, see \cite[vol.1, p.185]{Lobachevsky:1946-51}. In another footnote Lobachevsky marks the end of this borrowed fragment, see \cite[vol.1, p.207]{Lobachevsky:1946-51}. The borrowed fragment, as we can read it in Lobachevsky's \emph{Foundations}, includes among other things a clear statement concerning the possibility of Non-Euclidean geometry: 
    
``The sum of [internal] angles of rectilinear triangle cannot be $> \pi$. So one can assume this sum to be either $= \pi$ or $< \pi$''. Neither of the two assumptions entails a contradiction in what follows. This is why two [different] Geometries emerge: one \emph{usual} in its simplicity, which accords with [common] practical intentions, and the other \emph{imaginary} which is more general'' [$\dots$ \cite[vol.1, p.194]{Lobachevsky:1946-51}]. 
 
 Thus, unless Lobachevsky deliberately confused his readers, the \emph{Brief Exposition} of 1826 already comprised the key ideas of Non-Euclidean geometry. Since the aforementioned note \cite[p.222]{Modzalevsky:1948}, which constitutes a piece of hard material evidence,  gives a slightly different date of Lobachevsky's presentation of the \emph{Brief Exposition} at his department, namely February 11, 1826, it is commonly assumed that in 1829 Lobachevsky confused the date writing by memory \cite[vol.1, p.411]{Lobachevsky:1946-51}.  
 
 Be it as it may, February 11, 1826 was fixed as the official birthdate of Non-Euclidean geometry by V.F. Kagan and A.P. Kotelnikov and other Soviet mathematicians and historians of mathematics \cite[p.183]{Kagan:1948},  \cite[vol.1, p.412]{Lobachevsky:1946-51} who tried to push this date as far back as possible in order to secure Lobachevsky's priority. During Stalin's rule in 1926 and later in 1952 this date was used for organising pompous celebrations in Kazan and elsewhere which included panegyrical speeches by Kagan and many others  \cite{Parfentiev:1927}, \cite{Norden:1952}. More recently some historians casted doubts that on February 11, 1826 the  \emph{Brief Exposition} was indeed publicly presented \cite[p. 60]{Polotovsky:2015}. Indeed the note \cite[p.222]{Modzalevsky:1948} does not fully support such a conclusion and allows for the possibility that the meeting in the Kazan University merely dealt with the question of how to review the submitted manuscript. But Lobachevsky in his 1929 footnote \cite[vol.1, p.185]{Lobachevsky:1946-51} says explicitly that he \emph{read} the manuscript at this meeting, and we cannot find a sufficient reason to question his credibility in this case.

\subsection{ \emph{On foundations of Geometry} of 1829-1830} (\cite[vol.1, p.187-405]{Lobachevsky:1946-51})\label{sfoundations}

In 1829-1830 Lobachevsky's  \emph{Foundations} have been published in the \emph{Kazanski Vestnik} (\selectlanguage{Russian} Казанскiй В{\cyryat}стникъ \selectlanguage{English}) in five instalments. \footnote{\cite{Lobachevsky:1829a}, \cite{Lobachevsky:1829b}, \cite{Lobachevsky:1829c}, \cite{Lobachevsky:1830a}, \cite{Lobachevsky:1830b}.  \emph{Kazanski Vestnik} was an all-purpose journal published in Kazan University in years 1821-1933. In 1934, by Lobachevsky's initiative as the Rector of Kazan university,  the journal was reformed and changed its name to \emph{Utchenye Zapiski } (\selectlanguage{Russian} Ученые Записки \selectlanguage{English}, Scholarly Notes). Under this latter name the journal continues to be published today.} It is save to assume, in our opinion, that Lobachevsky's appointment to the position of Rector in 1827 greatly helped this publication to happen. 

On August 19, 1832 Lobachevsky arranged for a formal request made by the Scientific Council of Kazan University to the Imperial Academy of Sciences, to review his recently published \emph{Foundations} accompanied by a copy if this text \cite[p. 320]{Modzalevsky:1948}. The request was formally responded on November 7 of the same year by academician Mikhail Ostrogradsky (1801-1862).  In his report Ostrogradsky says that he does not understand the most of Lobachevsky's work, and that the small part that he does understand is erroneous. Ostrogradsky concludes that Lobachevsky's work does not deserve an attention of the \emph{Academy} \cite[p. 332-333]{Modzalevsky:1948}. 

In addition, in 1834  Lobachevsky's \emph{Foundations} were anonymously reviewed in conservative journal for general public \emph{Son of the Homeland} (\selectlanguage{Russian} Сынъ отечества  \selectlanguage{English}) published in Saint-Petersburg. The reviewer suggested that a ``Caricature of Geometry'' would be a better title for Lobachevsky's book and mocked on it in many other ways. Interestingly, the reviewer showed their acquaintance with the aforementioned Ostrogradsky's report even if it was public at that point. There are many suggestions about the identity of the author of this pasquinade in the literature but it still remains unknown.\cite{CC:1834}, \cite[245-253]{Kagan:1948}. According to one popular conjecture, the author of the pasquinade is Fedor Ivanovich Busse  (\selectlanguage{Russian} Федор Иванович Буссе \selectlanguage{English}) (1798-1859), a professor of mathematics education from Saint-Petersburg \cite[p. 65-66]{Polotovsky:2015}.  

Lobachevsky's \emph{Foundations} like his earlier \emph{Geometry} is written in the textbook format. The author once again does not separate here axioms, statements of theorems and their proofs (if any). The \emph{Foundations} is a long continuous mathematical narrative (without author's chapters or sections) that combines in various proportions a text in natural language (Russian), formulas and diagrams. Its content can be partitioned and analysed in a number of different ways of which we propose the following:

\begin{enumerate}
\item The \emph{first part} contains the core of both Euclidean (Flat) and Hyperbolic geometry \cite[vol.1, p.185-206]{Lobachevsky:1946-51}. This part, as the author tells us, is borrowed from the  \emph{Brief exposition} of 1826 \ref{succinte}. The key mathematical result presented in this part can be anachronistically described as the establishment of basic identities of Hyperbolic Trigonometry (in terms of the usual trigonometric functions), which allows Lobachevsky to solve triangles in the Hyperbolic case just like this is done in the Euclidean case using the usual trigonometric formulas. An improved version of the same core theory Lobachevsky later presented in his \emph{Geometrical Researches} of 1840, see \ref{researches}. This initial part of the \emph{Foundations} is very difficult to read without an extensive technical commentary that fills computational leaps and provides necessary proofs (such as the commentary provided in \cite{Lobachevsky:1946-51}). It is a little surprise that neither Ostrogradsky nor any other mathematician whom Lobachevsky could reach at this point of his career was capable to decrypt this text. It goes without saying that Lobachevsky's \emph{Foundations} could not possibly be used for teaching geometry in the university or elsewhere.

\item The \emph{second part} of \emph{Foundations} treats the questions of whether or not the physical space is Euclidean  \cite[vol.1, p.207-210, footnote 32]{Lobachevsky:1946-51}. Lobachevsky tackles this question using astronomical data borrowed from a contemporary astronomical publication \cite{Assas-Montardierd:1828}; he designs a sophisticated computational method allowing one to evaluate the sum of angles of a triangle in the physical space having as its vertices two opposite positions of the Earth on its solar orbit and a distant star.\footnote{Given the absence of non-congruent similar triangles in the Hyperbolic space, such a method is not a matter of course. In fact, the method used by Lobachevsky in the \emph{Foundations} allows only for computing the upper bound of \emph{defect} (i.e., of the difference between $\pi$ radians and the sum of angles of a given astronomic-scale triangle) but does not allow one to tell between Euclidean and Hyperbolic metric of the ambient space. See \cite{Idelson:1949}, \cite{Norden:1958} for details. \label{astro1}} Later historical studies revealed that these astronomical data were not reliable; in addition, Lobachevsky made a trivial computational error of the two orders of magnitude  \cite[vol.1, p.286, footnote 32]{Lobachevsky:1946-51}.  In spite of the erroneous premises Lobachevsky rightly concluded that the available observational data did not allow him to tell the observable physical space from the Euclidean space. In this context Lobachevsky conceives of the possibility that ``only Euclidean assumptions are true [in the physical space] albeit they will remain unproven forever'', and remarks that ``[I]f the new [Hyperbolic] Geometry [$\dots$] does not exist in the Nature it may nevertheless exist in our imagination.'' He further suggests that the new geometry ``makes possible new mutual applications of Geometry and Analysis even without being used in the real [physical] measurements''  \cite[vol.1, p.209-210]{Lobachevsky:1946-51}.\label{imaginary}  

\item In the  \emph{third part} of \emph{Foundations} Lobachevsky realises the Hyperbolic version of his d'Alembertian program: he develops, in sequence, analytic methods for computing lengths \cite[vol.1, p.210-219]{Lobachevsky:1946-51}, surfaces \cite[vol.1, p.219-232]{Lobachevsky:1946-51} and, finally, volumes \cite[vol.1, p.232-249]{Lobachevsky:1946-51} in a Hyperbolic space. In this part Lobachevsky demonstrates his virtuous analytic skills empowered by his novel Non-Euclidean geometrical intuition. Once again, Lobachevsky's computations in this part of his work are quite difficult to follow; they require a genuine decrypting like one provided by Alexandre P. Kotelnikov (1865-1944) in his useful commentary \cite[vol.1]{Lobachevsky:1946-51}. In view of John Milnor, this is the most interesting direction of Lobachevsky's research \cite{Milnor:1982}. Referring to Milnor's reception of Lobachevsky's work, Athanase Papadopoulos remarks that ``the real substantial progress on the subject itself [i.e., on the subject of computing areas and volumes and the Hyperbolic space and its generalisations] is due to [William] Thurston, and it came about 150 years after Lobachevsky's first works appeared in print.''\cite[p. 298]{Lobachevsky:2010}.\footnote{Relevant Thurston's results are presented in \cite[Ch. 7]{Thurston:2002} (or  \cite[Ch. 7]{Thurston:2022}), which is authored by Milnor.}  

\item The  \emph{fourth part} is a continuation of the third but here Lobachevsky changes a perspective: instead of using the Integral Calculus as a tool that helps to calculate lengths, surfaces and volumes in the Hyperbolic space, he now uses the \emph{Imaginary Geometry} as a tool that helps to compute difficult integrals \cite[vol.1, p.249-260]{Lobachevsky:1946-51}. Lobachevsky's techniques in such cases amounts to computing the same geometrical magnitude, say, a volume of a prism, in the Hyperbolic space in two or more different ways, which provides helpful equalities between integrals. In 1836 Lobachevsky published these and some other similar results in \cite{Lobachevsky:1836}, \cite[vol.3, p.175-407]{Lobachevsky:1946-51} as a supplement to his \emph{Imaginary Geometry}. 

\item The  \emph{fifth part} is Lobachevsky's \emph{Conclusion} (this time separated from the rest of this work by the author himself) where he provides a formal meta-theoretical argument showing the consistency of his new geometrical theory (relatively to the Euclidean geometry) \cite[vol.1, p.260-261]{Lobachevsky:1946-51}. Lobachevsky's argument is based on the observation (independently made by Johann Lambert back in 1766 \cite{Papadopulos&Theret:2014}) that a formal replacement of radius $R$ of a given sphere (a real number) by complex (imaginary) magnitude $\sqrt{-1}R$ transforms trigonometric formulae of Spherical geometry on that sphere into the corresponding formulae of Hyperbolic geometry.\footnote{It remains a controversial issue whether or not this Lobachevsky's argument constitutes a proof of logical consistency in the modern sense, see \cite[194-196]{Braver:2011} for a recent discussion}. Characteristically, after presenting this argument Lobachevsky makes the following remark: ``It remains to study the effect of application of the \emph{imaginary} Geometry in Mechanics. Do commonly accepted and certain [mechanical] notions rule out the mutual dependence of [lengths of] lines and angles [which is characteristic for Hyperbolic geometry]?'' \cite[vol.1, p.261]{Lobachevsky:1946-51}. Lobachevsky suggests that there are no such mechanical notions. These Lobachevsky's  remarks suggest that logical consistency was not, in his view, a sufficient condition of validity of his new geometry; the ``physical soundness'' of this mathematical theory was in his view equally and perhaps even more important.

\end{enumerate} 

\subsection{\emph{New foundations of geometry, with a complete theory of parallels} of 1835-1838} (\cite[vol.2, p. 147-588]{Lobachevsky:1946-51})\label{new}

Remarkably, after the poor reception of his \emph{Foundations} Lobachevsky did not give up his attempts to communicate to the mathematical community he discoveries. In 1835 he began to publish in journal \emph{Utchenye Zapiski} (\selectlanguage{Russian} Ученые Записки \selectlanguage{English}, Scientific Memoirs) a new version of the same fundamental work titled this time \emph{New foundations of geometry, with a complete theory of parallels} (\selectlanguage{Russian} Новые основания геометрии с полной теорией параллельных  \selectlanguage{English}). The \emph{New Foundations} was published between 1835 and 1838 in six instalments\footnote{\cite{Lobachevsky:1835a}, \cite{Lobachevsky:1836a}, \cite{Lobachevsky:1836b}, \cite{Lobachevsky:1837a}, \cite{Lobachevsky:1838a}, \cite{Lobachevsky:1838b}}. The \emph{Utchenye Zapiski} started in 1834 as a reformed version of the \emph{Kazanski Vestnik} published earlier in the same university. The reform was conceived and realised by Lobachevsky himself in his capacity of Rector of Kazan University; it brought the journal closer to the European academic standards existing at this time. 

Like the \emph{Foundations} the \emph{New Foundations} comprise no axioms, no formally specified definitions and no formally specified theorems. Nevertheless in the new work Lobachevsky significantly improved the structure of his book, so it became more readable. The \emph{New Foundations} is divided \underline{by its author} into the \emph{Introduction} and thirteen chapters, which are further divided into paragraphs; the paragraphs are numbered end-to-end continuously (there is 203 paragraphs in total). Some of these paragraphs can be qualified as definitions, theorems (in some cases with proofs, and in some other cases without proofs), or commentaries. Thus even if the \emph{New Foundations} has no explicit axiomatic structure its form of presentation comes closer to what one could expect to see in a contemporary Elementary geometry textbook\footnote{like \emph{Kurganov:1765} where a similar end-to-end continuous numbering of blocks is present.}. At the price of a significant increase of the volume (by the factor of four with respect to the   \emph{Foundations}) the \emph{New Foundations} closed (or rather narrowed) many logical gaps found in the \emph{Foundations}. The last two chapters of the \emph{New Foundations} (Chapters 12 and 13) contain a material not found in the \emph{Foundations}: in these additional chapters Lobachevsky elaborates on numerical methods of solving plane (Euclidean) and spherical triangles providing concrete numerical examples. As suggests  A.P. Kotelnikov, these chapters were supposed to provide a theoretical background for a new attempt to evaluate the geometry of physical space on the basis of astronomical data. The  \emph{Introduction} to the \emph{New Foundations} contains important philosophical and methodological remarks as well as a detailed review of earlier attempts to prove Euclid's Postulate Five. New English translation of a relevant fragment of this  \emph{Introduction} is found below in \ref{apH}.

\subsection{\emph{Imaginary Geometry} of 1835/1837 ( \cite[vol.3, p. 16-138]{Lobachevsky:1946-51},  \cite{Lobachevsky:1837}) and the \emph{Application of imaginary geometry to certain integrals} of 1836} ( \cite[vol.3, p. 175-407]{Lobachevsky:1946-51}) \label{simaginary}
Simultaneously with his work on the  \emph{New Foundations} Lobachevsky made an attempt to reach out the international mathematical community by writing research papers in foreign languages and submitting them to mathematical journals published abroad. Early in 1835 (or perhaps even earlier) Lobachevsky sent a manuscript of his \emph{Imaginary geometry} written in French (\emph{ La g\'eom\'etrie imaginaire}) to \emph{Crelle Journal} (Journal f\"ur die reine und angewandte Mathematik); the article appeared in this journal in 1837 as \cite{Lobachevsky:1837}. Soon after the submission Lobachevsky published an extended and partly corrected Russian version of the same paper in the newly founded  \emph{Scientific Memoirs} \cite{Lobachevsky:1835}.\footnote{The fact that Lobachevsky first prepared the French version of this paper, and later himself translated it into Russian (correcting an eventual error and providing some extensions) is evidenced by an author's footnote in the Russian edition, see \cite[vol.3, p.171]{Lobachevsky:1946-51}}. 

In the \emph{Imaginary geometry} Lobachevsky opts for a mode of presentation of the Hyperbolic geometry, which significantly differs from that used in his \emph{Foundations}, the \emph{New Foundations} and the \emph{Pangeometry}. Without trying to deduce the basic identities of Hyperbolic Trigonometry he simply postulates here these identities inviting the reader to \emph{imagine} triangles that might be governed by these formulas in a way similar to which the usual Euclidean triangles are governed by the usual trigonometric formulas. In this analytic presentation of Hyperbolic geometry its intuitive geometrical background no longer plays a foundational role, and can be even wholly ignored. In the beginning of the \emph{Imaginary geometry} Lobachevsky reiterates and reinforces his ``metamathematical'' argument that hinges on the possibility to formally convert  the Spherical trigonometry into the Hyperbolic trigonometry by replacing radius $R$ of a given sphere with $\sqrt{-1}R$ \cite[vol.3, p.22]{Lobachevsky:1946-51}. 

One can think of three independent reasons behind Lobachevsky's choice of the analytic form of presentation in this case. First, as indicates Ostrogradsky's report (and likely other similar reactions that Lobachevsky might receive from his peers), the analytic part of his proposal appeared less controversial than the first synthetic part. Second, Lobachevsky himself believed that the Analytic part of his achievement was more important that the intuitive Synthetic part \ref{sum5}. Third, Lobachevsky's failure to prove the reality of Hyperbolic geometry using astronomical data forced him to conceive of the possibility that his new geometry will always remain a purely theoretical exercise.\footnote{See \ref{sfoundations} above. Compare: ``The assumption of the usual [Euclidean] Geometry [i.e. Euclid's Postulate Five or an equivalent statement] should be taken \emph{as if} it were proven. At the same time one should be convinced that this assumption does not follow from our conception of bodies and thus cannot be proved without the appeal to experience.'' \cite[vol.3, p.26-27]{Lobachevsky:1946-51}.} In that case, in Lobachevsky's view, as we can reconstruct it, the theory had to be given the form of \emph{calculus} rather than geometry proper. 

In 1836 in the \emph{Scientific Memoirs} appeared the \emph{Application of imaginary geometry to certain integrals} ( in Russian) \cite{Lobachevsky:1836},  \cite[vol.3, p.175-407]{Lobachevsky:1946-51} which is a supplement to the  \emph{Imaginary geometry}.  Here Lobachevsky uses his novel geometry as a tool for computing a significant number of earlier unknown integrals. This new Lobachevsky's paper allowed a sceptical reader to wholly dismiss Lobachevsky's claims concerning the alleged new geometry and still appreciate his solutions of open problems in the Integral Calculus. 

\subsection{\emph{Geometrical researches on the theory of parallels} of 1840} (\cite[vol.1, p. 79-159]{Lobachevsky:1946-51}) \label{researches}
The \emph{Researches} is another Lobachevsky's attempt to gain an international recognition by publishing a concise version of his novel theory outside Russia in a foreign language. It is a summary of the first several chapters of the \emph{New Foundations} that fills many of its remaining logical gaps. The \emph{Researches} were first published in 1840 in German as a separate booklet \cite{Lobachevsky:1840}. Jules Ho\"uel claims (without providing any evidence) that the manuscript of the \emph{Geometrical researches} was earlier submitted to the \emph{Crelle Journal} under the title \emph{Beitr\"age zu der Theorie der Parallellinien} (Contributions to the theory of the parallel lines) but was rejected \cite{Houel:1870}. This claim, to the best of our knowledge, is not confirmed by an independent source.

The same year 1840 in Berlin appeared an anonymous negative review to Lobachevsky's \emph{Researches} \cite[p. 412-413]{Modzalevsky:1948}.  Nevertheless it became very soon clear to Lobachevsky that the publication of his \emph{Researches} in Berlin was, in fact, a great success. As evidences Friedrich Gauss' letter to Encke dated January 20, 1841, by this date Gauss had already read the \emph{Researches} as well some other unnamed Lobachevsky's Russian papers, which Gauss managed to obtain from Kazan. \cite[p.415-416]{Modzalevsky:1948}. Gauss claimed that he himself first conceived of the theoretical possibility of the Non-Euclidean geometry back in the 1790s; throughout his carrier he developed this idea in his private correspondence and in private notes but avoided to make these ideas public \cite{Houzel:1992}. He was clearly the best possible judge among those very few Lobachevsky's contemporaries who were capable to understand and appreciate his work. Without giving up his self-imposed ban to talk about the Non-Euclidean geometry in public, Gauss initiated the nomination of Lobachevsky as a member (by correspondence) of G\"ottingen \emph{Royal Society of Sciences} (K\"onigliche Gesellschaft der Wissenschaften). The nomination was formally accomplished on November 11(23), 1842.  A diploma, a letter signed by Gauss in his capacity of the Director of this Society, and by its secretary, were followed by Lobachevsky's formal letter of acceptance addressed to Gauss. Copies of these letter were later found by Modzalevsky in the archives of Kazan university and published  \cite[p.453-454, 451]{Modzalevsky:1948}. Since in Russia the academic recognition of Lobachevsky's mathematical achievements was effectively blocked by Ostrogradsky's authority\footnote{Remarkably, the lack of recognition on the part of the Russian mathematical community did not make any harm to Lobachevsky's very successful administrative career in his university.} the public support of Friedrich Gauss and the  \emph{Royal Society of Sciences} was morally very important for Lobachevsky even if it did not have any immediate social effect.\footnote{Some historians and mathematicians sometimes speculate about the possibility that Gauss and Lobachevsky could exchange more letters. But what we have so far is only these two formal letters without any mathematical content.}

In his letter to Gerling dated February 8, 1844, Gauss notices that in the \emph{Researches} Lobachevsky achieves a conciseness and precision, which makes a sharp contrast with his earlier works \cite[483-484]{Modzalevsky:1948}. So the \emph{Researches} convinced Gauss to take Lobachevsky's works very seriously. Even if Gauss never pronounced about the Non-Euclidean geometry publicly, his formal recognition of Lobachevsky's contribution was the first decisive step toward the wider acceptance of Lobachevsky's ideas in Europe, in Russia, in the U.S., and elsewhere. 

For today's reader the  \emph{Researches} remains the most accessible and useful of all Lobachevsky's writings and thus may serve as a good entry point to studying his works.\footnote{For a modern mathematical commentary to Lobachevsky's \emph{Researches} which provides a step-by-step assisted reading of this mathematical work we highly recommend \cite{Braver:2011}.} 

\subsection{\emph{Pangeometry} of 1855 and 1856}  (\cite[vol.3, p. 435-533]{Lobachevsky:1946-51}) \label{pangeom}

\emph{Pangeometry} has been written by Lobachevsky in Russian in 1855 at the occasion of the 50th anniversary of Kazan university. The anniversary celebrations were, however, postponed until the next year. This explains why the first Russian version of \emph{Pangeometry} was published in a regular way in the \emph{Utchenye Zapiski} \cite{Lobachevsky:1855}.The next year a slightly revised French version of this work appeared  in the collection of papers dedicated to the 50th anniversary of Kazan university\cite{Lobachevsky:1856}. Preparing Russian and later French versions of \emph{Pangeometry} Lobachevsky was already ill and nearly blind. He was assisted by his pupils N.I. Burno and  I.A. Bolzani \cite[vol.3, p.431]{Lobachevsky:1946-51}. In February 1856 Nikolai Ivanovitch Lobachevsky died in the age of 63, and thus the French version of the \emph{Pangeometry} turned out to be the last Lobachevsky's publication during his lifetime. 

The fact that that  \emph{Pangeometry} was available in French (in the author's translation from Russian) greatly helped the dissemination of Lobachevsky's ideas in Europe. The first German translation of  \emph{Pangeometry} appeared as early as 1858 \cite{Lobachevsky:1858}; soon followed translations to some other European languages and new editions of the French text  \cite[vol. 3, p. 534-535]{Lobachevsky:1946-51}. In 2010 appeared a new English translation of \emph{Pangeometry} (from French) with very valuable introduction and commentaries by  Athanase Papadopoulos \cite{Lobachevsky:2010}. 

As suggests Vasiliev, the \emph{Pangeometry} can be read as the second chapter of the  \emph{Researches} \cite[p. 145]{Vasiliev:1992}. Arguably the two works comprise the best concise presentation of Lobabachevsky's ideas achieved by the author himself. It should be borne in mind, however, that the  \emph{Researches} and the \emph{Pangeometry} jointly present a distilled version of Lobachevsky's life-long work leaving behind his philosophical motivations and many mathematical ideas beyond the Hyperbolic geometry.

\subsection{Other mathematical works} 
Lobachevsky's interests and his mathematical publications did not limit to geometry.  When Lobachevsky submitted his \emph{Geometry} to publication in 1823 he was already working on a textbook in Algebra; in 1934 Lobachevsky's \emph{Algebra or Calculus of Finites} was published as a separate book, see \cite[vol.4, p. 23-365]{Lobachevsky:1946-51}). Lobachevsky's conception of algebra is interesting because it includes original finitary versions of Differential and Integral Calculi, see Chapter XV of the \emph{Algebra}. The \emph{Introduction} to Lobachevsky's  \emph{Algebra} is found below as Appendix G \ref{apG}

Lobachevsky's published mathematical output also includes a number of other papers in algebra (treating a special case of binomial equation), analysis (on the convergence of trigonometrical series and integration), statistics, hydrodynamics and astronomy, all of which are found in the two last volumes of his \emph{Collected Works} \cite[vol.4-5]{Lobachevsky:1946-51}. Overviewing these papers is out of the scope of the present work.

 \section{The Sources} \label{sources}

Lobachevsky works and archival materials were published in the USSR in several instalments. During the period 1946 \textendash\ 1951 (the late period of Stalin's rule) there appeared the five volumes of Lobachevsky's \emph{Complete Works} under the general editorship of V.F. Kagan \cite{Lobachevsky:1946-51}, which contained all earlier published Lobachevsky's works provided with valuable commentaries. In 1948 L.B. Modzalevsky published a volume of archival materials related to Lobachevsky \cite{Modzalevsky:1948}. The volume includes some of Lobachevsky's letters, memories of his students and his children and other valuable documents. According to the initial plan, Lobachevsky's  \emph{Complete Works} had to include one more volume containing earlier unpublished manuscripts. After the death of V.F. Kagan in 1953  the work on the sixth volume of the Collected Works stopped. Nevertheless in 1976, a part of relevant archival documents was published in a separate volume under the general editorship of P.S. Alexandrov and B.L. Laptev \cite{Lobachevsky:1976}. In 1979 A.G. Karimullin and B.L. Laptev published lists of titles borrowed by Lobachevsky from Kazan university library and bought by him for this library in his capacity of Rector \cite{Karimullin&Laptev:1979}. Finally, in 1988 another significant portion of archival documents related to Lobachevsky life and work was published by B.V. Fedorenko \cite{Fedorenko:1988}. It complemented Modzalevsky's 1948 collection and reflected the progress in the Lobachevsky scholarship made during the four decades between 1948 and 1988. The present work is made on the basis of the aforementioned publications and does not claim new archival findings.  

For our analysis we selected Lobachevsky's texts where he expressed his general views on mathematics and its role in the sciences. Lobachevsky never published separate books or journal articles on such philosophical topics.Nevertheless he expressed his views on mathematics in writing on many occasions, which provide us with a solid textual evidence sufficient for an accurate historical reconstruction. The following list comprises relevant Lobachevsky's texts put in the chronological order of their writing; English translations of these texts (or more precisely, of their relevant fragments) are found below in the Appendix \ref{ap}. 

\begin{enumerate}
\item A methodological introduction to Lobachevsky's popular 1823 paper on acoustics  \emph{Origin and Diffusion of Sounds in the Air} \cite[p. 396-399]{Lobachevsky:1976}, \ref{apA}. This popular article was first published anonymously in 1923 in \emph{Kazanski Vestnik} (Kazan Messenger).  In the 1950s it was rediscovered by I.N. Bronstein, who attributed it to Lobachevsky on the basis of extant archival evidence, including relevant preparatory materials. The article was republished first in 1956 as \cite{Bronstein:1956}, and then again in 1976 as \cite[p.396-399]{Lobachevsky:1976}. In the Appendix \ref{ap} we include only an introductory part of this article, where Lobachevsky expresses his general thoughts about mathematics and its role in physics and other natural sciences.

\item Teaching plans of the academic year 1824/25 prepared by the request of Kazan University administration \cite[173-185]{Modzalevsky:1948},  \ref{apB};  

\item  Teaching plans of the academic year 1825/26  prepared by the request of Kazan University administration \cite[206-219]{Modzalevsky:1948}, \ref{apC}; 

 \item Public speech \emph{On the most important aims of education} delivered on July 5, 1828 \cite[p. 321-327]{Modzalevsky:1948}; \ref{apD}. This public speech was delivered by Lobachevsky at the occasion of the end of Academic year, which this time marked the end of the first year of Lobachevsky's rectorship.  

\item \emph{Instructions for teachers of mathematics and physics in gymnasia} composed in 1830 in Lobachevsky's capacity as a Rector of Kazan University  \cite[p. 526-531, 534-537]{Lobachevsky:1976},  \ref{apE}; 
 
 \item Introduction to the \emph{On Foundations of Geometry} (1829-30), \cite[vol.1, p. 185-406]{Lobachevsky:1946-51},  \ref{apF};
 
 \item Introduction to \emph{Algebra or Calculus of Finites}  (1834) \cite[vol.4, 23-27]{Lobachevsky:1946-51},  \ref{apG};
 
 \item Introduction to the \emph{New Foundations of Geometry with Complete Theory of Parallels} (1835-38),  \cite[vol.2, p. 148-167]{Lobachevsky:1946-51} and the following Chapter 1 \emph{The First Concepts of Geometry}, \cite[vol.2, p. 168-179]{Lobachevsky:1946-51},  \ref{apH}.
  
\end{enumerate}

These texts and their fragments translated in the Appendix contain a significant amount of repetitions, which we did not always try to avoid.  At the same time, we didn't include into this publications Lobachevsky's fragments drawing on administrative, pedagogical and some other details, which can be interesting for historian but appear to be irrelevant to our present work.

\section{Summary and Discussion}\label{summary}
The present Section contains a summary of  Lobachevsky's views of mathematics, and in particular of geometry, made on the basis Lobachevsky's documents  found in the following Appendix \ref{ap}.  It also contains a preliminary discussion concerning Lobachevsky's motivations, influences and sources. We divide this Summary into ten thematic blocks even if in Lobachevsky's own presentation all of those themes are tightly intertwined. 

We don't claim that the system of views presented in this Section is fully coherent and moreover formally consistent. On the contrary, we identify in it some significant internal tensions which can be easily presented as contradictions. Since these tensions adequately characterise Lobachevsky's way of thinking we don't try to eliminate them in our historical reconstruction.

\subsection{Mathematics and the Natural Sciences}\label{sum1}
Lobachevsky often refers to the ``pure mathematics'' but in his view the mure mathematics hardly enjoys an epistemic autonomy. He consistently sees mathematics as an integral part of the natural science. He argues that the ``unification with mathematics'' marks the progress of any (natural) science. He notices that certain natural sciences are more mathematised than some other sciences, and predicts that in the future chemistry will become mathematical like physics.\ref{apA} . Unless a given science is properly mathematised, argues Lobachevsky, ``we are unable to follow the Nature intellectually'' \ref{apE}. Elsewhere Lobachevsky says that ``physics is everywhere supported by mathematics and owes to mathematics its very existence''\ref{apC}. 

In Lobachevsky's view, mathematics only recently acquired the right method of inquiry, which was indicated by the ``great [Francis] Bacon'' who ``told us that we should stop labouring in vain trying to extract the wisdom from the mind alone'' and ``urged us instead to ask nature''. \ref{apD}.\footnote{On the importance of the historical dimension in Lobachevsky's thoughts about mathematics see \ref{sum8}.} According to Lobachevsky, the modern scientific method pushed forward by Francis Bacon is applicable in the pure mathematics as well as in physics, chemistry and all other natural sciences. This is why including in his \emph{Foundations} of 1829 an astronomical part (bearing on the question whether or not the physical space is Euclidean at the astronomical scales) does not immediately pose to Lobachevsky any methodological problem. 

If the pure mathematics is a part of the natural science then it is appropriate to ask about its specifics \emph{within} the natural sciences. Here is a relevant Lobachevsky's quote:

``Foundations of physics can be hypothetical but foundations of the pure mathematics should be unquestionable truths, [i.e.,] be our first concepts of Nature [$\dots$] which make part of every mental representation and serve as the ultimate foundation of any [human] judgement. Foundations of geometry should meet the same requirements.'' (\ref{apB}). 

So the pure mathematics like any other natural science accounts for the truths of Nature but mathematical truths are somehow more reliable; unlike physical truths, mathematical truths are not merely hypothetical. Since all mature Science, according to Lobachevsky, should become mathematical at certain point of its development, mathematical truths apply across all specific scientific domains. Notice that this view is consistent with the idea of treating Euclid's \emph{Fifth Postulate} (or Playfair's \emph{Axioms of Parallels}) as an uncertain physical hypothesis which should be avoided in a properly mathematical treatment of geometry. 

More light on Lobachevsky's understanding of the place and the role of mathematics (including geometry) in physics and other natural sciences is shed by the following quote:

``Everything in the physical world exists under the necessary condition of being measurable; hence every physical entity is subject to mathematical laws.''\ref{apA}

the above quote shows that Lobachevsky's conception of Natural Science is somewhat instrumentalist. He considers the \emph{measurement} as a joint that connects mathematics to (the rest of) the natural science. But at the same time Lobachevsky does not buy the Kantian idea according to which mathematical truths are \emph{a priori} and determine the scope of possible experience independently of any actual experience.\footnote{As it was already mentioned above in \ref{context} , a plausible source of Lobachevsky's Anti-Kantian stance were published philosophical lectures by T.F. Osipovsky \cite{Osipovski:1807},  \cite{Osipovski:1813}.} Instead, Lobachevsky develops a form \emph{sensualist} account of mathematical knowledge that qualifies as a form of empiricism about mathematics. As Lobachevsky conceived of physical measurements as genuine physical interactions, one may speculate that he could be open to revising at the foundational level the current mathematics in view of newly discovered ways of physical measurements like those discovered in the 20th century Quantum physics.\footnote{ We'll come back to Lobachevsky's sensualism in \ref{sum2} Specific features of \emph{spatial} measurements that lay in the foundation of Lobachevsky's conception of geometry will be discussed in \ref{sum4}. }

In \ref{apB} one finds a trace of Lobachevsky's hesitation between a realist and an instrumentalist interpretation of the mathematical physics. He accepts the reality of physical forces but doubts the reality of other mechanical concepts like time, distance and velocity (as these concepts are construed mathematically in the Classical mechanics). He conceives of the possibility that ``our computations are strictly correct but their [received] foundations are false: our concepts [like extension and velocity] can be artificial but for some unknown reason still capable to replace the true concepts.'' Lobahcevsky notices that these philosophical reflections cannot be useful in teaching science, and he does not return to them in his later writings. Nevertheless these Lobachevsky's thoughts are important for better understanding of his continuing attempts to find ``true'' foundations of geometry, see \ref{sum3} below.  They explain how Lobachevsky applied the Baconian scientific methodology in mathematics. It does not amount to replacing mathematical proofs by physical experiments but rather amounts to reaching an agreement between the foundations of mathematics and the foundations of physics. Unlike Immanuel Kant Lobachevsky did not believe that such an agreement could be possibly forced with an appropriate experimental design\footnote{See for example \cite[B xiii]{Kant:1998}}; he rather believed (referring to Francis Bacon, see \ref{apD}) that the agreement can be reached if mathematics beginning with its foundations ``follows the Nature'', i.e., designs its foundations in a close view of its applications in physics and other natural sciences. This view on mathematics shares with Kant's view the idea that the pure mathematics is a part and parcel of the natural science. But it nevertheless distinctively different, and it motivates a different directions of mathematical research. Lobachevsky's project of redesigning the foundations of geometry in order to make this mathematical discipline better serve the needs of current physics certainly goes beyond the Kantian way of thinking about mathematics.    

In the \emph{Introduction} to his \emph{New Foundations} of 1835 Lobachevsky provides an interesting speculative argument, which brings Lobachevsky's  
reflections on the physical meaning of Non-Euclidean geometry beyond the question of whether or not the physical space is Euclidean, which he thoroughly discusses his \emph{Foundations} of 1829-1830 \ref{sfoundations}), see \ref{apH}. Lobachevsky stipulates here that the core physical concept is that of \emph{motion} and that motions, generally, are caused (or otherwise regulated) by certain physical \emph{forces}. He further stipulates that `` geometrical concepts are produced in our minds artificially on the basis of properties of motion''. It follows that ``space does not exist for us separately and by itself'' but depends on involved physical forces and ways of measuring distances.\footnote{For further details on Lobachevsky's conception of space see \ref{sum4} below }. Lobachevsky concludes that ``the assumption according to which some forces in Nature follow one geometry while some other forces follow another specific geometry, cannot produce a contradiction in one's mind''. This striking conclusion, which for today's reader echoes Albert Einstein's later achievements, is illustrated by Lobachevsky with the example of Newton’s inverse-square law of gravitation, which he relates to the Euclidean character of the Newtonian space. Lobachevsky speculates that different laws describing physical interactions at microscopic scales might make the Hyperbolic geometry physically adequate at this scale. Remarkably, after the failure of finding an empirical confirmation of the reality of Hyperbolic geometry at the astronomical scales Lobachevsky turns to the possibility to find such a confirmation at microscopic scales.\footnote{Lobachevsky's attempts to provide his new geometry with a physical meaning are further discussed in \cite{Daniels:1975}}.

\subsection{Senses} \label{sum2}

 Lobachevsky justifies his choice of \emph{motion} as the core physical concept by saying that ``in the Nature we properly grasp only motion without which no sensual perception is possible''. This is an epistemological rather than ontological justification. The idea that the five human \emph{senses} are in the origin of all human knowledge including the pure mathematics (and thus including geometry), is in the core of Lobachevsky's epistemology; he makes strong claims to this effect repeatedly and consistently, see \ref{apC}, \ref{apE}, {apF}, \ref{apH}. Moreover, Lobachevsky insists that ``mathematics $\dots$ should never go beyond the requirements of our senses'' because it ``aims at the actual [physical] measurement'' (rather than at some unrealistic imaginary measurement)  \ref{apC}. This letter condition created significant tensions with Lobachevsky's contemporary mathematical practice which already involved certain ``imaginary'' objects like $\sqrt{-1}$ as well as with his own notion of ``Imaginary'' geometry. We'll discuss these tensions separately in \ref{sum6}. 
 
 A probable source of Lobachevsky's sensualism is the \emph{Treatise on the Sensations} by \'Etienne Bonnot de Condillac first published (in the original French) in 1754 \cite{Condillac:1754}. This historical conjecture remains so far unsupported by a direct evidence but it is nevertheless strongly supported by the fact that in his putative foundations of geometry Lobachevsky secures a special role for the sense of \emph{touch}(see  \ref{sum4} below), which also plays a special role in Condillac's \emph{Treatise}.\footnote{I thank Paul Rusnock for this observation}  Further parallels with Condillac are found in Lobachevsky's philosophy of Language, which will be discussed in \ref{sum6}. Lobachevsky could develop his interest in Condillac via Osipovsky who in 1805 published a Russian translation of Condillac's \emph{Logic} \cite{Condillac:1780},\cite{Condillac:1805}.

In \ref{apF} Lobachevsky provides, however, an argument, which demonstrates the failure of naive forms sensualism in sciences. While the claim that the concepts of natural number and geometrical magnitude have been formed on the basis of everyday human sensual experience is plausible (as defended by Candillac in his  \emph{Treatise}), the claim that fundamental \emph{mechanical} concepts and mechanical laws have been discovered similarly is not plausible. By mechanics Lobachevsky means, of course, his contemporary Newtonian mechanics, not the outdated Aristotelian mechanics, which better agrees with common intuitions about the locomotion (but at the same time fails to account for many common phenomena like the motion of a thrown ball). \footnote{On the Aristotelian mechanics see \cite{Leeuwen:2017}}. Lobachevsky rightly stresses the fact that the law of inertial motion which belongs to the foundation of the Newtonian mechanics cannot be directly observed or otherwise straightforwardly experienced; in the real history of science the principles of Classical mechanics required the genius of Galileo and Newton to be conjectured, properly formulated and finally empirically confirmed, see \ref{innate}. Continuing this Lobachevsky's argument one can remark that fundamental concepts of Classical mechanics such as \emph{mass}, \emph{force} and \emph{energy} do not straightforwardly correspond to anything found in the common human experience either (even if they may have an equally strong intuitive appeal for educated and non-educated individuals). These mechanical concepts can be properly introduced only along with the theory of Classical mechanics itself. It is important to realise that Lobachevsky's project of building new foundations of geometry on naturalistic grounds is supposed to follow the modern pattern of Galileo and Newton but does not reduce to a naive attempt to develop geometry on the basis of everyday spatial experience. To repeat, were Lobachevsky alive today, he would most certainly have found research programs such as the development of \emph{Quantum geometry}  \cite{Torma:2023} very appealing.

\subsection{Naturalised foundations of Geometry} \label{sum3}
Let us now see how Lobachevsky applies the above general principles to his project of rebuilding the traditional foundations of geometry. According to Lobachevsky's assessment of 1824 ``[T]he Foundations of Geometry remain obscure. $\dots$ [S]o many things in Geometry do not stand against a rigorous critical analysis''. Lobachevsky's critical analysis is not merely logical but rather epistemological:

``[T]he conceptual obscurity [in the traditional foundations of Geometry] is produced by the abstractness, which is redundant in the actual measurements and  hence also unnecessary in the theory itself. Surfaces, lines and points as they are [usually] defined in Geometry exist only in our imagination while the [actual] measurement of surfaces and lines is done using bodies. This is why it is appropriate to talk about surfaces, lines and points only as so far as these concepts are used in the actual measurement. In this case we shall stick to concepts which are familiar to our imagination, and which can be verified in the Nature straightforwardly without using other artificial and irrelevant concepts.''  \ref{apH}

In order to better understand Lobachevsky's project it is appropriate to recall the context of his work \ref{context}. Along with Osipovsky, Legendre, and a number other mathematicians both in and outside Russia, Lobachevsky is a strong proponent of the ``modern'' approach in teaching geometry described by d'Alembert in the \emph{Encyclopaedia} \cite[vol. 7 (1757), p. 629A-639B]{Diderot&DAlembert:1751-65}. In the section  titled \emph{On the subject-matter of Geometry} d'Alembert writes the following:

``We start with considering bodies with all their sensible properties; then we step-by-step separate and abstract away those various properties in our mind. [In this way] we arrive to the conception of body as a penetrable, divisible and figured portion of the [spatial] extension. Thus a geometrical body is nothing but a portion of the extension, which is delimited in all senses. First of all, as a general assumption, we consider this extension to be three-dimensional. But later, in order to easier determine its properties, we consider first just one dimension, that is, the length, then two dimensions, that is, the surface [i.e., the air], and finally the three dimensions together, i.e. the solid [i.e., the volume]. So the properties of lines, surfaces and solids constitute the subject-matter of Geometry and its natural division [into the Longimetry, the Planimetry, and the Stereometry].

We can consider lines without width and surfaces without depth via a simple mental abstraction. Thus Geometry represents bodies in an abstract state, in which, in reality, they are not. The truths concerning these bodies, which are discovered and proved in Geometry, are therefore purely abstract and hypothetical truths. But these [geometrical] truths can be nevertheless useful.  [$\dots$] Even if mathematical theorems do not hold exactly in the Nature, they are nevertheless practically useful allowing one to measure with a sufficient precision distances between inaccessible places  [$\dots$]  and predict astronomic phenomena. In order to prove truths concerning the form of certain bodies with the full rigour one is obliged to consider these bodies in a perfectly abstract state, in which they are really not.''  \cite[vol. 7 (1757), p. 632B]{Diderot&DAlembert:1751-65}
\footnote{Nous commen\c{c}on par consid\'erer les corps avec toutes leurs propri\'et\'es sensibles ; nous faisons ensuite peu-\`a-peu \& par l'esprit la s\'eparation \& l'abstraction de ces diff\'erentes propri\'et\'es ; \& nous en venons \`a consid\'erer les corps comme des portions d'\'etendue p\'en\'etrables, divisibles, \& figur\'ees. Ainsi le corps g\'eom\'etrique n'est proprement qu'une portion d'\'etendue termin\'ee en tout sens. Nous consid\'erons d'abord \& comme d'une v\^ue g\'en\'erale, cette portion d'\'etendue quant \`a trois dimensions ; mais ensuite, pour en d\'eterminer plus facilement les propri\'et\'es, nous y consid\'erons d'abord une seule dimension, c'est \`a-dire la longuere, puis deux dimensions, c'est \`a-dire la surface, enfin les trois dimensions ensemble, c'est \`a-dire la solidit\'e : ainsi les propri\'et\'es des lignes, celles des surfaces \&celles des solides sont l'objet \& la division naturelle de la G\'eometrie.

C'est par une simple abstraction de l'esprit, qu'on considere les lignes comme sans largeur, \& les surfaces comme sans profondeur : la G\'eometrie envisage donc les corps dans un \'etat d'abstraction o\`u ils ne sont pas r\'eelement ; les v\'erit\'es qu'elle d\'ecouvre \& qu'elle d\'emontre sur les corps, sont donc des v\'erit\'es de pure abstraction, de v\'erit\'es hypoth\'etiques ; mais ces v\'erit\'es n'en son pas moins utiles.  [$\dots$] Mais si les th\'eor\`emes math\'ematiques n'ont pas exactement lieu dans la nature, ces th\'eor\`emes servent du-moins \`a trouver avec une pr\'ecision suffisante pour la pratique, la distance inaccessible d'un lieu \`a un autre, [$\dots$], \`a pr\'edire les ph\'enomenes c\'elestes. } 

Legendre (in \cite{Legendre:1794}) and Osipovsky (in \cite{Osipovski:1801}) follow the above d'Alembert's  directions only partly; both introduce the Cartesian notion of \emph{extension} (French \emph{\'etendue}, Russian  \selectlanguage{Russian} \emph{протяжение} \selectlanguage{English}) and the idea of measurement in the very beginning of their courses but very soon they return to the traditional Euclidean concepts of line and surface; Legendre does not even avoid reproducing the traditional Euclid's definition of line as \emph{length without width} \cite[p. 1]{Legendre:1794}.) Similar remarks can be made about other contemporary geometry textbooks. Lobachevsky, according to our reading, intends to go further in the same direction and implements d'Alembert's ideas in his geometrical theory at a deeper theoretical level. 

In order to achieve this end, Lobachevsky attempts to provide  d'Alembert's discourse concerning the development of basic geometrical concepts via abstraction with a concrete geometrical meaning. This is why Lobachevsky's geometrical theory  in its mature form (after the \emph{Geometry} of 1823) begins with considering three-dimensional \emph{bodies} (solids) and only afterwards comes to surfaces, lines, and points. \footnote{In this way Lobachevsky qualifies as an early representative of a movement in mathematics education, which later in the 19th century was called \emph{fusionism}. Fusionists reject the traditional order of teaching and learning geometry which begins with the Plane geometry and only much later proceeds to the Stereometry. Instead, they ``fuse'' these two parts of the subject from the very beginning of their courses \cite{Borgato:2006}. For Lobachevsky, however, the reason of this fusion is not merely pedagogical.} Lobachevsky's theoryexplains in mathematical terms how abstract geometric concepts like that of line and surface are obtained from the fundamental concepts of \emph{body} and  \emph{touch}.  This theory qualifies a form of \emph{Dimension theory}. We discuss it separately providing more details in \ref{sum4}.

\subsection{Spatial Measurements and the Dimension Theory. Touch and Space} \label{sum4}
The starting point of Lobachevsky's foundations of geometry, which connects this mathematical discipline with his sensualist epistemological views, is the concept of \emph{touch}, which Lobachevsky probably borrows from Candillac \ref{sum2}. We find this concept already in Lobachevsky's \emph{Plan} of 1824  \ref{apB}, written before Lobachevsky developed his ideas about the Non-Euclidean geometry.  According to Lobachevsky, the \emph{touch} is not just one primitive geometrical concept among others but the defining property of the very subject-mater of geometry:

 ``I believe that geometry studies the property of natural bodies, which is called the \emph{touch} [$\dots$].'' \ref{apB}\footnote{According to Condillac, the sense of touch alone allows humans (and other animals) to make a distinction between one's body and other bodies, and on this basis for form the ideas of space and spatial extension. In order to demonstrate his point Condillac imagines a human-like Stature \textemdash\ today we would rather call it a Robot \textemdash\ that is endowed only with the sense of touch:

 ``Since our Statue lacks the sense of smell, hearing, and taste, and has only the sense of touch, it exist  first of all with the feeling that it obtains from the action of some parts of its body onto some other parts such as the movements caused by the respiration. This is the minimal degree of sensation, which we can attribute to the Statue. I will call it the ``fundamental sensation'' because the animal life begins with such a mechanical activity.''
 
( ``Notre statue priv\'ee de l'odorat, de l’ou\"{\i}e, du go\^{u}t, de la vue, et born\'ee au sens du toucher, existe d'abord par le sentiment qu'elle a de l'action des parties de son corps les unes sur les autres, et surtout des mouvements de la respiration : voil\`a le moindre degr\'e de sentiment, o\`u l'on puisse la r\'eduire. Je l'appellerai \emph{sentiment fondamental} ; parce que c'est \`a ce jeu de la machine que commence la vie de l'animal : elle en d\'epend uniquement.'' \cite[Part 2, Chapter 1 (p. 204-205)]{Condillac:1754})

(\emph{Italic} is Condillac's). In the following chapters of his \emph{Treatise on Sensations} Condillac elaborates how the Stature can develop the concepts of space and spatial extension. 
}
 
 The concept of \emph{touch}, according to Lobachevsky, is closely related to another basic geometrical concept, namely that of (spatial) \emph{measurement}:

`` In Geometry, [physical] bodies are measured only in respect of their [mutual] touch, so this property should be the subject-matter of this mathematical discipline.'' \ref{apB}

The notion of spatial \emph{measurement} as the principle subject-matter of geometry is common d'Alembert-style geometry textbooks.\footnote{Cf. in Legendre: ``Geometry is a science that studies measuring of spatial extension.'' (La g\'eom\'etrie est une science qui a pour objet la mesure de l'\'etendue. \cite[p.1]{Legendre:1794})} Lobachevsky's idea that spatial measurement is grounded in \emph{touch} is original, and it doesn't reduce to epistemological consideration but has a certain mathematical content, which we refer here as Lobachevsky's \emph{Dimension theory}. This theory is systematically presented in the \emph{Introduction} to his \emph{Foundations} of 1929 \ref{apF} and  (with some additional details) in \emph{Chapter 1} of the \emph{New Foundations} of 1835 \ref{apH1}. The concept of dimension was earlier discussed by d'Alembert in the \emph{Encyclopaedia}\footnote{See references in \cite[p. 29]{Johnson:2025}} but little advanced were made at this point of history towards its mathematical treatment. Lobachevsky's is an early original attempt to develop a version of Dimension theory in the top-down way, i.e, beginning with the concept of body and then revealing its three-dimensional character. In this work we restrain from giving any mathematical assessment of this theory from the point of view of today's theories of dimension. The reader is invited to make their own assessment of this theory after reading \ref{apF} and \ref{apH1}. In any event Lobachevsky's Dimension theory (unlike his works in the Non-Euclidean geometry) had no posterity and made no impact on the development of the Dimension theory in the late 19th and the 20th century. It is interesting to observe, however, that in the first half of the 19th century there was another mathematician in Europe, namely, Bernard Bolzano (1781-1848), working on Dimension theory in the intellectual isolation; his works pertaining to this theory remained almost unknown util very recently \cite[Chapter 4]{Johnson:2025}. It goes without saying that Bernard Bolzano and Nikolai Lobachevsky worked on the Dimension theory wholly independently and with very different motivations. 

This is how Lobachevsky defines in his \emph{Foundations} \ref{apF} the notions of  \emph{measure} and \emph{distance} using his Dimension theory based on the primitive concept of touch. 

``To \emph{measure} a body means to count the congruent parts into which this body is divided in three dimensions using translational sections along with another body taken for the unit. ''  \ref{apFM} 

\emph{Section} in Lobachevsky's parlance is an operation dual to \emph{touch}: composed body $AB$ is dissected into its parts $A,B$ just in case these parts (bodies) $A,B$ mutually touch each other. Lobachevsky distinguishes between translational and rotational sections as follows: a part $P_{i}$ obtained via a \emph{translational} section of a given body touches at most two adjacent parts $P_{i-1}$  and $P_{i+1}$; a rotational section may divide a given body into any number of parts all of which remain (pairwise) in touch. 

In order to define the concept of \emph{congruence} Lobachevsky first introduces the concept of \emph{ambient body}: body $A$ is ambient body of body $B$ when (i) the two bodies are in touch and (ii) body $B$ cannot touch any third body. In case the composed body $AB$ is the \emph{space}, i.e., the totality of all bodies, body $A$ (always satisfying conditions (i)-(ii)) is also called the \emph{ambient space}. Following Aristotle  Lobachevsky also calls the ``void'' taken in $A$ by body $B$ a \emph{place}. Now he calls two given bodies $X,Y$ congruent when they \emph{can} fill the same place. The character of this geometric modality needs, of course, a further analysis and clarification. 

These definitions reveal Lobachevsky's hesitation between the absolutist conception of space defended by Newton and his followers and the relational conception dating back to Aristotle and later defended by Leibniz \cite{Barbour:1998}. When Lobachevsky describes the composition of a body and its ``ambient'' complement, which is another body, he apparently follows the Aristotelian \emph{no vacuum} principle and leans towards the relational conception of space. But when Lobachevsky introduces on the top of the space concept the concept of \emph{place} describing it as a ``void inside the space'' he leans towards the Newtonian absolutist conception. Aristotle avoids this difficulty by defining place as a ``[internal] boundary of containing [i.e., ambient] body'', see his \emph{Physics} $\Delta$, 212a6-7.  

Saying that the measured body is divided into congruent parts ``along with the unit'' Lobachevsky probably points to the possibility to divide the given unit into smaller units.   

After introducing the general concept of measure Lobachevsky makes a further effort towards distinguishing between different kinds of geometrical measures beginning with the concept of linear measure, that is, the concept of \emph{distance}. Surfaces, lines and points Lobachevsky identifies with different kinds of \emph{touches} between bodies (providing obvious examples but hardly a theoretical justification of this distinction): thus he distinguishes between \emph{surfacial}, \emph{linear} and \emph{point-like} touches. In this way the traditional Euclidean concepts of point, line and surface become specific sorts of touch in Lobachevsky's theory.\footnote{Given Lobachevsky's intellectual ambition one might expect that he could come up with a wholly new set of basic geometrical concepts replacing the standard Euclidean set. This did not happen during Lobachevsky's lifetime. Will the project be eventually accomplished by the future generations of mathematicians?} Then he proceeds as follows:

``If [i] two bodies $A, B$ touch a third body $C$ in points [$P_{A}, P_{B}$ correspondingly], and [ii] $A, B$ are connected through certain body $D$, which does not touch body $C$, then the relative position of [$A$ and $B$] is determined, which is also called the \emph{distance} between bodies $A$ and $B$. The distance so determined is not affected by eventual changes in $A, B, D$ related to taking out from or adding some new parts to these bodies provided these new parts do not touch body $C$. [More generally,] the distance is not affected by changes in $A, B$, which do not affect the touch of these bodies to body $C$. Thus a compass serves for assigning distances [between points].'' \ref{apFD} 

Condition [i] reduces the notion of distance between \emph{bodies} $A, B$ to that between their \emph{points} $P_{A}, P_{B}$ determined by touches of $A, B$ to certain third body $C$. Lobachevsky needs here the third body $C$ in order to introduce points $P_{A}, P_{B}$ (since points are determined, in his approach via a touch of two bodies). According to this Lobachevsky's stipulated definition the distance between certain bodies $A, B$ reduces to the distance between some \emph{points on the boundaries} of these bodies. Notice that in the given setting body $C$ serves as a \emph{unit of length} for measuring distances between bodies like $A, B$.    

The role of condition [ii] is less clear. Why Lobachevsky needs to apply in his theory of distance and measurement one more additional body $D$? Lobachevsky's pointer to the example of compass gives us a hint. Think of $A, B$ as two legs of a compass which are connected by their joint $D$. Thus in addition to aggregate body $ACB$ that helps Lobachevsky to introduce the concept of distance (between points $P_{A}, P_{B}$) we get another aggregate body $ADB$ (the compass) which apparently may play a similar role. But notice now that the concept of measurement involves a relation between \emph{what} is measured (an object) and an  \emph{instrument} used for measuring. The two things can be indeed structurally similar or even identical: think of a stick measured by another stick. Since Lobachevsky aims here at a theory of \emph{measured} distance (rather than some kind of ``absolute'' human-independent distance) he should take into account this relational character of measurement and consider (1) the measuring instrument  $ADB$,  (2) the measured object $C$, and finally (3) the parts $A, B$ of the instrument which interact with the measured object directly via a touch.

In the \emph{New Foundations} Lobachevsky presents yet more elaborated and more technical version of the same theoretical development, which we reproduce in the following Appendix only partly \ref{apH}.

\subsection{Analysis and Synthesis}  \label{sum5}
This distinction is central in Lobachevsky's reflections throughout his continuing attempts to reform the foundations geometry; we find it  in his \emph{Teaching Plans} of 1824/1825 and 1825/1826 \ref{apB},\ref{apC}, in his \emph{Algebra} \ref{apG}, and in the most developed form in the \emph{New Foundations} \ref{apH}. A short remark on Analysis and Synthesis is also found in the \emph{Introduction} to Osipovsky's textbook  \cite[p.2]{Osipovski:1801}, so it is possible that Lobachevsky's interest to this topic was triggered by Osipovsky. But while Osipovsky refer to Analysis and Synthesis only as to very general methods of acquiring knowledge, Lobachevsky discusses the same distinction in the context of different styles of geometrical reasoning, namely, the Euclid-style constructive and intuitive \emph{synthetic} geometrical reasoning and the \emph{analytic} reasoning empowered by symbolic algebraic methods (including infinitesimals methods), which since the beginning of the 18th century constituted the mainstream of mathematical progress in Europe and elsewhere.\footnote{Since Leonhard Euler joined Russian Academy of Science in 1727 and moved to Russia a significant part of these new developments was geographically located in Saint-Petersburg.} Introducing these mathematical contexts Lobachevsky nevertheless does not abandon the traditional notion of Analysis as ``ascending to the first principles'' and the corresponding notion of Synthesis as ``proceeding from the first principles to conclusions'' \ref{apG}. 

Lobachevsky's thinking about Analysis and Synthesis has a historical dimension or, more precisely, the dimension of \emph{historical epistemology}:  

``Synthesis, which is  [historically] an earlier [human] invention, remains applicable only in the foundations [of mathematics]. Later [in the history] Synthesis gave way to Analysis and recognised its superiority.'' \ref{apC}

Saying that Synthesis is an ``earlier invention'' Lobachevsky, if course, refers to the fact that Analytic geometry (in other words, the applications of symbolic algebra in geometry) emerged much later than the Euclid-style geometry, which provided to Lobachevsky the paradigm of Synthesis in mathematics.\footnote{Analytic geometry emerged in 16-17th centuries in mathematical works of Fran\c{c}ois Vi\`ete \cite{Viete:1591} and Ren\'e Descartes \cite{Descartes:1637}.} 

Lobachevsky's claim of the ``superiority'' of Analysis in the above quote is easy to understand in historical terms as a pointer to the obvious fact that the emergence and later development of the Analytic geometry supported the great progress made in mathematics during the two centuries preceding the date of writing (1825). But Lobachevsky's claim appears to be stronger. He doesn't simply say that new analytic methods discovered in the 16th-17th centuries allowed mathematics to achieve a significant progress. He also claims that mathematicians at certain point of history changed their epistemic priorities and ceased to prioritise the Euclide-style Synthesis as the standard way of presenting their achievements. They ceased to think about Analysis only as a useful \emph{way} to obtain mathematical knowledge but began to think about it as a form of mathematical knowledge on its own. 

Let us stress that the idea of epistemic primacy of Analysis is at odds with the popular notion according to which a piece of ready-made mathematical knowledge is supposed to be presented in the form of \emph{axiomatic theory}  \textemdash\  that is, in the \emph{synthetic} form where one begins from first principles (axioms) and draws certain conclusions (theorems) via the logical inference. Lobachevsky associates this axiomatic way of thinking about mathematics with Euclid and judges it outdated. Since in the 20th century the axiomatic way of reasoning in mathematics and elsewhere has been largely rehabilitated by David Hilbert \cite{Hilbert:1918} and became anew popular, our understanding of the earlier critique of this method by Lobachevsky and his contemporaries, requires today a special hermeneutical effort.   

An important breakthrough in our current understanding of Analysis in terms of the 20th century logic occurred back in 1974 when Jaakko Hinitkka and Unto Remes published their monograph \cite{Hintikka&Remes:1974}. More recently Carlo Cellucci proposed a fundamental revision of the common 20th century ideas about Logic and Knowledge, which gave more epistemic significance to the concept of Analysis \cite{Cellucci:2013},\cite{Cellucci:2017}.
More relevant literature can be found in \cite{Beaney&Raysmith:2024}.This literature greatly helps one to avoid a widespread prejudice according to which mathematical knowledge always exists (or \emph{should} always exist for some timeless deontic reason) in the form of axiomatic theories. A general discussion on this topic would lead us too far afield but the reader will find below some reconstructed historical arguments supporting the epistemic primacy of Analysis over the axiomatic Synthesis.   

Indeed, Lobachevsky like many of his contemporaries did not think about mathematics, and in particular about geometry, axiomatically.\footnote{This point has been quite rightly stressed by Jeremy Gray back in 1979 \cite{Gray:1979}} As we have already suggested in \ref{context}, Lobachevsky's decision to avoid the axiomatic organisation of his geometry textbooks could be influenced by d'Alembert \cite[vol. 7 (1757), p. 629A-639B]{Diderot&DAlembert:1751-65}. Another possible influence that led Lobachevsky to this decision was  Condillac's  \emph{Logic}, which Lobachevsky could read both in the original French \cite{Condillac:1780} and in Russian translation by Osipovsky \cite{Condillac:1805}. Condillac provides in this work a long argument according to which all forms of traditional synthetic organisation of mathematical and scientific knowledge must be definitely abandoned and replaced by an analytic form of presentation using artificial symbolic languages. He prises Algebra as an artificial language for Science and rejects the logical apparatus of the traditional synthetic Euclid-style geometry including traditional geometrical definitions, and what he calls the ``method of doctrine'' cite[p.168]{Condillac:1780}. Condillac makes explicit some sources of his inspiration pointing to works by Euler and Lagrange as examples of the new analytic style of presentation of mathematical knowledge \footnote{See cite[p.171-173]{Condillac:1780}. Condillac does not provide exact references but since he mentions ``Euler's Elements'' he probably means his course of Differential and Integral Calculi \cite{Euler:1755}, \cite{Euler:1768-1770}. } Condillac applies his epistemic principles to his own work warning the reader that he is not going to begin this work ``with definitions, axioms and principles''; instead he begins his logical treatise ``with certain lessons given us by the Nature'' \footnote{``Nous ne commencerons donc pas cette logique par des d\'efinitions, des axiomes, des principes : nous commencerons par  observer les le\c{c}ons que la nature nous donne.''\cite[p.8-9]{Condillac:1780}} 

Lobachevsky similarly sees mathematics (including geometry) primarily as a \emph{language} that allows one (i) to formulate fundamental physical concepts (like the concepts of mass, force and velocity in the Classical mechanics) and (ii) solve physical and engineering problems. In this way mathematics helps to exercise ``knowledge as power'' in the Baconian vein. From this point of view the most valuable part of mathematics are symbolic calculi (including the traditional Arithmetic, Elementary Algebra and the Differential/Integral calculi), which allow one to solve problems using symbolic computations rather than traditional Euclid-style constructive reasoning. In this context Lobachevsky emphasises the traditional \emph{heuristic} aspect of Analysis, i.e., its capacity not only to express and justify a ready-made solution but also to provide an assistance in finding this solution at the first place. The heuristic aspect of Analytic geometry (including the Differential geometry) is at least as important for Lobachevsky as the validity of results that this mathematical discipline allows one to obtain.\footnote{In order to appreciate the heuristic power of the standard Analytic methods with which Lobachevsky was familiar the reader is invited to compare Archimedes' computation of the area under parabola \cite{Archimedes:2009} with the modern solution of this problem found in any today's textbook on Calculus of the undergraduate level. The solution of this problem, which during Archimedes lifetime required a lot of his mathematical ingenuity, since the mid-18th century has been routinely obtained by generations of students by a simple symbolic computation. It is, of course, a wholly different question whether this development made any good to the mathematics education.} 

Lobachevsky's thoughts about Analysis and Synthesis do not, however, reduce to praising the power of Analysis. \emph{Contra} Condillac Lobachevsky argues that Synthesis is indispensable in geometry and mechanics. Saying this, Lobachevsky limits the scope of Synthesis to  \emph{foundations} of these disciplines. Thus the traditional roles of Analysis and Synthesis in Lobachevsky's account are reversed: Synthesis serves as a preparatory step for Analysis but not the other way round. Let us now see how it works providing more details.

By ``foundations'' Lobachevsky understands primarily the conceptual foundations, i.e., the set of basic concepts in a given scientific discipline. We are taking  here the liberty of using our own example in order to explain Lobachevsky's point. Consider the case of trigonometric functions which play a major role in the Analytic geometry. Theoretically, trigonometric functions can be introduced and further used in pure mathematics and beyond by purely analytic (symbolic) means without any reference to geometry and to triangles, namely as infinite sums which involve only elementary arithmetical operations (Taylor series). But such treatment of trigonometric functions will leave them without any intuitive support and will not explain where these functions come from. It will also obstruct applications of trigonometric functions in physics.  To avoid mentioning the geometrical background of trigonometric functions in mathematical courses altogether would be, in Lobachevsky's view, among other things a pedagogical error: 

``Synthesis should prepare a given discipline to become a perfect subject to Analysis as this happens in geometry and mechanics. A systematic teaching requires a strict separation of such [synthetic] parts of mathematics [from analytic parts] in order to make explicit their specific features and their sources.'' \ref{apC}

But the role of synthetic foundations on Lobachevsky's account is not \emph{merely} pedagogical:   

``Only Synthesis can provide science with its basic concepts, which serve as foundations of all further judgements.  These further judgements are derived from the first data [i.e., the first principles provided by Synthesis]; thus they extend the horizon of our knowledge in all directions unlimitedly. The first data are, beyond any doubt, the concepts obtained in the Nature via our senses. ''  \ref{apH} 

Thus, in Lobachevsky's view, only the synthetic foundations are capable to provide the crucial link between geometry and the Nature. Without such a link the power of Analysis would be in vain. Recall from \ref{sum3} that in Lobachevsky's view the traditional Euclid-style synthetic geometry does this job very poorly being full of unnecessary abstractions like depth-less surfaces, breadth-less lines, part-less points and incommensurable magnitudes. Hence Lobachevsky's will to replace these traditional foundations of geometry by more ``natural'' ones.   

One might ask Lobachevsky why in his view a link between the Analysis and the Nature needs to be mediated by the intuitive geometrical (or mechanical) Synthesis? Why such an epistemic link could not be established between the Analysis (i.e., an appropriate symbolic analytic calculi) and the Nature directly via the Measurement? It appears that Lobachevsky is not wholly opposed to this line of argument. He quite sympathetically writes about Lagrange's effort to formulate and teach mechanics in purely analytic terms. Only Lagrange's theory of functions formulated analytically without mentioning its geometrical origin (and Lobachevsky could also say, without its geometrical foundation) meets Lobachevsky's objections \footref{lagrange}.  Lobachevsky does not further elaborate on this point but apparently he simply cannot figure it out how, say, the concept of trigonometric function could be found and successfully applied in Science bypassing its traditional geometrical representation. This is why Lobachevsky is not ready to fully embrace Condillac's view according to which a properly designed \emph{language}, and particularly an artificially designed precise symbolic language of Algebra, may by itself constitute a vehicle for scientific knowledge which may not need auxiliary means such as foundations of geometry (or foundations of physics) and synthetic geometric (physical) theories.\footnote{Cf. Condillac's position:

``Since languages which have been formed according to human analysing activities became themselves analytic methods, we feel it natural to think accordingly to our linguistic habits. We think with languages which provide rules for our judgements. Languages produce are knowledge, our opinions, our prejudgements. They produce everything what is good and what is bad in this domain.''

(Puisque les langues, form\'ees \`a mesure que nous analysons, ont devenues autant de m\'ethodes analytiques, on conçoit qu'il nous est naturel de penser d'apr\`es les habitudes qu'elles nous ont fait prendre. Nous pensons par elles; r\`egles de nos jugemens, elles font nos connoissances, nos opinions, nos pr\'ejug\'es : en un mot, elles font en ce genre tout le bien et tout le mal. \cite[p. 142]{Condillac:1780}

``Algebra is an astonishing proof that the progress in sciences depends only on the progress in languages.''

([L]'alg\`ebre est une preuve bien frappante que les progr\`es des sciences d\'ependent uniquement des progr\`es des langues.\cite[p. 189]{Condillac:1780} 
}  

A similar argument applies to Lobachevsky's brand new Hyperbolic trigonometric functions designed by Lobachevsky on the basis of his intuitive synthetic geometrical considerations that led him to his discovery of Non-Euclidean geometry. Lobachevsky might reduce the synthetic part of his new geometrical theory and present it in the form of a purely analytic calculus as he did in his \emph{Imaginary Geometry} of 1835. But he nevertheless believed that only the intuitive synthetic foundations of his geometrical theory provide the full epistemic access to its content and demonstrate its scientific significance. 
 
An interesting further remark concerning mathematical Synthesis is found in Lobachevsky's \emph{Teaching Plan} of 1925:

``Synthesis will always be a fruitful source for [creative] mathematicians but only a Genius is allowed to open and use it. Teaching mathematics should not draw upon this source. '' \ref{apC}

The timing of writing suggests that at this point Lobachevsky started to seriously think about the possibility of rebuilding the foundations of geometry and considered the possibility of including some of his ideas in his university courses. As we can read in the above note, he quite reasonably ruled out this option which would be hardly welcomed by his colleagues in Kazan. The available textual evidence suggests that Lobachevsky never included the Non-Euclidean geometry in his regular teaching courses and only privately discussed it with some interested students\footnote{Memories of many Lobachevsky's students were collected by Modzalevsky and published in \cite[p. 611-662]{Modzalevsky:1948}}. 

Lobachevsky's basic pedagogical position is that the synthetic foundations of geometry should be learned in a simplified Euclidean form in gymnasia (that is, before students enter the university), while at the mathematical courses of the university should focus on analytic methods referring to the synthetic background. Studying foundations of mathematics as a serious research subject Lobachevsky reserves to mathematical Genius. Fortunately for the posterity Lobachevsky had enough self-esteem to pursue such a study during all of his adult life and career in spite of the lack of recognition on the part of colleagues in Kazan, Saint-Petersburg or elsewhere in Russia.

\subsection{Language}  \label{sum6}
Lobachevsky never extensively wrote on matters of  linguistics but in his  public speech \emph{On the most important aims of education} (1828) \ref{apD} and the \emph{Instructions for teachers of mathematics and physics in gymnasia} (1830) \ref{apE} he makes some relevant remarks, comparing natural languages and the artificial symbolic language of Algebra. As one may expect Lobachevsky stresses advantages and the importance of Algebra in this context:

``What is the reason of recent brilliant successes of mathematical and physical sciences, which are the glory of our age and the triumph of the human mind? Without any doubt these successes are due to their artificial language. For how else all those symbols of various calculi can be called if not a special and very concise language, which expresses extensive concepts with a single stroke without tiring unnecessary our attention.'' \ref{apD}

All these Lobachevsky's remarks are wholly in line of Condillac's \emph{Logic} \cite{Condillac:1780}, \cite{Condillac:1805}

\subsection{Infinity, Irrationals and Imaginaries}  \label{sum7}
Here we are coming to the most controversial part of Lobachevsky's views of mathematics. On the one hand, as we have already stressed in \ref{sum1},  \ref{sum2}, and \ref{sum3} Lobachevsky thinks about mathematics as a proper part of the natural science and wants to revise its foundations to the effect of rejecting all mathematical concepts which are detached from sensual experience and hence from the contact with the Nature:

``[A]ll those concepts, which could not be acquired by our senses \textemdash\ for example, the infinity of space and time  \textemdash\ should be rejected. Those [mathematicians] who attempted to introduce such like concepts in mathematics, did not have their followers.'' \ref{apINF}\footnote{The last sentence of this quote is, of course, rhetorical, see \footref{followers}.}

Lobachevsky extends this argument to irrational numbers which, in his opinion, ``exist only in symbols of the Analytics but not in the Nature''; on this basis Lobachevsky concludes that the traditional theory of incommensurable magnitudes is  ``dry and absolutely unnecessary'' \ref{apB}.

Similarly (but this time less explicitly) Lobachevsky rejects infinitesimals and attempts to develop finitary versions of the Differential and Integral Calculi.

Notice that the above Lobachevsky's critique of ``ideal'' mathematical concepts that are not straightforwardly instantiated by corresponding natural objects perceivable by human senses is in line of his critique of traditional geometrical concepts (like that of line and point), which in Lobachevsky's view are similarly too abstract and too remote from the sensual experience \ref{sum2}, \ref{sum3}. 

But when Lobachevsky comes to the imaginary numbers (like the square root of minus one) he changes his tone:

``\emph{Imaginary quantities} which are so called because they comprise the square root of the negative unit, greatly help Trigonometry by shortening [some computations] and making solutions of certain problems easy and direct. If these quantities would be thrown away from Analysis because of the non-reality of their values then much [of the power of this Analysis] would be lost. [$\dots$] Analysis is, so to speak, only a game of symbols, which does not limit itself to actual (that is, rational) values but  involves [symbolic] expressions, which do not admit for any [real] value. Using signs for non-extractable roots [like $\sqrt{-1}$].  Analysis  extends its operations to numbers, which strictly speaking do not exist.'' \ref{apC}

This passage from Lobachevsky's \emph{Teaching Plan} of 1825 reveals an internal tension in Lobachevsky's thinking. Following the pattern of rejecting ``ideal'' and ``unnatural'' mathematical concepts  Lobachevsky had to reject the imaginary numbers along with  the irrational numbers. He doesn't do it because he realises that using the imaginary numbers Analysis becomes more powerful. In fact, the complex exponential function appears in the very beginning of Lobachevsky's theory of Non-Euclidean geometry; it is an essential element of the analytic calculus that Lobachevsky has designed for this new geometry, namely, the Hyperbolic Trigonometric calculus. Had Lobachevsky forbid himself to apply imaginary (complex) numbers and complex functions he would hardly make his geometrical discoveries. But the same passage shows that the acceptance of imaginary numbers comes for Lobachevsky with a price. Accepting  the imaginary numbers (and preserving his general naturalistic stance) he (dis)qualifies Analysis as a mere ``game of symbols''. Which is at odds at least with Lobachevsky's later Condiliacian notion of Analysis (and, more specifically, of Analysis in the form of algebraic symbolic calculi) as \emph{the} principal form and vehicle of human knowledge about the Nature. 

 Lobachevsky might ask himself: How does it happen that a symbolic trick, that corresponds to no observation and no measurement in the Nature provides Analysis with an extra power in the natural domain of its application, that is, in the Physical geometry? Were Lobachevsky alive in the 20th century he might illustrate the same question by pointing to Quantum mechanics. Lobachevsky did not have any good answer to this question just as we still don't have it today.  
 
Lobachevsky's ambiguous attitude to \emph{imaginary} mathematical objects and constructions can be also seen in how he applied this title to his Non-Euclidean geometrical theory. On the one hand, he hoped that his theory would apply in astronomy \cite[vol.1, p.207-210]{Lobachevsky:1946-51}; when he realised that astronomical evidences supporting his theory in its intended role of mathematical theory of physical space are lacking, he speculated that it could be applied instead in Molecular physics at the microscopic scales \ref{apH}. Nevertheless he also conceived of the possibility to think of his new geometry as an \emph{imaginary} theory presented in the form of analytic calculus (which we call today Hyperbolic Trigonometry). This is how he presented it in his 1835/1837 paper  \cite[vol.3, p. 16-138]{Lobachevsky:1946-51}. Even if the notion of \emph{Imaginary Geometry} could not fully satisfy Lobachevsky's ambition as a scientist, and did not fully represent Lobachevsky's geometrical ideas, he coined this title, which became popular and helped other people to learn about the possibility of Non-Euclidean geometry.

\subsection{History of Mathematics and History of Science}   \label{sum8}
Lobachevsky did not make any significant contribution to historical studies but in his thinking about mathematics and science historical aspects played an important role. He makes historical digressions in the \emph{Foundations} of 1829-1830 \ref{apF}, \emph{Algebra} of 1834, \emph{New Foundations} of 1835-1838 \cite[vol.2, p. 147-588]{Lobachevsky:1946-51} and in some other places; in his \emph{Teaching Plan}  of 1825 he stresses the importance of teaching the History of Science to students. 

 Lobachevsky's basic historical picture is that he's living through the time of continuing scientific revolution started by Galileo, Descartes, Bacon and Newton. 
 In his view this revolution occurred not only in physics and other natural sciences but also in mathematics. But this revolution, in Lobachevsky's opinion,  somehow left aside the foundations of geometry: 
 
 ``[S]ince Newton and Descartes all of mathematics became Analytics and progressed very rapidly leaving  behind the [Euclidean] doctrine, which it no longer needed. Since then this doctrine attracted little attention [on the part of the mathematical community]. This is why Euclid's \emph{Elements} preserved [until today] all their original flaws $\dots$.'' \ref{apF}

 Notice that this historical setup does not imply, by itself, the need to reform the Euclidean ``doctrine''. If the Analytic geometry were self-sustained then the old synthetic geometrical doctrine could be simply abandoned rather than replaced. As we have seen in \ref{sum5} Lobachevsky knew about this epistemological position and was even sympathetic to it but ultimately he disagreed. He admitted the epistemic supremacy of Analysis over Synthesis but he still believed that Synthesis was indispensable in geometry. Hence his zeal to rewrite Euclid's \emph{Elements} but not simply leave them behind.

\subsection{The Problem of Parallels}  \label{sum9}
As we have already stressed Lobachevsky's research on parallels and his discovery of Non-Euclidean geometry were byproducts of his more ambitious project aiming to reforming the foundation of geometry. He mentions the issue of parallels in his \emph{Teaching Plan} of 1825 and explains why, in his view, it worths to be studied:

``The parallelism of [straight] lines presents a difficulty of a different sort, which remains so far unresolved. Nevertheless this difficulty involves sensible truths, which, beyond any doubt, are so important for Science that they cannot be bypassed. '' \ref{apC}

According to our reading, in this passage Lobachevsky justifies his attention to the problem of parallels by pointing that this problem ``involves sensible truths'',  i.e., facts perceivable by humans senses. In Lobachevsky's eyes this feature makes the problem of parallels important \textemdash\ unlike, say, the problem of incommensurable magnitudes, which Lobachevsky disqualified as a pseudo-problem  on the same grounds \ref{sum7}.  Saying that the problem of parallels  ``involves sensible truths'' Lobachevsky means that the hypothesis about the Euclidean character of the physical space allows for an empirical verification and falsification.

In the \emph{Introduction} to his \emph{Foundations} of 1829-1830 Lobachevsky mentions the treatment of parallels among other features of the traditional Euclid-style geometry, which he judge unsatisfactory, and which in his opinion should be fixed. Recall from \ref{sfoundations} that a significant part of this work is reserved for discussing astronomical data pertaining to the question of whether or not the physical space is Euclidean at the astronomical scales. Since Lobachevsky considers his geometrical theory as a part of the Natural Science such an intrusion of empirical argument into geometrical reasoning does not present to him a particular methodological  trouble.\footnote{Lobachevsky's astronomical argument was later criticised by Henri Poincar\'e \cite[Ch.3]{Poincare:1902}  and Ernst Cassirer  \cite[Ch. VI]{Cassirer:1921} to mention only few. This critique is briefly overviewed in the end of the next Section \ref{hilb}.}

The problem of parallels progressively attracted more of Lobachevsky's attention during his career and became central in Lobachevsky's late  writings. 
While in the \emph{Foundations} of 1829-1830 this problem is mentioned as only one defect of the traditional Euclid-style geometry among its other defects (such as the abstract character of its basic concepts, see \ref{sum3}), the \emph{New Foundations} of 1835-1838 already begins with a discussion of the problem of parallels and an overview of earlier attempts to solve it. The \emph{Geometrical Researches} of 1840 and the two editions of \emph{Pangeometry} of 1855 and 1856 are wholly focused on the problem of parallels. This is why we don't find in these late Lobachevsky's writings his  Dimension theory and any other significant trace of general ideas of reforming the foundations of geometry in an empiricist vein.

\subsection{Mathematics and Philosophy}  \label{sum10}
To conclude this Section it is appropriate to add that many of his general thoughts about science and mathematics Lobachevsky would likely not qualify as a  \emph{Philosophy}.  In his \emph{Teaching Plan} of 1825 he says that

``Mathematics must be wholly independent of Philosophy. One can say that philosophy ends where mathematics begins.'' \ref{apC}

But in his writings it is not always clear where philosophy ends and mathematics begins. In particular, it is not clear which parts (if any) of his \emph{Introductions} to \emph{Foundations} and  \emph{New Foundations} (\ref{apF}, \ref{apH}) are in Lobachevsky's own eyes philosophical rather than mathematical.

It should be also born in mind that in the early 1800s the contraposition of science and philosophy (which we also find, for example, in Condillac's writings), was a part of the European intellectual climate of the time.  Remarkably the first volume of the \emph{Course of Positive Philosophy} by August Comte, which appeared in 1830, was titled the \emph{Mathematical Philosophy} \cite{Comte:1830}. Lobachevsky ideas about science and mathematics perfectly fit into this early wave of the Positivist movement.

\section{Hilbertian and Lobachevskian Perspectives on the Hyperbolic geometry} \label{hilb}

As it is well-known, Lobachevsky's and Bolyai's pioneering works in Non-Euclidean geometry published in the 1830s were followed by Riemann's \emph{Habilitation thesis} of 1854 (first published posthumously in 1867 \cite{Riemann:1854}), two Beltrami's papers of 1868 \cite{Beltrami:1868a}, \cite{Beltrami:1868b}, Klein's papers of 1870's \cite{Klein:1871}, \cite{Klein:1873}\footnote{See also his lecture course on Non-Euclidean geometries given in 1889-1890 \cite{Klein:1928}} and, finally, Hilbert's \emph{Foundations of Geometry} first published in 1899 which forged the way in which both Euclidean and Non-Euclidean geometries are usually taught today at the university level \cite{Hilbert:1899}.\footnote{For the history of the 19th century geometry see \cite{Gray:2007}.} How exactly Lobachevsky's publications and his singular approach performed during this development is an interesting and so far little studied question which, however cannot be covered in the present work. In this Section we only try to describe more exactly the difference between Lobachevsky's and Hilbert's views of Non-Euclidean geometry, and understand its origins. 

A the first look, Hilbert's view of geometry appears very similar to Lobachevsky's. Like Lobachevsky Hilbert thinks of geometry as a natural science which builds its foundations on empirical grounds:

``Geometry is the science that deals with the properties of space. $\dots$ Space is not a product of my reflections. Rather, it is given to me through the senses. I thus need my senses in order to fathom its properties. I need intuition and experiment, just as I need them in order to figure out physical laws, where also matter is added as given through the senses.'' (From the Introduction to Hilbert's lecture course of 1891, quoted by \cite[p. 84]{Corry:2004})

Moreover, in 1893-1894  Hilbert like Lobachevsky believed that the Axiom of Parallels can and should be empirically tested.
\footnote{See \cite[p. 87-88]{Corry:2004}. Hilbert's referred to what he believed were Gauss' attempts to empirically check the Euclidean character of physical space using geodesic measurements. Whether or not Gauss made these measurements for that purpose remains today a matter of historical controversy \cite[p. 88, footnote 22]{Corry:2004}. Apparently Hilbert was unaware at this point about Lobachevsky's attempt to falsify Euclidean Axiom of Parallels using astronomical data. Given the linguistic barrier and the fact that Engel's German translation of relevant Lobachevsky's works \cite{Engel:1898} appeared later in 1998 this is not surprising.}  
But already during the 1890s Hilbert's view of geometry also takes a very distinctive \emph{axiomatic} turn:
 
``Geometry also [like mechanics] emerges from the observation of nature, from experience. To this extent, it is an experimental science. [$\dots$] But its experimental foundations are so irrefutably and so generally acknowledged, they have been confirmed to such a degree, that no further proof of them is deemed necessary. Moreover, all that is needed is to derive these foundations from a minimal set of independent axioms and thus to construct the whole edifice of geometry by purely logical means. In this way [i.e., by means of the axiomatic treatment] geometry is turned into a pure mathematical science.''
(From Hilbert's lecture notes of course of mechanics delivered in 1892-1893, quoted by \cite[p. 90]{Corry:2004})
 
 There are several points in the above passage, which mark significant divergences between Lobachevsky's and Hilbert's ways of thinking. First, Hilbert's claims that experimental foundations of geometry \textemdash\ unlike experimental foundations of mechanics, for example \textemdash\ \footnote{Cf. ``In mechanics it is also the case that all physicists recognise its most basic facts. But the arrangement of the basic concepts is still subject to a change in perception $\dots$ and therefore mechanics cannot yet be described today as a pure mathematical discipline, at least to the same extent that geometry is. We must strive that it becomes one.'' (\emph{ibid.}) }, are ``irrefutably acknowledged''. So in Hilbert's view what distinguishes geometry from mechanics, physics, chemistry and any other natural science is the alleged advanced stage of development of this discipline, which shifts its focus from empirical observations and experimentations to logical inferences from empirically established axioms; these axioms are supposed to be so well established that it becomes unnecessary to return to the question of their justification over and over again. This passage demonstrates that for Hilbert the axiomatic form of scientific knowledge is an epistemic ideal that Hilbert hopes to achieve and implement not only in geometry and other mathematical disciplines but also in physics and other natural sciences.  And we know that Hilbert indeed pursued this ideal throughout his career spending a significant part of his time for the project of axiomatising physics and in 1900 included this problem in his famous list of 23 open problems \cite{Hilbert:1900}) .
 
 It is not immediately clear how Hilbert's axiomatic epistemic ideal and his claim of the empirically irrefutable character of geometrical axioms combine with Hilbert's interest to empirical verification of the Euclidean Axiom of Parallels. It is equally unclear how the idea of geometrical axiom as empirically irrefutable truth combines with Hilbert's notion of axiom as a propositional scheme that may admit for multiple interpretations only some of which turn this scheme into a true proposition, that he developed during the same early period of his career. \footnote{See letter exchange between Hilbert in Frege in \cite[p. 31-51]{Frege:1980}}). Just like in Lobachevsky's  case, we should rather admit that Hilbert did not have a fully consistent set of ideas and beliefs about these general matters but often pursued contradictory epistemic goals and made contradictory statements about general epistemological matters. 
 
The internal tensions in epistemic views of Lobachevsky and Hilbert did not prevent these mathematicians from designing research programs according to their philosophical motivations. In spite of some common background these philosophical motivations were very different, and the corresponding research programs pointed to different directions. While for Hilbert the axiomatic architecture of mathematical theories represents an epistemic value which he tries to extend to physics, Lobachevsky, following d'Alembert, is very suspicious about the axiomatic method, and avoids it in his mathematical work. Lobachevsky's rejection of axiomatic method goes along with his critique of traditional geometrical concepts like (breathless) lines and zero-dimensional points, and he wants to replace them by some other concepts that could make more physical sense. So at least in this respect Lobachevsky's project of revising traditional Euclide-style foundations of geometry is more ambitious and more radical than Hilbert's. Both mathematicians want to narrow the gap between mathematics and physics. But they pursue this common goal in different ways and expected to obtain different results. Hilbert aims at the new axiomatic physics designed after the model of the axiomatic geometry. Lobachevsky, on the contrary, finds inspirations of his new geometry in the modern experimental science of Bacon, and aims at the complete integration of geometry and other mathematical sciences in the form of analytic symbolic calculi with physics and other natural sciences. While Hilbert values the logical purity of ``pure mathematical sciences'' and hopes to logically purify physics after the example of his axiomatic geometry, Lobachevsky along with Bacon and Condillac stress the importance of ``learning from the Nature'' using observations and experiments. In \ref{sum3} we showed how exactly Lobachevsky applies this strategy to the foundations of geometry. While Hilbert's conception of logic is broadly Kantian \footnote{By ``broadly Kantian'' conception of logic we mean here any conception which takes logical rules and principles to be a priori, that is, not learned from experience and not revisable on the basis of empirical data.}, Lobachevsky's conception is rather Condillacian: it construes logical analysis as a natural cognitive capacity of humans reinforced by natural and artificial languages including symbolic algebraic calculi, see \ref{sum5}, \ref{sum6}). Lobachevsky's Condillacian conception of logic is that of ``impure'' logic of scientific inquiry ready to revise its principles under the weight of new observations and experiments. 

Hilbert's axiomatic rendering of Euclidean and Non-Euclidean geometries  was a great success which in the 20th century triggered many important developments. The axiomatic method in its modern Hilbertian form still serves today's mathematics professors for teaching their students about Non-Euclidean geometry. It is regrettable, however, that the power of Hilbert's axiomatic thinking obscured other views on that subject including Lobachevsky's. Here are two examples.

In 1902, very soon after Hilbert's \emph{Foundations of Geometry} were first published, Henri Poincar\'e published his famous collection of essays \emph{Science and Hypothesis} where he among other things presents Lobachevsky's achievements and criticises his attempt to use astronomical data for the empirical verification of Euclidean Axiom of Parallels. \cite[Ch.2]{Poincare:1902}  

In this work Poincar\'e wholly misrepresents Lobachevsky's way of thinking by explaining his geometrical discoveries using a Hilbert-style  axiomatic presentation of Hyperbolic geometry; he supplies this axiomatic presentation with his original (Euclidean) models which bear today Poincar\'e's name.\footnote{When and how exactly Poincar\'e first learned about Non-Euclidean geometry remains unknown. In 1870s as a student of \'Ecole Polytechnique he could read French translation of Lobachevsky's \emph{Researches} published by Ho\"uel back in 1866, see \cite[p. 277-278]{Gray:2007}. But in his \emph{Science and Hypothesis} Poincar\'e's presentation of Hyperbolic geometry is Hilbert-style axiomatic; references to Hilbert found in this work make it clear that by 1902 Poincar\'e was already familiar with Hilbert's \emph{Foundations of Geometry} of 1899 \cite{Hilbert:1899}.} 
 
 This misrepresentation has consequences for Poincar\'e critique of Lobachevsky's astronomical argument. Poincar\'e argues as follows. Suppose that measuring the angle sum of a triangle formed by three remote stars an astronomer gets a values smaller than $pi$ radians.\footnote{Lobachevsky failed to find such an example in his contemporary astronomical literature but this fact is irrelevant to Poincar\'e's argument.} Then the astronomer has two options: (i) to explain this empirical observation by the fact that the physical space at this scale is Hyperbolic or (ii) to explain the same observation that the light rays, which usually are taken to be straight, at the large astronomical scales bend (in the way, which mimics the Hyperbolic behaviour) but the underlying space remains Euclidean. In Poincar\'e's opinion, the last option is preferable because it allows one to stick to the familiar Euclidean geometrical concepts and the associated spatial intuitions. 

Poincar\'e's argument is smart but it doesn't go through as long as one admits Lobachevsky's conception of geometrical space. Poincar\'e's conception assumes that geometrical properties of physical space are fixed first (as Poincar\'e's suggests, ``by convention'') and physical laws concerning objects inhabiting this space are established afterwards on empirical grounds. Lobachevsky's conception of geometrical space disallows this order of theory-building. In Lobachevsky's view, geometrical properties of spaces are determined by forces which are involved in the physical interactions of bodies forming these spaces \ref{sum1}. Thus Lobachevsky, unlike  Poincar\'e, does not allow for an unmotivated choice of a background geometrical theory that could be later used in a physical theory. In Lobachevsky's view, geometry is itself a part of physics. The notion of bent light ray might make sense to Lobachevsky only if it were both physically and geometrically grounded like in Einstein's General Relativity.

In 1921 Poincar\'e's critique of Lobachevsky was praised by Ernst Cassirer in his philosophical analysis of the General theory of Relativity and Non-Euclidean geometry \cite[Ch. VI]{Cassirer:1921}. Once again Cassirer (as a good disciple of Kant) takes it for granted that one's conception of space is fixed a priori, whether this conception is Euclidean or not. On this ground Cassirer qualifies Lobachevsky's empirical considerations as an obvious methodological fallacy. Once again, this is a sheer misunderstanding of Lobachevsky's naturalistic way of thinking about the concept of geometrical space. Cassirer might be right or wrong rejecting Lobachevsky's naturalism but at any rate Lobachevsky's epistemological ideas are not naive and deserve more philosophical consideration.

\section{Lobachevsky's Views of Geometry and the Future of Mathematics} \label{future}

Neither Hillbert, nor Poincar\'e, nor Cassirer would agree without reservation with Marshall Stone's and Jean Dieudonn\'e's  bold claims to the effect that mathematics is  ``entirely independent of the physical world'' \ref{intro}. But Hilbert's axiomatic rendering of Non-Euclidean (and Euclidean) geometry was an essential ingredient of the 20th century intellectual history, which helped this view of mathematics to emerge and make it respectful and popular in the early 1960s. The synthetic axiomatic building of a mathematical theory \textemdash\ no matter whether it proceeds in the old-fashioned geometric Euclid-style or in the modern abstract Hilbert-style \textemdash\ always allows one to conceive of the given theory as emerging \emph{ab ovo} and producing its peculiar ``world out of nothing'' \cite[the title]{Gray:2007} which may interact with other similar worlds in intellectually fruitful ways. Questions concerning the historical and cognitive origins of the involved mathematical concepts can be in this case either wholly ignored or safely placed in the \emph{context of discovery} \cite{Reichenbach:1938} and left to historians, philosophers and other non-mathematicians. A similar rhetorical technique allows one to remove out of the scope of the pure mathematics all questions concerning its applications in the natural sciences and engineering. This is the context in which the ``effectiveness of mathematics in the natural sciences'' begins to appear ``unreasonable'' and even miraculous \cite{Wigner:1960}.\footnote{Notice that declarations by Stone \cite{Stone:1961} and Dieudonn\'e \cite{Dieudonne:1961}, on the one hand, and Wigner's famous paper \cite{Wigner:1960} appeared almost simultaneously. In our view the historical dimension is crucially important for understanding the ``Wigner problem''.}  

No similar move is possible with a methodology of mathematics that primarily relies on Analysis  \textemdash\ not only as a way to mathematical results but as a self-standing form of mathematical knowledge that allows one to apply it in the natural sciences and elsewhere. Unlike Synthesis (in the form of axiomatic theories) Analysis does not allow one to begin a mathematical reasoning \emph{ab ovo}. Analysis is always an analysis of something already \emph{present} \textemdash\ be it sense-data, practical and theoretical problems, certain prior beliefs, opinions or prejudges. It does not produce ``worlds out of nothing'' but helps us to solve problems in the world in which we live.

As we already stressed claims of autonomy of the pure mathematics made in the early 1960 by some leading mathematicians were immediately objected by some other mathematicians and physicists. More recently some mathematicians urged for a new ``reunion'' of mathematics and physics \cite{Kapustin:2010}, \cite{Bah:2024}. There are reasons to believe that recent developments in the foundations of mathematics will help to realise this project \cite{Rodin:2021}. 

In this paper we did not try to provide arguments pro or contra the claim of mathematical autonomy. Our goal was more modest: we tried to correct a historical error. Recall how Marshall Stone grounded his claim of the mathematical autonomy using Lobachevsky's legacy. He apparently assumed that making his discoveries Lobachevsky already somehow leaned to this view of mathematics. Perhaps Stone believed that Lobachevsky's discoveries simply could not be made otherwise. We have shown that this assumption is erroneous: Lobachevsky discovered Non-Euclidean geometry trying by all means to merge mathematics and physics more tightly but not to divorce them. Seeing Lobachevsky's achievements exclusively through the lens of Hilbert's axiomatic reasoning and the following Structuralist trend in mathematics that in the middle of the 20th century eventually led to the Bourbaki Manifesto \cite{Bourbaki:1950} not only impoverishes and distorts our understanding of history but also narrows the conceptual horizon for the future development of Science. Rather than being ``lost in transition''\footnote{\emph{Lost in Transition: Topics Forgotten in the Turn to Modern Mathematics} is the title of conference held in Seville on March 5-6, 2026 where an earlier version of this paper was presented. I thank Jos\'e Ferreir\'os and other organisers of this conference for their hospitality and valuable comments.} Lobachevsky's views of geometry may help us today to build mathematics and science of the future.

 \appendix  \label{ap}

\section{Introduction to \emph{Origin and Diffusion of Sounds in the Air} (1823)} \label{apA}

 \cite[p. 396-397]{Lobachevsky:1976}

\label{physmath1}
A mathematically illiterate reader, or a reader familiar only with basic Mathematics, very likely does not understand the power of Mathematics and the extension of its subject-matter. Such a reader will probably doubt that small experiments in a physical laboratory and computations can show us how sounds propagate in the air and allow us to estimate their speed. In what follows the reader will find a proof of this fact. Everything that has a magnitude belongs to the science of number. And what in the physical world does not have a magnitude? Everything in the physical world exists under the necessary condition of being measurable; hence every physical entity is subject to mathematical laws. This is why all Natural Sciences try to attain the high degree of perfection that will be followed by their unification with Mathematics; such a unification enhances their progress. Such a unification [with Mathematics]  already occurred in Physics, and more recently in Mineralogy. We expect that in the future this will also happen in Chemistry. 

When we discover important truths being led by Mathematics, and return to the very origin of these great and unexpected achievements, we can see that the rigorous science is also the science of common sense. The ultimate foundation of science is a fair judgement about things:  it invariably leads a mathematician through his computations. After that there remains no natural phenomenon that he could not account for. There is no phenomenon he could not predict and determine precisely, both their time and their measure. It appears that one's notion of things, and one's fair judgement on things, should not depend on computations. But in fact, a mind familiar with computations goes far beyond the limit where a [uneducated] mind, which is not initiated into the mysteries of the science of number, stops.

[$\dots$]

\section{From \emph{Teaching Plans for the Academic Year 1824/25} (Summer 1824)}\label{apB}

Teaching of pure mathematics 

1. Teaching method in general

[On Analysis and Synthesis]  \cite[173-175]{Modzalevsky:1948}

Recent brilliant successes in Mathematics are due to Analysis. This excellent invention of human mind makes it possible to determine everything numerically, to express all qualities and connections symbolically, to represent all relations with equations, and on this basis solve any problem. A different, geometrical method amounts to representations of everything with lines, surfaces or solids, and solving a given problem by finding relations between lines on a diagram. Most often, the geometrical method is a Synthesis, that is, a composition: knowing a truth [in advance] one invents a geometric construction for its proof.  The analytic method [on the contrary], leads directly to the discovery of truths.  It always involves the same [kind of] swoop towards the solution, the solutions themselves are wide [and general], and the equations, which express the mutual dependences of quantities, contain everything one needs for answering the [proposed] questions, and thus subsume the process of problem-solving to uniform, straight and brief actions. 

The main difficulty of Analysis is due to the wide [i.e. general] and abstract character of its concepts. When one applies Analysis, one's imagination doesn't stop on something unique [say, on some unique geometrical representation] but is forced to embrace many objects at once; the analytic judgement should apply to all these objects taken together. One may rightly fear that the choice of a special case may lead to an error because of a possible confusion of its specific features with the general case. Another disadvantage of Analysis is that  it proceeds only with numbers. This is why analytic conclusions always need an interpretation: some numbers should become bodies, some other numbers should become times or natural forces. The characters [i.e., mathematical symbols used in Analysis] should demonstrate their quality and their [mutual] dependences. Such difficulties increase with widening [= generalsation] of analytic solutions but they can be mitigated with the art of teaching and exercising. Particularly useful in this respect are examples of special cases that can serve as models for all other examples. 

In any event, Mathematics has its present highly developed state due to Analysis; the [synthetic] geometrical method is no longer sufficient [to support the progress of Mathematics]. For this reason in the university [teaching] it is necessary to follow the analytic [method] from the very beginning. This helps [students] to familiarise themselves with the analytic method and makes teaching more uniform. Indeed, [student's] capacity of abstraction may increase progressively by a continuous exercising. 

[$\dots$]

Analysis presupposes that bodies, times and forces have been already measured or at least that such measurements are possible and the appropriate methods of measurement are known. Indeed, Analysis involves only numbers which represent a measure of every magnitude but not these magnitudes themselves.\footnote{By word ``magnitude'' I translate Lobachevsky's word \selectlanguage{Russian}``коликое''\selectlanguage{English} which does no longer belong to today's Russian. In the given theoretical context the term could be alternatively translated as ``quantity'' but this would make Lobachevsky's point more obscure. Today's Russian speaker can easily identify a common root in the outdated \selectlanguage{Russian}``коликий''\selectlanguage{English} and today's \selectlanguage{Russian}``количество''\selectlanguage{English} and \selectlanguage{Russian}``сколько''\selectlanguage{English}. In other fragments Lobachevsky also uses word \selectlanguage{Russian}``величина''\selectlanguage{English}, which is the modern Russian equivalent of English ``magnitude''.} So there should be parts of mathematics that do not admit for Analysis. Such [synthetic] parts must be separated and not extend beyond the point where measurements are determined; beyond this point the supremacy of Analysis must be recognised. A strict application of this principle [in teaching] preserves a large part of [the traditional synthetic account of] foundations of geometry and mechanics. 

[$\dots$]

Analysis consists of the Higher Analysis and the Elementary Analysis. The former involves the Differential Calculus; the latter involves only [elementary] algebraic operations and for this reason can be called Algebraic Analysis. [Augustin Louis] Cauchy calls it in this way, see his \emph{Cours d'analyse alg\'ebrique}\footnote{Lobachevsky's reference is obviously \cite{Cauchy:1821}.} assuming that Algebra is the part of Analysis that does not include the Differential Calculus. The supremacy of the Higher Analysis is well known; what is easy to do using the Differential Calculus is often very difficult or even impossible to do using only Algebra. 

We cannot deny that the geometrical method  in certain cases provides simpler and shorter solutions than the analytic method.  But this occurs only in singular cases, which are, once again, identified with the Analysis. Such cases must be mentioned in teaching. 

[$\dots$]

III Academic disciplines

\cite[177-178]{Modzalevsky:1948}

The \emph{geometrical part} [of the teaching program]  will include [i] Foundations of Geometry, [ii] Linear and Spherical Trigonometry, [iii] Analytic Geometry. 

The Foundations of Geometry appear to some [teachers] so easy that they teach this subject before Arithmetic and Algebra. If we assume, however, that mathematics needs to be rigourous then we should say that Foundations of Geometry is a difficult subject, which should be studied only when [students'] minds, are sufficiently prepared by studying other mathematical disciplines. Ignoring these difficulties in teaching [the Foundations of Geometry] would be an important omission.  [$\dots$] 

The triple measure [i.e. the three spatial dimensions] of bodies is usually treated poorly.  It is impossible to provide a clear conception of length, width and thickness of bodies beginning the [presentation of] geometry from these very notions. Even if in Geometry our senses protect us from false conclusions, it is nevertheless desirable to save a mathematical discipline [namely, Geometry] from the accusation of being unrigorous, obscure and insufficient in its foundations. In addition [to the former argument we can remark that] we don't know which truths are concealed from us in what we don't [properly] understand.\label{measure1} 

 \label{touch1}I believe that geometry studies the [special] property of natural [i.e., physical] bodies, which is called the \emph{touch} (\selectlanguage{Russian} прикосновение  \selectlanguage{English}). [Two] bodies are geometrically identical when they occupy the same place, [i.e.] when they touch the ambient space similarly. In Geometry, [physical] bodies are measured only in respect of their [mutual] touch, so this property should be the subject-matter of this mathematical discipline. Let us now see what the Foundations of Geometry should consist of. 

A distinctive feature of exact sciences is that they put in their foundations certain notions, from which the rest is produced by the power of our judgement. Foundations of Physics can be hypothetical but foundations of the Pure Mathematics should be unquestionable truths, [i.e.,] our first concepts of nature, which being once obtained, are preserved [in our mind] forever, which make part of every mental representation and serve as the ultimate foundation of any [human] judgement. Foundations of Geometry should meet the same requirements. Further, primitive concepts directly apply to the nature. In this respect they differ from complex [concepts], which require the existence of some other concepts wherever they might originate from. Surfaces and lines do not exist in the nature but exist only in the [human] imagination; as a consequence, they [i.e., geometrical surfaces and lines] presuppose a property of [physical] bodies, which had to produce in us the concepts of surface and of line. Nobody so far attempted to reach those origins [of geometrical concepts], and [as a result] the Foundations of Geometry remain obscure. Now it is clear why so many things in Geometry do not stand against a rigorous critical [conceptual] analysis (\selectlanguage{Russian}разбор\selectlanguage{English}). It is appropriate to confess here  that the development of [human] intellect, if we may call by this name the acquirement of the ability to judge on firm foundations, shares the fate of all things human, i.e., it is imperfect. The more it advances the more it requires an assistance and, finally, a help of all our enlightenment.    

\label{3D1} I believe that the foundations of geometry should be derived from the property of bodies, which we discover when we imagine the [physical] space along with the triple division in this space.  First, [imagine] a series [$p^{1}_{1}, p^{1}_{2},\dots, p^{1}_{k-1}, p^{1}_{k}, p^{1}_{k+1}\dots$] of adjacent [bodily, 3-dimensional] parts, which do not touch each other over another such part [i.e., such that $p^{1}_{i}$ touches $p^{1}_{j}$ only if $j = i\pm 1$]. Then in one of such parts [$p^{1}_{i}$] we consider a similar series of new parts   [$p^{2}_{i1}, p^{2}_{i2},\dots, p^{2}_{ij},\dots$], each of which touches two parts [$p^{1}_{i-1}$ and $p^{1}_{i+1}$] of the first division. Finally, in an element [$p^{2}_{ij}$] of this latter series [$p^{2}_{i1}, p^{2}_{i2},\dots, p^{2}_{ij},\dots$] one can imagine a further series  [$p^{3}_{ij1}, p^{3}_{ij2},\dots, p^{3}_{ijk},\dots$], each element [$p^{3}_{ijk}$] of which touches upon parts [$p^{1}_{i-1}$, $p^{1}_{i+1}$ and $p^{2}_{i(j-1)}$, $p^{2}_{i(j+1)}$] of the first and the second divisions. One cannot, however, imagine further divisions [of the given body] into parts, which would satisfy the same conditions. In short, it is impossible to assign to a body [like a cube] more than six faces when we imagine that this body is carved out of the [ambient] space using such a [triple] division [as described above]. Each division will reveal two faces, and after the third division the generation [of the body] will be complete\footnote{See corresponding diagrams in \ref{apH}}.

Further, one can consider three [kinds of] touches between bodies.  Consider, for example, the touch of a prism to a flat body [i] through its base, [ii] through its edge and [iii] through its vertex. \label{measure2}The measurements of bodies similarly can be of three different kinds accordingly to their touch: there is one [kind of] complete measurements and two [kinds of]  incomplete measurements (one kind of measurement is yet less complete than the other).  In incomplete measurements one is allowed to ignore those parts of two touching bodies, which do not [directly] touch another body.\footnote{Think of measuring a table using a ruler. The only relevant feature of the ruler concerns the intervals between its marks, which touch the table during the measurement. The width and the thickness of the ruler in the context of this measurement are ignored.} One may ignore such parts in one's thoughts as useless and reach in one's imagination the thinness of a paper sheet, sewing thread or a point marked by pen on a paper. This is how we usually imagine [geometrical] surfaces, lines, and points, which only simplifies our reasoning about [spatial] measurements relatively to the touch. When a [agricultural] field is measured by touching it with a square board the thickness of the board is not taken into considerations but the measurement proceeds, nevertheless, with bodies [but not with abstract geometrical entities]. 

[My new prospective] Foundations of Geometry will comprise only properly geometrical principles, i.e., only those principles, which lead directly to general rules of measurement for lines, surfaces and bodies. Everything else will be left to application of Analytics [i.e. to the Analytic Geometry].  

[$\dots$]
      
[$\dots$] 

\cite[182-184]{Modzalevsky:1948} \label{apDIFF}

Differential Calculus

[$\dots$] 

The Differential Calculus is introduced [in my course] into Geometry and Mechanics in the form of calculus of [finite] increments provided that higher degrees of increments are discarded. This new [finitary] form does not make the differential calculus less rigorous because there exists, indeed, a passage from one calculus to the other [i.e., from the standard differential calculus to the proposed calculus of finite increments]. Considering the aforementioned rule of discarding higher degrees of increments it is easy to see that the computation of [finite] increments and the application of the Differential Calculus will always bring the same result. Thus the rigorous character of the Differential Calculus is not compromised [in its geometrical and mechanical applications].  We should blame [of the lack of rigour in the Differential Calculus] not the computations but rather our [ineffective] concepts of geometrical and mechanical magnitude, which we cannot measure either in our mental representation [i.e., in the imagination] or in the actual practice via approximation. This is why Analytics puts limits of all such approximations aside of the true numerical values, so that these values may agree with any actual measurement: the more precise is the measurement, the better [is the agreement]. Moreover the true values should agree with the nature itself, which is unimaginably refined in its foundations. 

But the nature itself tends to avoid our computations. This certainly happens in Geometry. In Mechanics, however, one of the following two options is possible: either all our [mechanical] concepts are false, or, like in Geometry, we are able to compute [in Mechanics] only approximately. Indeed, in the nature there exists a motion of attracted bodies. But do things that we call velocity and [spatial] extension, which can be compared with repeated hits over a body, really exist? \footnote{It is not quite clear what Lobachevsky means here by ``repeated hits'' (\selectlanguage{Russian}повторяемые удары\selectlanguage{English}) but apparently he points in this passage to the problem of measuring time and to Zeno paradoxes.} The agreement between our computations [and the empirically observed mechanical phenomena] forces us to believe that this is the case [i.e., that velocities and spatial extensions exist] but this agreement [between theory and experiment] may have also different possible causes [unrelated to the reality of motion and spatial extension]. In this case it is possible that our computations are strictly correct but their [received] foundations are false: our concepts [like extension and velocity] can be artificial but for some unknown reason still capable to replace the true concepts.  

Anyway, such speculations, provided that they are approved [by the university administration], could be useful only for those [students] who have already studied Mathematics and Mechanics sufficiently and [passed the final exams] successfully [$\dots$]. Such speculations are not appropriate in teaching when the most important value is clarity. I talk to students about [indivisible] elements of bodies and forces as if they were actually existent. I remark that these elements must be very small, so they can be represented with differentials, while a whole composed of these elements should be obtained via the integration. Everybody knows how much clarity brings the introduction of [indivisible] elements to Geometry and Mechanics, how simple and intuitively clear every problem  becomes in this case. To neglect this simplicity in teaching pursuing an imaginary [chimeric] rigour and try, without any need, to liberate oneself from the geometrical point of view, would mean to deprive students from an important benefit.\footnote{Compare in d'Alembert: ``Un inconv\'nient peut-\^etre plus grand que celui de s'\'ecarter de la rigueur exacte que nous y recommandons, seroit l'entreprise chim\'rique de vouloir y chercher une rigueur imaginaire. Il faut y supposer l'\'etendue telle que tous les hommes la con\c{c}oivent, sans se mettre en peine des difficult\'s des sophistes sur l'id\'ee que nous nous en formons, comme on suppose en m\'echanique le mouvement, sans r\'epondre aux objections de Zenon d'El\'ee.'' [Perhaps a greater drawback than deviating from the exact rigour we recommend would be the chimerical undertaking of seeking an imaginary rigour within it. We must assume extension as all people conceive it, without concerning ourselves with the difficulties the sophists have with the idea we form of it, just as in mechanics one assumes motion without addressing the objections of Zeno of Elea.] \cite[vol.7, p. 635B]{Diderot&DAlembert:1751-65}. } I would like to add that the stipulation of [indivisible] elements in nature helps to dissolve the incomprehensible issue of incommensurability. In this case the inexpressible (\selectlanguage{Russian}невыразимые\selectlanguage{English}) [i.e., irrational] numbers become artificial numbers, which exist only in symbols of the Analytics [i.e., as mathematical symbols] but not in the Nature. This is why I always find it useless to talk about the incommensurability. This is a dry and absolutely unnecessary doctrine for Analytics; it is unnecessary in applications provided one follows my method of teaching\footnote{\label{noirr} In his critique of the concept of incommensurability, which underlies Euclid's theory of geometrical proportions \cite[Books 5-6]{Euclid:2008} and which later gave rise to the modern concept of irrational number, Lobachevsky goes further than d'Alembert who describes the concept of incommensurability as a necessary evil. Like Lobachevsky d'Alembert remarks that this concept significantly complicates the very first chapter of the Elementary geometry, namely, the \emph{Longimetry} which takes care of measuring distances. But unlike Lobachevsky d'Alembert believes that this concept is indispensable in geometry, and thus he exempts it from the sin of the ``chimeric rigour'', see  \cite[vol.7, p. 634A-B]{Diderot&DAlembert:1751-65}.}.

[$\dots$]  \cite[186-187]{Modzalevsky:1948}

Teaching Analytic Mechanics

 [$\dots$]
\label{physmath2}
Only the Differential Calculus is able to embrace  Statics and Mechanics and represent all [relevant] problems with a single equation. Once one knows forces which act upon some given bodies, the solutions of that equation always bring knowledge of equilibrium and motion of the bodies. Many words have been written on the Differential Calculus, and mathematicians made a lot of efforts to prove its correctness. The Analytics, which is a proper creation of human mind, has its sources in the indisputable primitive and simple concept of magnitude, which is impressed on our mind when it develops. These concepts always serve us as a firm foundation of all [our] judgements. There cannot be anything contingent in this [fundamental] concept of magnitude: everything here is rigorous, clear and definite. But as soon as one turns from human speculations to the nature itself one notice that his best concepts do not [fully] apply in this domain but serve only as approximations which are quite limited because of the limits of our senses. In [physical] bodies there is no points, no lines and no surfaces. The [physical] substance does not fully fill the void occupied with [physical] bodies.\footnote{This is probably because there remain multiple voids between the hypothetical atomic particles, which, in Lobachevsky's view, compound this substance. The idea that physical substances are composed of atoms \emph{and} the hidden void, dates back to Lucretius and his Greek sources, see \cite[Liber Primus: 417ff]{Lucretius:2008}}. [The visible bodies are] aggregates of particles, which are not accessible for our senses. Opposite [and thus counter-balancing] forces preserve [constant] distances between those particles and do not allow them to approach and touch one another. \label{increments1825}Thus in the nature there is no continuous variation of mutually dependent magnitudes, which is assumed in the Analytics. It follows, strictly speaking, that that the nature admits only for the calculus of [finite] increments but not for the Differential Calculus. But as far as the increments diminish, the two calculi become closer to each other. In nature the difference between the two calculi becomes inconspicuous.  Thus the Differential Calculus takes here the form in which it has been first conceived by its inventors \footnote{Here Lobachevsky likely points to the early form of Differential Calculus found in Leibniz's and Newton's pioneering works. Notice that in 1824 when Lobachevsky wrote these lines  the modern Cauchy-Weierstrass foundations of Mathematical Analysis were still abscent. This allows us to interpret Lobachevsky's proposal to replace his contemporary Differential Calculus with a Calculus of Finite Increments as an attempt to to develop for this calculus a new foundation that would both better serve practical needs and be more mathematically rigorous.} 

Whether or not the new way of reasoning proposed by Lagrange in his \emph{Calculus of Functions} [\cite{Lagrange:1806}] constitutes an improvement in teaching Analytics, it gives no advantages in teaching Mechanics. In Mechanics it would be certainly wrong to deviate from the first primitive [mechanical and geometrical] concepts and thus deprive students from a clearer representation, which allows one so easily to judge about equilibria and motions of every minute part of a given body. I dare to say that Lagrange in his \emph{calculus of functions} made an unsuccessful  attempt to design a new method of teaching [Analytics].  Lagrange didn't want to apply the new method in his \emph{Analytic Mechanics} [\cite{Lagrange:1788-89}] which I use [in my teaching practice].\footnote{\label{lagrange} Lagrange's \emph{Analytic Mechanics} of 1788-89 (in two volumes) is renown for its formal abstract treatment of Mechanics. As the author stresses in the \emph{Preface} to this work ``[T]here is no diagrams in this book. The methods that I applied don't require either constructions or geometrical or mechanical reasoning but only uniform and regular algebraic operations.''  \cite[vol.1 p. vi]{Lagrange:1788-89}. In his  \emph{Calculus of Functions} of 1806 \cite{Lagrange:1806} Lagrange goes even further in the same direction and  applies the same formal approach in the foundations of Differential Calculus (including the Calculus of Variations). Interestingly, Lobachvesky approves on the former Lagrange's textbook but not on the latter. I understand here Lobachevsky's position as follows. On the one hand, Lobachevsky is very enthusiastic about the power of Analysis, and sees Analysis as a higher form of mathematics than the traditional ``synthetic'' form of mathematical reasoning exemplified by Euclid. But on the other hand, Lobachevsky believes that the traditional intuitive forms of geometrical and mechanical reasoning perform in mathematics and natural science a different epistemic function, which is equally or even more important than that of Analysis: to ground mathematics on ``natural'' concepts or at least to forge primitive mathematical concepts in the most \emph{natural} way, so they could better serve for the needs of Natural Science and Engineering. This is why Lobachevsky embraces Lagrange's formal analytic treatment of Mechanics but doesn't want to extend the same approach to the foundations of Differential Calculus where, in Lobachevsky's view, geometrical and mechanical intuitions play a significant epistemic role. See \ref{sum5}}

\section{From \emph{Teaching Plans for the Academic Year 1825/26} (Summer 1825)}\label{apC}

[On Analysis and Synthesis]  \cite[201-203]{Modzalevsky:1948}

I. Teaching method in general

Beyond any doubt, the best method of teaching mathematics is the analytic method, which is practiced in the Kazan University \textemdash\ except those parts [of mathematics] where this method does not apply as, in particular, in the Foundations of Geometry.  Synthesis, which is  [historically] an earlier [human] invention, remains applicable only in the foundations [of mathematics]. Later [in the history] Synthesis gave way to Analysis and recognised its [epistemic] superiority. This is why at the university level Mathematics is mostly taught by Analytic methods. 

The Analytic method [of mathematical reasoning] amounts to expressing relationships between magnitudes with [algebraic] equations. It has the following advantages: [i] the uniform approach to solving all sorts of questions, [ii] generality [of obtained solutions]; the most important advantage of this method is [iii] that the equations, which express mutual dependences of magnitudes, contain everything that is needed for solving the given question.  Disadvantages of the Analysis are [i] difficulties in understanding, which are due to the generality and abstractness, and, finally, [ii] that Analysis operates with [unsorted] numerical values. When by magnitudes [involved in the Analysis] one understands things, times or natural forces, conclusions of Analysis need yet to be interpreted, i.e., converted from numbers to actual magnitudes. Sometimes this procedure is very difficult and prone to errors. An example of such error is the dispute between d'Alembert, Euler and Lagrange concerning the continuity of the curvature of vibrating strings.\footnote{About this dispute see \cite{Wheeler&Crummett:1987}, \cite{Oliveira:2020}.}  

Synthesis does not have the advantages of Analysis but it does not have its disadvantages either: it is clear, sensible and much more persuasive for beginners. It is necessary, however, to retain [elements of] Analysis in foundations [combining them with Synthesis] in order to preserve the uniformity of teaching and to prepare [teaching] of higher parts of Mathematics. Advantages of Analysis and advantages of Synthesis can be combined by proceeding from several special cases to general cases, and by enriching one's teaching with [concrete] examples. This is a rule of the local Gymnasium (\selectlanguage{Russian}Гимназия\selectlanguage{English})\footnote{In the early 19th century Russia Gymnasia where schools preparing their students to entering an university. Lobachevsky refers here to Kazan Gymnasium from which he graduated in 1806. }

Nevertheless, there are certain parts of mathematics where Synthesis is the only available method. Synthesis should prepare a given discipline to become a perfect subject to Analysis as this happens in Geometry and Mechanics. A systematic teaching requires a strict separation of such [synthetic] parts of mathematics [from analytic parts] in order to make explicit their specific features and their sources. At the same time, Foundations of Geometry [built according to the above condition] would be too short, and too far remote from Analysis, where Geometry is, in fact, widely applied. 

This is why I find it natural to divide the course of pure mathematics into two parts: one of which is taught in the Gymnasium and the other, the complete one, which is taught in the University. The \emph{gymnasium} teaching program comprises the Foundations of Algebra, the Synthesis of Geometry, and some applications of Analysis in Geometry including various special cases and examples from the common life. Thus this course has a double benefit [including Synthetic Geometry, elements of Algebra, and applications of the latter to the former]. The \emph{university} program once again starts with the Foundations of Mathematics considering them [this time] from a different point of view. It embraces the Foundations in their full scope, separates Synthesis from Analysis, and systematically explores operations of Analysis in their natural order which conducts towards a higher [degree of] artificiality and a higher  [degree of] generality. 

The simplicity and the easiness of Synthesis in certain cases needs to be stressed.  [For example,] the integration of [differential] equation where the ratio of the squares of two differentials of two [given] variables is equal to the ratio of their sinuses can be easily accomplished via the consideration of a spherical triangle while all of Euler's and Lagrange's [technical] wit was needed in oder to find a solution [of the same equation] using Analysis alone\footnote{A similar geometrical technique is systematically by Lobachevsky in his \emph{Application of imaginary geometry to certain integrals}  \cite[vol.3, p. 175-407]{Lobachevsky:1946-51}) }. 

Synthesis will always be a fruitful source for [creative] mathematicians but only a Genius is allowed to open and use it. Teaching mathematics should not draw upon this source. \footnote{Notice the apparently contradictory Lobachevsky's description of the \emph{Synthetic} Geometry. On the one hand, in Lobachevsky's view, this part of geometry constitutes the most primitive layer of geometrical thinking, which is surpassed by the Analytic Geometry. This is why the \emph{Synthetic} Geometry is taught mostly in the primary schools and gymnasia while the university-level mathematical education should focused on the Analytic Geometry. But at the same time, paradoxically, working in the \emph{Synthetic} Geometry is a privilege of a mathematical genius while regular mathematicians typically used universal analytic methods. This double-faced character of Synthetic Geometry anticipates Lobachevsky's latter attitude to his \emph{imaginary} geometry, which is roughly the following. The analytic face of this theory (presented by Lobachevsky separately in \cite{Lobachevsky:1837}) is public: it can be more easily communicated and probably even taught in the university (albeit we don't know whether or not Lobachevsky ever made such attempts). It can be also useful in some external contexts \cite{Lobachevsky:1836}. But the synthetic face of this theory requires much more geometrical ingenuity and imagination and for this reason is reserved only for few. Indeed, in Lobachevsky's lifetime the synthetic presentation of his new geometry was well understood (through \cite{Lobachevsky:1840}) only by Carl Friedrich Gauss.}      

  [$\dots$]

 \cite[203-206]{Modzalevsky:1948}
 
A. Foundations of Geometry \label{apINF}

Here we shall talk about concepts laying in the foundation of mathematical sciences; building conceptual foundations for Geometry is very important. \label{senses1}Without any doubt all our concepts concerning [physical] bodies have their origins in our senses. This is confirmed [i] by the fact that our judgement halts when senses fail to lead us and [ii] that we abstract away from bodies those concepts, which are suggested to us by senses disregarding the essence of those things. Think, for example, about straight and curve lines and surfaces, which do not exist in nature. Our imagination possesses these ideal entities obtained in the shortage of senses. This is why all of our knowledge based on concepts borrowed from nature is true only relatively to our senses.  This [i.e., attaining truths relatively to our senses] is, however, the only purpose of mathematical sciences as far as they remain mathematical, i.e., as far as they deal only with counting and numbers.

 It follows that  [i]  all concepts, whatever they are, which are acquired from the nature, can be admitted  as foundations of mathematical sciences, and [ii] that mathematics [built] on such foundations can be truly called an exact science. On the contrary, all those concepts, which could not be acquired by our senses \textemdash\ for example, the infinity of space and time   \textemdash\ should be rejected. Those [mathematicians] who attempted to introduce such like concepts in mathematics, did not have their followers.  This was the destiny of Kant's Phoronomy and more recently also of the [concept of being] infinite in Analysis.\footnote{About Kant's Phoronomy see \cite[ch. 1]{Friedman:2013}. \label{followers}Lobachevsky's factual claim according to which the notion of infinite space is unpopular is difficult to evaluate from a historical point of view but it is nevertheless plausible as long as it refers to Lobachevsky's contemporaries. But from the anachronistic perspective of today's mathematics this claim appears clearly erroneous. After the invention of Set theory by Georg Cantor in the end of the 19th century, and the wide acceptance of this theory as a standard foundation of mathematics in the 20th century, the notion of infinite geometrical space also became fairly standard.}, of the heterogeneity of lines and angles \footnote{\label{hetero}Compare in \emph{Geometry}: ``Some mathematicians wanted to postulate the impossibility to determine lines with angles as a foundation of Geometry. But such a foundation is defective because heterogeneous magnitudes can be mutually dependent.''  \cite[vol.2, p. 69]{Lobachevsky:1946-51}. Lobachevsky points here to Legendre who attempted to show that the negation of the Euclidean Axiom of Parallels is inconsistent with other principles of Euclidean Geometry because it leads to the ``absolute unit of length'', i.e., to the dependence of lengths from angles (so one can define the ``absolute unite of length'', for example, as \emph{the} length of a side of the equilateral triangle with the sum of internal angles equal to $\frac{\pi}{2}$) \cite[p. 287-294 (Note IV)]{Legendre:1794}. The fact that Legendre's argument is erroneous can be seen from the fact that the dependence of lengths from angles also takes place in the spherical geometry. For further details see V.F. Kagan's commentaries in  \cite[vol.1, p. 64-66]{Lobachevsky:1946-51}. \label{homogen} } Lagrange in the title of his book on the theory of analytic functional presents as its major value the fact that it ``distances itself from all considerations of infinitely small [aka infinitesimal] [...] and reduces Analysis to Algebraic Analysis of finite quantities''.\footnote{See the cover of \cite{Lagrange:1797}. Lobachevsky quotes Lagrange's title in the original French.} This is why all foundations of mathematics, which some mathematicians attempt to produce from the reason alone independently of [real] things in the world, will be always futile in mathematics, and often not even mathematically justified. The uniformity of primitive concepts of all things, their simplicity and their small number, demonstrate that they necessarily follow from the essence of things with respect to the human nature. This is why these primitive concepts will forever be a solid foundations for sciences. 

Uneducated people believe that the fact that bodies fall dawn is a necessary condition of the very existence of these things. In other words, they take this quality of bodies for a primitive concept [without which no body can be thought of]. A little bit of education compels, however,  these people to see here only an appearance and think about its cause. This effort of human mind towards knowing the cause and, as it were, its desire to produce everything from itself [i.e., from the human mind alone], compels mathematicians to reduce the number of concepts to minimum; it is fair to say that this strategy in mathematics has been well rewarded with great successes. 

But the closer concepts stand to each other, the more difficult it becomes to distinguish [artificially] composed concepts from [naturally] acquired [primitive] concepts. This difficulty is not yet overcame in Geometry. In the nature we know only bodies. As a consequence, the concepts of lines and surfaces are derivative [artificial] concepts but not [naturally] acquired concepts, and thus should not be taken as a foundation of the mathematical science. But what distinguishes bodies from other magnitudes which we know in nature? How bodies give rise to the doctrine of lines and surfaces? No geometry can answer this question so far.\footnote{In this last passage Lobachevsky formally contradicts himself by saying, first, that in the nature we know only bodies, and, second, that in addition to bodies, we know in the nature some ``other magnitudes''.  By ``other magnitudes'' Lobachevsky apparently means natural forces, compare \ref{physmath4}. Lobachevsky says in this passage that the concept of body is the only \emph{geometrical} concept that has a physical meaning and for this reason deserves a foundational status in Geometry. Since the ``other magnitudes'' do not belong to Geometry Lobachevsky does not consider them here.}

\label{3D2}\label{measure3} The triple  measurement  [i.e., the number of metric dimensions equal to three] of bodies, which is presently taken for a foundation of Geometry, would not comprise any concept if it would not point to what we sense (without being able to provide an account for it). But until [new] firm and properly mathematical Foundations of Geometry are laid down, which will completely change the mathematics education, [as a matter of pragmatic compromise] the concepts of lines and surfaces should be treated \emph{as if} they were acquired from the nature of things, and be used in the [current] Foundations of Geometry without further clarifications.

\label{parallels1}The parallelism of [straight] lines presents a difficulty of a different sort, which remains so far unresolved. \label{senses2} Nevertheless this difficulty involves sensible truths (\selectlanguage{Russian}ощутительные истины\selectlanguage{English}), which, beyond any doubt, are so important for Science that they cannot be bypassed. It remains to reduce all such truths to one, which could easily convince [anyone] in spite of the lack of rigour with the help of some clarifications.\footnote{\label{P5} The problem of parallels, which is not mentioned in the 1824/1825 teaching plan, appears in the 1825/1826 version of the plan as a side problem in the context of Lobachevsky's project of building new science-friendly foundations of geometry. At first glance Lobachevsky points here to a desired proof of the Axiom of Parallels (treated as a theorem) by using an equivalent or a stronger statement, which has more intuitive appeal. This interpretation of Lobachevsky's words is based on the assumption that the Axiom of Parallels and its known consequences (such as the existence of non-congruent similar triangles) belong to the class of ``sensible truths'', that is, to the class of truths known via the five human senses. But this assumption is not taken by Lobachevsky for granted. He doesn't rule out a possibility that, say, the perceived similarity of certain non-congruent triangles may be only approximate, and that human vision empowered with astronomical instruments may produce new Non-Euclidean foundations of geometry. At least in his later \emph{Foundations} Lobachevsky describes such an attempt \label{foundations}. Thus Lobachevsky says in this passage that the \emph{problem of parallels} has important empirically observable implications and for that reason must be tackled in the foundations of geometry. Lobachevsky does not tell us in this passage, however, whether the new prospective foundations will justify the Axiom of Parallels or refute it. See also the above discussion about Lobachevsky's lost 1826  paper referred to in the \emph{Foundations} \ref{brief1826}.}

This concludes our list of major difficulties in Geometry, which so far could not be resolved. To increase the number of basic concepts (\selectlanguage{Russian}основных понятий\selectlanguage{English}) in Geometry would mean to change the [received method of] teaching [geometry], which is so simple and so natural, without gaining any benefit. This is why one should not follow those who wish to accept as a foundation [of geometry] the principle of similarity, the principle of the heterogeneity\footnote{see \footref{hetero}}, and moreover those who believe that the infinity \label{inf} is a subject-matter of human judgement. I find it useless to focus on such difficulties in teaching or try to approach them philosophically. Mathematics must be wholly independent of Philosophy. One can say that Philosophy ends where Mathematics begins. Once one takes the way of Mathematics [one should stop philosophising because doing otherwise] would mean to go back.\footnote{Claiming the independence of Mathematics from Philosophy Lobachevsky stress its interdependence with Natural Science. Lobachevsky's take on Mathematics, Science and Philosophy shares many common features with the account of mathematics given by Auguste Comte in his \emph{Course of Positive Philosophy}. This cannot be a direct influence, however, because Lobachevsky wrote these lines in 1825 while the first volume of Comte's fundamental work, which treats mathematics, appeared only in 1830  \cite{Comte:1830}. But the two authors rely upon the same sources in their contemporary mathematics and science, and often come to similar conclusions.} 

One can find two kinds of obscurities in Geometry textbooks: 1) obscurities caused by the [authors'] failure to follow the rule according to which everything [in Geometry] needs to be determined by its \emph{measure}; 2) obscurities caused by the desire to preserve the  \emph{ideal character} of Geometry while the true Geometry, which reaches its [epistemic] goal, does not need such an ideality. \label{measure4}The whole of Mathematics is a science of measuring; everything that exists in nature is conditioned by being measured. As a consequence, differences between magnitudes  should relate to different kinds of measurements and to numbers which represent those magnitudes. All other concepts [are and] will always be obscure and insufficient. Thus I begin [the presentation of] Geometry with [i] the measurement of straight lines and [ii] the measurement of arcs with respect to their circles. Then I give the definition of angles, discuss the measurement of parts of the surface of a sphere with respect to the whole of this surface, and draw from that the definition of plane and solid angles. 

The second kind of obscurity [in geometry textbooks], namely, the ideality, occurs in the [definition of the] ratio of [straight] lines, in the length of a curve, and the magnitude of a curve surface. In order to avoid the obscurity of this sort it is sufficient to stick to the following principles.  \label{senses3}First, Mathematics aims at the actual measurement (\selectlanguage{Russian}действительное измерение\selectlanguage{English}) and for this reason it should never go beyond the requirements of our senses.  Second,  concepts, which are defined for certain objects, can be arbitrarily extended to other objects which are measured similarly. For example the magnitude of a curve with respect to a straight line cannot be understood along the pattern of the magnitude of [two] straight lines with respect to one another. Nevertheless nothing prevents us from understanding by the length of a curve the sum of those straight lines, which replace the parts of the curve, with the proviso according to which the smaller are the parts of the curve, the closer the sum of the straight lines approximates the length of this curve. This is what is [commonly] done in practice. Such a measurement is useful both theoretically and practically because it gradually approaches certain limit with the diminishment of the parts of the curve; this limit should be taken for the true length [of the given curve]. This definition [of length of a curve] is useful not because a curve admits for a concept of length but because the actual measurement [of a curve] should approach this limit as much as possible. The definition of such limits belongs to Analysis, which can take care of the ideality because it possesses [symbolic] means for representing ideal curves with equations. This is why the measurement of curve lines, curve surfaces and bodies bounded with curve surfaces should not be introduced in the Foundations of Geometry \label{foundations}; [In my teaching] I follow this rule and leave a treatment of this topic in all its generality to the Analytic Geometry. By the Foundations of Geometry I understand a part of Geometry, in which one cannot avoid considering geometrical magnitudes in the space. In this part of Geometry one proceeds from special cases to more general cases; Analysis proceeds in the opposite order beginning with general cases and only then considering special cases. At the present state of the art in Geometry the way in which one learns the foundations of this science is not a subject to any general rule;  probably it will always remain in this state in the future. The Foundations of Geometry extend until the point where one finds general rules for determining spatial positions and measuring triangles and pyramids; at this point the \label{ansyn} Synthesis ends and the Analytic Geometry further develops this science in an algebraic form.
$[\dots]$

B. Trigonometry

$[\dots]$,  \cite[208-209]{Modzalevsky:1948} 

Geometrical considerations are necessary in the Foundations of Trigonometry as far as they help to discover the specific property of trigonometric functions, namely, the value of the sum of two angles expressed in terms of sinuses and co-sinuses of those two angles separately.\footnote{\label{sinus} $sin(\alpha + \beta) = sin\alpha \cdot cos\beta + sin \beta \cdot cos\alpha$} But from that point onward Trigonometry becomes fully independent of Geometry and acquire all advantages of Analysis. \emph{Imaginary quantities} \label{imaginary} [in Analysis] which are so called because they comprise the square root of the negative unit, greatly help Trigonometry by shortening [some computations] and making solutions of certain problems easy and direct. If these quantities would be thrown away from Analysis because of the non-reality of their values then much [of the valuable content] would be lost. Indeed, the imaginary quantities lead to conclusions which are as valid as the possible [i.e., real] quantities. Thus we can see that the generality of Analysis extends not only to all possible magnitudes but also to imaginary ones. In this respect Analysis makes a contrast with [the synthetic] Geometry and all other parts of mathematics, which account [only] for real magnitudes (\selectlanguage{Russian}действительные величины\selectlanguage{English})\footnote{In the Russian traditional mathematical terminology there exist two interchangeable terms for real numbers none of which translates word ``real'' literally: (\selectlanguage{Russian}действительные числа \selectlanguage{English} (literally ``actual numbers'') and (\selectlanguage{Russian}вещественные числа \selectlanguage{English} (literally ``material numbers''). Lobachevsky uses here the former term.} of nature  and at every step require not only a [theoretical] possibility of measurement but also a specification of \emph{how} this measurement [can and] should be done.  Analysis is, so to speak, only a game of symbols, which does not limit itself to actual (that is, rational) values but  involves [symbolic] expressions, which do not admit for any [real] value. Using signs for non-extractable roots \footnote{Like $\sqrt{-1}$.}  \label{extend} Analysis  extends its operations to numbers, which strictly speaking do not exist. With the help of imaginary magnitudes Analysis expresses the dependence of one variable from another variable leaving it unknown how the former variable is transformed into the latter one.\footnote{Here Lobachevsky points to the capacity of Analysis to deal with real functions (like the trigonometric functions), which do not reduce to a finite number of elementary arithmetical operations.} This is achieved with signs which show [formal] properties of a given function.\footnote{This general point concerning the epistemic role of symbolic calculus can be demonstrated by the example of $sin$ function: rather than defining $sin$  geometrically one can define it formally as a function that satisfies certain algebraic equations. In this case the domain of a given function often can be extended. For example, the geometrical definition of $cosin(\alpha)$ as a ratio of two sides of a rectangular triangle that which form angle $\alpha$ implies that  the domain of the  $cosin$ function is a real interval $(0; \frac{\pi}{2})$, i.e., that $0 < \alpha < \frac{\pi}{2}$. But if the same function is introduced ``analytically'' via its formal properties then its domain can be extended to all real numbers and further to all complex numbers. An interesting fact that Lobachevsky does not discuss in this passage is that such formal extensions of traditional mathematical concepts may admit for new intuitive geometrical interpretations. In particular, the extensions of the domains of basic trigonometric functions to all real numbers is supported by interpreting their arguments as \emph{angles of rotation} (assuming, as usual, that clock-wise rotations are negative and counter-clock-wise rotations are positive). See \cite{Rodin:2010} for a further discussion.\footref{sinus}} When Analysis is granted more generality, the qualities of magnitudes, which are the subject-matter of this Analysis, should be treated with the same degree of generality.  Once imaginary magnitudes are admitted in the Analysis, using real values of [symbolic] expressions should be disallowed; one cannot reason about real values unless imaginary values disappear.\footnote{The passage is rather obscure but in my understanding Lobachevsky says here that the domains of functions should be carefully specified, and eventual changes of these domains should not remain unnoticed.}

This is my opinion about this branch of Analysis [i.e., the branch that involves imaginary magnitudes]. Sticking to this opinion I feel obliged to introduce my students to the imaginaries [i.e. to imaginary numbers and akin ``ideal'' mathematical objects] and to demonstrate to my students the fruitfulness of this source. There always exist, however, other means allowing one to reach the same conclusions but this time without using imaginaries \footnote{In other words, in Lobachevsky's view, imaginary quantities are dispensable. Indeed, a complex number can be represented with a pair of real numbers just like a fraction can be represented with a pair of integers. But such representation may dramatically increase computational costs and be unpractical otherwise.} Thus one should learn to use imaginaries because it is an easy instrument that facilitates [mathematical] discoveries. The formula for sinus of the sum of two arcs [see \footref{sinus}] implies many other formulae, which should be learned by students because they are very important in the Trigonometry. 

$\dots$        
\cite[211-212]{Modzalevsky:1948}

\label{centesimal} I prefer the decimal division of the circle because it is the most advantageous from a computational point of view. But since this division is not yet so common I also provide examples for the [traditional] 60-parts division.

C. Analytic Geometry

This part of Pure Mathematics in my teaching [program] follows Trigonometry and draws on Trigonometry. The Analytic Geometry is perfect as far as its systematic organisation and the clarity of its methods are concerned. This is why teachers should follow the best writers on this subject such as [Gaspard] Monge who made a significant contribution to the dissemination of this discipline.  \footnote{Most likely Lobachevsky refers here to Gaspard Monge's textbook on Analytic Geometry \cite{Monge:1809}, which is based on Monge's lectures on this subject in the \emph{\'Ecole Polyt\'echnique}. The first edition of this textbook under a different title appeared in 1795.} The reason of this perfectness is that the Analytic Geometry comprises analytic operations, which aim at the determination of geometrical magnitudes. But the scope of Analytic Geometry cannot be  determined [once and for all] because this discipline expands its scope along with Analysis [itself] \ref{extend} . As it has been already stressed, teaching Analytic Geometry is unproblematic. But before this discipline becomes fully analytic one needs to make a transition from the [traditional language of] Foundations of Geometry to algebraic symbols and equations, which cover all possible [geometric] cases and thus make the geometrical representation unnecessary. Then one should establish rules allowing to make a backward transition from numbers obtained in the Analysis [back] to actual geometrical magnitudes. The involvement of Synthesis in this procedure sometimes causes difficulties. Those who wrote on that subject failed to attain the necessary rigour.  But it is nothing but a matter of attention: one should not forget to interpret analytic operations [in geometrical terms] when these operations apply to geometrical magnitudes.\footnote{It is hard to say whether this Lobachevsky's critical remark targets any particular author. In any event it is clear that Lobachevsky criticises in this passage the way in which some authors of his contemporary textbooks combine the traditional ``synthetic'' foundations of geometry with the new analytic apparatus. In Lobachevsky's view this problem can and should be fixed by the rebuilding of geometrical foundations rather than by adapting the analytic apparatus to the existing geometrical foundations.} Indeed, a simple notice would suffice to clarify the notion of positive and negative lines and angles: adjoining this notion to Foundations of Geometry gives rise to Analytic Geometry. After that it is easy to demonstrate the following important truth: when one applies geometrical considerations and uses diagrams in order to facilitate a [geometrical] representation, the obtained algebraic equations often turn out to be perfectly general in spite of the fact that a diagram may possibly represent only one special case. 

$\dots$    

D. Differential Calculus
$\dots$     \cite[p.215]{Modzalevsky:1948}

\label{increments1826}Recently significant attempts have been made by mathematicians to attain in the foundations of Differential Calculus the [degree of] rigour, which is appropriate in Analysis. It should be admitted, however, that these attempts are irrelevant when Analysis is applied in Geometry and in Mechanics. In these cases one needs to see differentials from a different point of view according to which the differentials are finite increments, which should  [under appropriate conditions] be discarded [in computations]. Consider a curve line represented in [a system of] coordinates and compare it with the [corresponding] curve in [synthetic] Geometry. Suppose that we want to measure the length of the latter [synthetic] line. It is easy to see that while the former [analytic] line is continuous, the latter [synthetic] line can be [properly] understood only as a composition of extremely small straight lines. A continuously accelerated motion can be understood only as a motion, which is uniform during very small time intervals. Since the real measurement of geometrical and mechanical magnitudes cannot be free of the representation of increments, Analytics should not try to avoid it either. On the contrary, this [finitary] approach to the concept of Differential Calculus agrees with the aims of Geometry and Mechanics and provides more clearness in teaching. While abstract numbers of Analysis require to make a series of [logical] conclusions, the representation of mechanical and geometrical magnitudes immediately brings to the imagination a whole nexus of objects through a single glance on these objects. Those who want to satisfy themselves with numbers neglect the advantages associated with the essential qualities of natural magnitudes.\footnote{See \ref{increments1825} above. Lobachevsky's idea of replacing differentials with finite increments is further discussed by him in his \emph{Algebra} of 1834, see \ref{apG}.}

I follow here far-reaching visions of great mathematicians and believe that I have grasped their ideas properly. But thanks to my long teaching practice I learned to distinguish between teaching advantages and the perfection of science.

$\dots$     \cite[217-218]{Modzalevsky:1948} 

Teaching Mathematical Physics in 1825/26. General overview. 

 $\dots$
\label{physmath3} 
 Physics is everywhere supported by Mathematics and owes to Mathematics its very existence. But the major applications of Pure Mathematics and Mechanics in Physics concerns so far the parts of Physics that treat weightless bodies such as the matter of heat, light, electricity and magnetism. The theory of equilibrium and motion of liquids and the theory of sound are two cases of mechanical theories, which deserve a spacial attention because of their importance. These theories qualify as physical because they involve a comparison of experimental measurements and [theoretical] computations; in these theories multiple observations help to discover properties of bodies, which, in their turn, help to interpret all other phenomena and represent them rigorously in numbers and magnitudes. In this discipline we no longer need the conceptual generality, which is characteristic for  Analytics and constitutes its great advantage. From all [possible] speculations about forces we take here only those, which are actually observed\footnote{i.e., which agree with what is actually observed} in nature, and which belong to the narrow domain of things that are at hand [i.e., are easily accessible for and can be controlled by humans].  \label{physast1}Other speculations about forces belong either [i] to Astronomy,  or [ii] to Chemistry, Mineralogy or Natural History where forces and their actions so far could not be computed, and where the human mind does not go beyond observations of laws of [chemical] affinity and laws of organism but merely lists  natural phenomena without being able to discover their first causes. These disciplines (of category [ii]) should be contrasted to Mechanics which derives all cases of equilibrium and motion from the concept of [mechanical] force. Separating from the Natural Science Astronomy, Chemistry and Natural History, we are left with the discipline called Physics.\footnote{\label{ptolemy} In this passage Lobachevsky makes a double contrast: 
 
 A)  between Mathematical Physics (which by the 1820s already included not only Classical Newtonian Mechanics but also Optics, theory of Electricity, Acoustics and some other disciplines), on the one hand, and various non-mathematised natural sciences such as Chemistry and Natural History (i.e., Geology and Biology), on the other hand. 
 
 B) between Mathematical Physics and Astronomy, which is mathematised like Mathematical Physics but nevertheless, in Lobachevsky's view,  does not fall under its scope. Mechanics is not qualified by Lobachevsky as a part of Physics either: Lobachevsky treats is rather as a part of Applied Mathematics that does not belong to the Mathematical Physics.  
 
In the end of the passage Lobachevsky identifies Physics with Mathematical Physics, so Astronomy, Chemistry, and Natural History do not qualify, in Lobachevsky's view, as parts of Physics.  

The concepts of ``organism'' and of chemical ``affinity'' were widely used in Lobachevsky's contemporary natural science but unlike the concept of ``force'' in Classical Mechanics they were not quantified, i.e., not measured, not represented with numbers and not theoretically computed \cite{Hornix:1988}. 

\label{ptolemy} Lobachevsky's treatment of Astronomy and Mechanics that separates these disciplines from Physics draws on a long tradition that dates back to classical works of Ptolemy, Archimedes and Aristotle. Book 3 of Newton's \emph{Mathematical Principles of Natural Philosophy} titled \emph{The System of the World}  treats astronomical phenomena on equal footing with human-scale mechanical motions which, by Lobachevsky's word, are found ``at hand''.  The abandoning of the traditional Aristotelian distinction between the ``sublunar'' and celestial  mechanical phenomena and treating these two kinds of phenomena on equal footing in Newton's \emph{Prinicples} was characteristic for the nascent Modern Science and determined the great success of the Newtonian Mechanics. So Lobachevsky's separate treatment of Mathematics Physics and Astronomy in the 1820s can be qualified as an anachronism. But in view of Lobachevsky's geometrical ideas this anachronism can be also seen as a preconception of later Einstein's discoveries who has shown that the Newtonian assumption according to which Mechanics apply uniformly at all scales is, after all, erroneous.} 

\label{history} In Physics [students] too easily abandon the difficult path which led the human mind to its [scientific] discoveries; too easily they accept the conceptions of bodies and forces established by the preliminary program of the Gymnasium, which are not always accurate; too quickly they assume the putative first causes of phenomena and want to derive all their actions. They imagine that Physics has already achieved its perfect state and can serve as a model for all natural sciences. In the reality, however, only Astronomy and some parts of Physics are presently found in this perfect state. Some other physical theories have only a [fancy] mathematical appearance but still lack firm [mathematical] foundations. This is a very significant flaw [of current physical theories] that the teacher must stress.  He must tell the truth without diminishing the [epistemic] value of mathematics. He should help students to understand the benefits of [applied] mathematics in spite of the fact that only with a great labor and after long continued efforts we presently manage to comprehend a few secrets of nature, which were earlier so deeply hidden in it.

\section{From public speech  \emph{On the most important aims of education} (1828)}\label{apD}

 \cite[p. 323-324]{Modzalevsky:1948}
 
\label{language1} It should be admitted that our [human] superiority over other animals is due to the gift of speech rather than to our mind. As anatomists tell us, even those animals which are otherwise the most close to humans lack organs that might allow them to pronounce complex sounds.\footnote{Recent findings seem to support this Lobachevsky's claim, see \cite{Nishimura:2022}.}  They are not allowed to communicate concepts to one another. Only humans are granted this right; only humans use this gift; only humans are ordered to learn, to master their minds, to seek truths collectively. Words disseminate and diffuse the light of learning like beams of their minds. The language of a people is an evidence of its education and a reliable proof of the degree of its enlightenment. 
 
 What is the reason of recent brilliant successes of mathematical and physical sciences, which are the glory of our age and the triumph of the human mind? Without any doubt these successes are due to their artificial language. For how else all those symbols of various calculi can be called if not a special and very concise language, which expresses extensive concepts with a single stroke without tiring unnecessary our attention. These achievements of mathematical sciences, which exceed achievements of all other sciences, rightly surprise us. They force us to admit that the human mind alone is capable to acquire truths of this kind, and that it probably is looking for truths of different kind in vain. 
 
 One should also admit that mathematicians discovered direct means for acquiring knowledge. We use these means since recently. They were indicated to us by the great [Francis] Bacon.  He told us that we should stop labouring in vain trying to extract the wisdom from the mind alone. He urged us instead to ask nature, which keeps all truths and which is ready to answer accurately and satisfactorily. Finally, [Ren\'e] Decartes' genius produced this happy change [in science and education], so thanks to his deeds we find ourselves in the age when the shadow of ancient scholasticism is hardly visible in the Universities. In this [Kazan] University the youth will not listen to empty words without any thought or sounds without any meaning. Here one is taught about what really exists but not about what was invented by an idle mind. Here we teach exact and natural sciences with the help of languages and of historical knowledge. Here professors teach disciplines which they study their whole lives since their youths when they first feel in themselves the desire [to study] and some gifts. What a pity that some [parents] still prefer the apparent advantages of home education to the true enlightenment [in a state-sponsored gymnasium or university]. One who wants to educate their children for the State should use means that only the State can provide. [In other words], they should educate their children in public institutions.

 \section{From \emph{Instructions for teachers of mathematics and physics in gymnasia} composed in 1830 in Lobachevsky's capacity as a Rector of Kazan University}\label{apE}
 
 \cite[p. 526-528]{Lobachevsky:1976}

 \label{language2} For [operating with] abstract and general concepts of magnitude, and for connecting magnitudes to each other people invented signs. Just like the gift of speech enriches us with other people's opinions, the language of mathematical signs serves as mean [of communication] which is more perfect [than natural languages], more concise and clear, and which allows one to communicate to another person one's concepts that they have conceived, their truth that they have comprehended, and the dependency between magnitudes that they have discovered. And just like an opinion may appear false when one interprets the meaning words not as it has been intended, every mathematical judgement stops when one no longer understands what the signs are signs of.  This is why it is necessary for a teacher to give perfectly determined and rigorous concepts along with [explanations of] using signs; on the top of that a teacher should provide examples which not only demonstrate rules [of using signs] but also prevent a mechanical use of those signs. 
 
 What has been said so far is sufficient for seeing how one should teach those part of mathematics, which are called \emph{analytic}. But in the very beginning [of teaching mathematics] there occur difficulties of a special sort, which should be well understood and overcame. These difficulties concern the transition of our mind from objects that immediately affect our senses to numbers to letters, which represent these numbers in a general manner. 
 
 \label{senses4} Those first concepts, which we obtain from the nature directly through our senses, serve mathematical sciences. Even our first judgements about objects falling under these concepts, are produced by a habit in our senses but not by an action of our mind, which embraces all possible cases from a general point of view. But the arithmetical calculi and their rules are invented with the help of such a general view on magnitudes, which are determined by measuring. This helps us to solve [arithmetical] problems in their full generality. 
 
A beginner is not able to understand where the general rules of Arithmetic come from. It is necessary for a beginner to replace [mathematical] judgements by sensual impresions, so they could proceed independently from these immediate impressions to the sphere of abstract concepts where the mind begins its operations. 
 
   $\dots$
   
  \cite[p. 528-529]{Lobachevsky:1976}
 
 In gymnasia teaching mathematics after Arithmetic is divided into teaching Algebra and teaching Geometry.  \label{measure5} Algebra assumes that the measurement has been already done, so all [geometrical] magnitudes are represented with numbers, and the numbers are denoted by letters for purposes of general reasoning. In Foundations of Geometry [on the contrary] we talk about the measurement itself. This is why the content of this part of mathematics is independent of Algebra and should be taught separately. It follows that the method of teaching purely geometrical contents should differ from the method of teaching Algebra up to the point where Geometry joins Algebra, which is called the \emph{application} of Algebra to Geometry. 
 
\label{parallels2}Major difficulties of Geometry are related to the first concepts of geometrical magnitudes and in geometrical constructions, which require the help of imagination.  \label{senses5}The initiation to Geometry should involve only those concepts which are obtained immediately from senses without further studies and without any external support. Such concepts are simple, and truths based on these concepts are sensual. These truths comprise some obscurity and vagueness but [this is acceptable because] a perfect rigour might involve [further] studies, which are inappropriate in gymnasium because they reveal to us the impossibility to learn all geometrical properties of bodies sufficiently.  Such further studies [into the Foundations of Geometry] also show that in order to determine a domain where the usual geometrical propositions are justified, one needs to use astronomical observations and apply other parts of mathematics.\footnote{It is clear that Lobachevsky refers here to the Problem of Parallels without calling it by its name. In his opinion whether or not the Axiom of Parallels is true, i.e., the geometry of the physical space is Euclidean or not, can be decided only on the basis of astronomical data. But Lobachevsky does not recommend to include such considerations in the courses of Elementary Geometry taught in gymnasia, and suggests teachers to accept this Axiom on an intuitive basis even if this basis is not sufficiently firm.} 

\label{measure6}\label{parallels3}The main goal of teaching Geometry is to give [pupils] general rules of measurement. These rules can be obtained only via a geometrical reasoning. They start with [the rules of measurement of] straight lines with straight lines, and the measurement of arcs of a circle with the circle. Then follows the doctrine of angles and of relative positioning of straight lines. Here of the major importance is the doctrine of parallels, which [at the given stage of mathematical education] cannot be yet rigorous but should be convincing on a sensual basis and by the simplicity of the assumption [i.e., the Axiom of Parallels] in spite of the fact that this assumption is arbitrary.

 \cite[p. 534-536]{Lobachevsky:1976} 
 \label{physmath4}
 Pure Mathematics takes care only of measuring the [physical] space assuming that all other measurements have been done earlier. It also assumes that every magnitude is expressed with a number. This is why Analytics is the main part of the Pure Mathematics. General assumptions concerning forces and their measurements (both theoretical and practical) belong to other mathematical sciences which are also known under the name of \emph{mixed mathematics} \footnote{ ``Mixed mathematics'' (Russian \selectlanguage{Russian}смешанная математика\selectlanguage{English}) is an outdated term, which by an large is equivalent to today's term  ``applied mathematics''  Lobachevsky might borrow this term from d'Alembert or another contemporary source. For the history of the term see \cite{Brown:1991}.} and which includes two major parts, namely, Mechanics and Physics. The latter discipline teaches us about natural forces, the laws of their actions, how they can be measured and which changes and phenomena they happen to cause. 
 
 Knowledge of basic Mechanics is necessary for studying Physics. Since in gymnasia Mechanics is not taught as a separate subject, the two disciplines [i.e., Physics and Mechanics] should be combined; [the program of such a combined course] should include first principles concerning forces and general rules of equilibria and motions. In  order to ensure the continuity of teaching, Mechanics, being an introductory discipline,  should meet Physics at certain point. Indeed, the general [geometrical and mechanical?]  properties of bodies, which include their capacity of motion, and which reveal to us the nature of things and constitute the [necessary] conditions of existence of those bodies, allow us to account with an appropriate degree of precision for laws and rules which govern [mechanical] forces and changes of bodies' positions after the disturbance of their equilibria. 
 
 These general properties of bodies can be seen as immediate consequences of hypotheses concerning their [material] composition including the hypothesis concerning the mode of existence of all things in the world as this world is seen by humans. Those hypotheses should serve as a foundation and a sufficient source for explaining away and proving possible actions  of hidden forces, which depend on the distribution of matter in the space. This foundation may not be sufficient, however, for computing details of all [physical] phenomena. As far as we satisfy ourselves with a mere interpretation of possible phenomena without being able to derive the very laws of motion and understand forces in their elementary actions, mathematics cannot be applied [in such a discipline].  Examples of such disciplines include the Natural  History and Chemistry where we merely describe phenomena, limit our knowledge to [what one may learn from] observations, and help ourselves by inventing a system. But in these sciences we are unable to follow the Nature intellectually and suggest a foundation from which such a system could be produced with the power of our judgement. 
 
 Thus Physics covers only those natural phenomena, which are sufficiently well known to us, so we are able to apply mathematics here. Such phenomena are usually multiple, i.e., are observed in big numbers rather than [only rarely] one by one. In such cases we are able to understand causes and discover the laws of acting forces because the bodies [associated to the observed phenomena] are multiple and do not avoid our senses. Forces through which those bodies interact in such cases can be observed and their dependence on the positions of the bodies can be known. After the fortunate discoveries of Newton it was universally recognised that all matters should be imagined as a collection of material points or centres of sort from which issue forces that alone cause all natural phenomena and all natural motion. Given these assumptions it is easy to see that the more narrow is the domain where act the forces, the more complex become their actions. This is why the sphere of Physics has been [first] discovered when scientists learned about forces which govern motions of the celestial bodies, which are placed at large distances [from human observers and from each other] in the almost empty space; since the sizes of celestial bodies are small comparatively to these distances, they act onto each other like points, and represent the laws of motion in all their simplicity. But because of its large scope, its importance, and its special goal this discipline is separated from Physics and constitutes a special science known under the name of Astronomy.\footnote{As the reader can see Lobachevsky's account of Physics and Astronomy in this document is somewhat ambiguous: on the one hand, he recognises Astronomy as a branch of Physics but, on the other hand, insists that the two disciplines are separate. Lobachevsky's historical remark that relates Astronomy to the emergence of Physics (or more precisely \emph{Mathematical} Physics) needs some comments.  Astronomy has been first established as a mathematical science quite independently of Physics and without using the concept of force \footref{ptolemy}. Premodern Physics typically did not draw on Astronomy either. But the birth of Modern  Physics, which culminated with Newton's achievements indeed systematically relied on the contemporary Astronomy and conceived of astronomical phenomena as physical, i.e., as a subject to universal physical laws (Newton's Laws). See also \ref{physast1}. As judged from today's anachronistic perspective Lobachevsky's emphasis on the simplicity of standard astronomical models which include only point-like massive bodies is quite justified. As it has been realised by the end of the 19th century, the complexity of such a mechanical system may explode already in the case when the number of such interacting point-like bodies exceeds two \cite{Barrow-Green:1997}. So in order to be more precise by today's standard Lobachevsky had also to limit the number of bodies in his argument by two.}
 
 \label{touch2} On the contrary, the most difficult is the case of forces acting at very short distances, that is, the case of touching bodies. Nevertheless the application of general laws of motion and equilibrium in this case also helped to make some [useful] conclusions and explain certain specific phenomena. Such phenomena which have been called \emph{capillary} also belong to Physics. 
 
After describing the subject-matter of Physics and its differences from other exact sciences it is appropriate to notice that teaching Physics has two different aims: [i] to discover causes and use them for explaining phenomena and [ii] to teach how to measure forces and their actions. The recent most important achievements of Physics and the present perfect state of this science are due to the satisfaction of the latter requirement (i.e., [ii]). Since this part of Physics  (i.e., [ii]) needs mathematical skills and thus should be taught in universities. But the first principles of Mechanics and the [physical] explanation of the most important natural phenomena (i.e., [i]) can be grasped with the [mere] common sense and be confirmed with experiments. It should be taught in gymnasia and belong to all the Russian educated youth.

\section{Introduction to the \emph{On Foundations of Geometry} (1829)}\label{apF}

\cite[vol.1, p. 185-190]{Lobachevsky:1946-51}

It appears that difficulties increase when one approaches fundamental concepts of nature just like they increase [when one moves] in the opposite direction where the [human] mind aims at acquiring new knowledge. This is why [all] difficulties in Geometry should belong primarily to this very discipline.\footnote{In this sentence Lobachevsky implicitly argues against a popular view according to which foundational difficulties in Geometry should be resolved in Philosophy. See \ref{sum9}.} [Mathematical] means that are needed for attaining the ultimate rigour in this case can hardly suit the purpose and the simplicity of [teaching] this discipline. Those [mathematicians] who tried to meet this requirement put on themselves such severe restrictions that all their efforts could hardly be rewarded. Finally, let us add that since Newton and Descartes all of Mathematics became Analytics and progressed very rapidly leaving  behind the [Euclidean] doctrine, which it no longer needed. Since then this doctrine attracted little attention [of the mathematical community]. This is why Euclid's \emph{Elements} preserved [until today] all their original flaws  in spite of their very old age and in spite of our brilliant successes in Mathematics. 

Indeed, one should admit that no mathematical discipline should begin from such obscure notion, from which we commonly begin [our presentations of] Geometry repeating Euclid. Nowhere in Mathematics one can tolerate such lack of rigour as one that we have to tolerate in the [Euclidean] theory of parallel lines.  It is true that the representation of [geometrical] objects in our imagination protects us from false conclusions and obscurities in the first general notions of Geometry. \label{astro2} It is also true that the validity of geometrical truths admitted without a proof is justified by their simplicity and by our experience including astronomical observations. Nevertheless this cannot satisfy a mind which is familiar with rigorous reasoning. We cannot neglect this problem of which we so far have no solution, moreover that this solution may help to solve some other problems.

$\dots$
\label{innate}
First concepts  from which begins a science should be clear and reduced to the least number [of concepts]. Only then they can serve as a firm and sufficient foundation of a science. The opinion according to which there exist innate concepts should not be trusted.\footnote{In the 18th century the doctrine of innate ideas was attributed to Descartes and Leibniz and severely criticised by John Lock, \'Etienne Candillac and their followers \cite{Samet:2019}. Lobachevsky aligns here with this critique.}  Nothing can be simpler than the concept [of counting numbers] which lays in the foundation of Arithmetic. We easily learn that everything in the nature can be measured and counted. But the [basic] propositions of Mechanics are different: one cannot obtain them from one's everyday experiences. The eternity and constancy of [inertial] motion measured by its [constant] velocity and the invariance of mass and other such like truths [could not be grasped immediately but] required for their comprehension a [lot of] time and the assistance of some other knowledge, so the discovery of this principle had to wait for a genius. \footnote{It is not wholly clear what Lobachevsky calls here ``inner senses''. But it is remarkable that Lobachevsky's view expressed in this passage diverges from the straightforward sensualism, to which he apparently adheres elsewhere, see \ref{senses1},\ref{senses2},\ref{senses3},\ref{senses4},\ref{senses5}. He quite rightly stresses the fact the principle of inertial motion aka the First Newton's Law is not supported by the everyday experience of observing and operating with moving objects. In oder to establish this principle and realise its fundamental role in Mechanics Galileo and other early modern scientists had to combine theoretical considerations with experimental and observational practices in very sophisticated ways. Even if Lobachevsky does not say this explicitly the context suggests that his following attempt to rebuild the Foundations of Geometry follows the pattern of Mechanics rather than that of Arithmetics where, in Lobachevsky's view, the relevant fundamental concepts including the concept of counting number are obtained ``directly through senses'', and don't require any further theoretical work. }

Among the properties, which are common for all bodies, only one deserves to be called [properly] \emph{geometrical}, namely, the touch.\footnote{See \ref{sum4} for a commentary}. It is impossible to fully communicate what we mean by touching using only words: this concept is obtained through senses, mostly through seeing, and for this reason can be understood only with these senses. Touch is a specific characteristic of bodies: we don't find it in forces, times or elsewhere in the nature. When, given a body, one abstracts away all its properties except touch, the body is called \emph{geometrical}. 

A touch of two bodies combines them into one body. All bodies are parts of the same body called the \emph{space}. A body [$A$] is [called] \emph{bounded}  when it is touched by another body [$B$]  (called the \emph{ambient} body  [of body $A$]) in such a way that  [$A$] cannot touch any other [third] body. If two bodies [$A$ and $B$] together constitute the whole space then [$B$] is called the \emph{ambient space} [of body $A$]. The void occupied by a body inside the space is called a \emph{place}. Two bodies  [$A$ and $B$]  are [called] \emph{congruent} when they can \emph{fill the same place} (i.e., to complement the same space) without undergoing any change. Two bodies [$A$ and $B$]  are only  \emph{equal} when, given body [$A$] fills certain place, the other body [$B$] can fill the same place if it is divided into [a finite number of] parts and then those parts are reassembled in a different order.\footnote{In this passage the adjective ``congruent'' translates Lobachevsky's Russian term \selectlanguage{Russian}одинаковый\selectlanguage{English}, which may have different meaning both in and outside mathematics. The adjective ``equal'' in the following sentence translates its standard Russian equivalent  \selectlanguage{Russian}равный \selectlanguage{English} and requires a historical comment. Following Euclid, Lobachevsky calls two given solids  $A$ and $B$ equal when they are \emph{equidecomposable}, i.e., when  $A$ can be cut into a finite number of pieces such that by reassembling these pieces one can build a new solid congruent to $B$.  Lobachevsky was, of course, not aware about the fact that polyhedra having the same volume are not always equidecomposable. This was first shown by Max Dehn in 1900 \cite{Boltianskii:1978}.}

An imaginary division of a given body into two parts is called a \emph{section}. Each of the two parts of a section is called a \emph{side} of this section.

We learn geometrical properties of bodies by dividing them into parts. These properties constitute the foundation of Geometry, which is the following:
\begin{enumerate}

\item Any body can be divided into parts [$P_{1}, P_{2}, \dots P_{n}$] such that a given part does not touch the part situated after the next part [i.e., such that for all $1 \leq i \leq n$, $P_{i}$ and $P_{i+2}$ do not touch each other]. We shall call such sections \emph{translational}  (\selectlanguage{Russian}поступательный \selectlanguage{English}). The number of such parts [of a given body] is unlimited.\footnote{Think about slicing a sausage.} 

\item  Any body can be divided into parts  [$P_{1}, P_{2}, \dots P_{n}$], which all [pairwise] touch each other in such a way that every new section brings about two new parts of the given body. We shall call such sections \emph{rotational}  (\selectlanguage{Russian}обращательный \selectlanguage{English}). The number of such parts [of a given body] is unlimited. \footnote{Think about dividing a round pie into slices by cutting it with a long knife through its centre. Then the slices come in even numbers, and every new cut increases the number of slices by two.}

 \item   \label{3D3} Any body can be divided with 3 sections into 8 parts, which all [pairwise] touch each other. But the further doubling of the number of parts with a new section is impossible. Such sections we shall call the \emph{three principal} sections.\footnote{Think about cutting a watermelon into 8 parts by three consecutive plane sections, which  pairwise orthogonal  and share a point located inside the watermelon.  Each such section doubles the number of parts (starting with 1 in the case of the whole watermelon) but in the 3-dimensional space the process stops after the third step.} 

\end{enumerate}        

Rotational sections can also increase the number of parts by one. This happens when a new part is added through division of an earlier obtained part. Thus the specific feature of rotational sections does not concern the number of parts but rather concerns their mutual touch. This, however, is not sufficient [for defining the concepts of translational and rotational sections appropriately]. It is necessary [to add the condition according to which] translational sections made on the top of two [series of] rotational sections applied to the opposite side of a given body produce parts which belong to different sides and do not touch each other.    
\footnote{Think about a pie cut in two slices and then of one of those slices cut into two thinner slices. As for the additional conditions suggested here by Lobachevsky, I'm not sure that I fully understand its intended meaning. But here is a possible interpretation. Cut a watermelon into two halves. The obtained semi-spheres are the two ``opposite sides''. Now slice each of the two semi-sphere in the pie-like manner (i.e., ``rotationally'') and finally chop each slice ``translationally'' into a number of pieces. Now Lobachevsky's condition is tantamount to saying that given two slices belonging to different semi-spheres there are always exist some slices that do not touch each other. Notice that slicing without further chopping does not always meet this condition.      

This whole passage shows the degree to which the text of Lobachevsky's \emph{Foundations} has been poorly edited. In the preceding paragraph Lobachevsky offers us what looks like his working definitions of ``translational touch'' and ``rotational touch''. But in the following paragraph he already revises these definitions.}

A similar remark applies to the three principal sections. Any two these sections qualify as rotational, and this is a sufficient condition for being principal. 

 Translational sections made on the top of the three principal section determine in a given body its \emph{three dimensions} (\selectlanguage{Russian}три протяжения \selectlanguage{English}).\footnote{Russian \selectlanguage{Russian}протяжение \selectlanguage{English} translates directly as English \emph{extension}. It is clear that Lobachevsky does not distinguish here between topological and metrical concepts as we do this today, and leans to a \emph{metric} concept of spatial dimension. This agrees with Lobachevsky's emphasis of the role of measurement in Geometry and science in general, see \ref{measure4}.}
  
\label{apFM} To \emph{measure} a body means to count the congruent parts into which this body is divided in three dimensions using translational sections along with another body taken for the unit. 

A concatenation of two bodies into one qualifies as a section of the latter [composed] body. A touch [of two or more parts of a given body] may belong either to a single [section], or to a number of rotational [sections], or to three principal sections. Correspondingly, this touch qualifies either as \emph{surficial}\footnote{In Russian original ``\selectlanguage{Russian}поверхностный \selectlanguage{English}''. This regular Russian adjective is usually translated by English adjective \emph{superficial}. In order to avoid irrelevant non-mathematical connotations I take the liberty to translate it in the present mathematical context by term  \emph{surficial} used in Geology.}, or as  \emph{linear}, or as  \emph{point-like}. Imagine a body divided into eight parts, which are obtained through three principal sections (and which could not be obtained otherwise). The first of those three sections produces a surficial touch of the two parts.  The second section produces [among other things] two parts [in fact, two pairs of such parts] which belong to the opposite sides of both sections [i.e., of the first and of the second section] and which touch linearly. The third section brings about two parts [in fact, two pairs of such parts] which are situated on the opposite parts of all the three sections, and which touch in [a single] point.\footnote{Compare \ref{3D1}}. 

A body is called a \emph{surface} when it touches another body surficially, and when one takes into consideration only the touch of the two bodies disregarding all parts of each of those bodies which do not touch the other body. So one of the three dimension is eliminated, and by leaving out the irrelevant parts of the surface one attains the thinness of a sheet of paper or as far as can reach one's imagination. The two mutually touching bodies in this case are [identical to] the two sides of the surface. 

A \emph{line} is a body that linearly touches another body provided that it is allowed to disregards parts which do not touch this other body. Thus one attains the thinness of a hear, of a line written by a pen on paper, etc. By turning a body into a line two dimensions are eliminated. This happens because a line is formed by two [principal] sections, and additional translational sections leave out the irrelevant extra parts [of the original 3-dimensional body]. 

A body is called a \emph{point} when one considers its point-like touch to another body disregarding all parts of the former body which do not touch the latter body. Thus one can reach the size of a grain of sand or that of a point produced by a pen on paper. 

A point is formed by three principal sections in the space provided that the irrelevant extra parts are left out by translational sections. This is why a point has not extensions [i.e., has dimension zero]. 

What matters in a surface, in a line or in a point is only the touch between two bodies [$A$ and $B$]. This means that all changes in one of these bodies [say, $A$] are irrelevant as far as they do not prevent parts of the other body [$B$] from touching [$A$] and do not produce new parts [of $A$] which touch the other body [$B$]. This is why it is allowed to replace all translational sections used for declaring extensions [i.e., for measurement] by their rotational counterparts. It follows that a line does not change the magnitude of a surface, and a point does not change the magnitude of a line.\footnote{In this passage Lobachevsky points to the possibility of changing a coordinate system in geometrical measurements. Consider on plane $P$ a Cartesian (orthonormal) coordinate system $S$ with axes $OX, OY$ where $O$ is the origin of $S$. The ``translational sections '' of axes $OX, OY$ determine distances from $O$ (in terms of some chosen unit of length). $S$ determines on $P$ a square grid which allows for measuring plane figures on $P$. Then consider another such coordinate system $S'$ with axes $OX', OY'$ sharing with $S$ the origin $O$, which can be obtained from $S$ by a rotation with center $O$. The measurements of areas in $P$ made with a new grid associated with coordinates $S'$ give same results: areas on $P$ do not depend on the choice of a coordinate system. Similarly, the measurement of distances on a given line $OX$ does not depend either of the choice of origin $O$ or of the choice of the positive direction on  $OX$. 

Let us now see how Lobachevsky describes this construction in terms of his theory of touches and sections. According to this theory plane $P$ emerges as a ``touch'' of two bodies $A,B$. Alternatively $P$ can be described as a section of the compound body $A\cup B$.  Lines $OX, OY$ and measuring marks on these lines emerge as results of further sections made on $P$ \textendash\ and hence also on $A\cup B$. These further sections, however, do not affect $P$ itself: compound body $A\cup B$ is still divided by plane $P$ into two ``opposite sides''  $A,B$. The  ``replacement of translational  sections used for declaring extensions by their rotational counterparts'', as far as I can decipher this Lobachevsky's expression, amounts to considering on $P$ a series of  ``rotational sections'' with center $O$ and choosing in this series a new pair of axes $OX', OY'$. It is less clear how to think of rotational sections made on a line but most probably in this case Lobachevsky means a central symmetry (which changes only the chosen ``positive'' direction without changing the origin). Lobachevsky does not say either that the involved surfaces are plane, or that the lines are straight.  But apparently this is the case that he has in his mind.}   It also follows that a line should belong to a whole system of  rotational sections, and that a line [$L$] is determined by two surfaces [$S_{1}, S_{2}$]. In this case the two surfaces are said to \emph{intersect} in this line [: $S_{1} \cap S_{2} = L$].  Each of those surfaces  [$S_{1}, S_{2}$] is divided into two parts [$S^{+}_{i}, S^{-}_{i}$], which are the two sides of the line [$L$]. 

A point [$O$] belongs not only to the three principle sections $OXY, OXZ, OYZ$ but also to all their associated rotational sections [i.e., to all surfaces obtained by rotation of surfaces  $OXY, OXZ, OYZ$ around axes $OX, OY, OZ$]. This is why a point [$O$] is determined by two lines [$L_{1}, L_{2}$], which are said to  \emph{intersect} in this point  [: $L_{1} \cap L_{2} = O$]. Each line [$OX$] is divided by point [$O$] into two parts  [$OX{+}, OX{-}$], which serve for assigning (\selectlanguage{Russian}назначение \selectlanguage{English}) two sides to the given point  [$O$].

A point has no magnitude and no extension, and thus does not admit for measurement.

\label{apFD} If [i] two bodies $A, B$ touch a third body $C$ in points [$P_{A}, P_{B}$ correspondingly], and [ii] $A, B$ are connected through certain body $D$, which does not touch body $C$, then the relative position of [$A$ and $B$] is determined, which is also called the \emph{distance} between bodies $A$ and $B$. The distance so determined is not affected by eventual changes in $A, B, D$ related to taking out from or adding some new parts to these bodies provided these new parts do not touch body $C$. [More generally,] the  distance is not affected by changes in  $A, B$, which do not affect the touch of these bodies to body $C$. Thus a compass serves for assigning distances [between points].\footnote{In this last paragraph Lobachevsky uses his Dimension Theory for grounding the geometric concept of measure. See \ref{sum4} above for a commentary.}

\section{Introduction to \emph{Algebra or Calculus of Finites} (1934)}\label{apG}
 \cite[vol.4, p. 23-27]{Lobachevsky:1946-51}
 
 According to Alexander of Aphrodisias \emph{analysis} ascends to the first principles while  \emph{synthesis} proceeds from the first principles to conclusions. This is how the Ancients wanted to understand analytic and synthetic methods in Geometry.\footnote{For Alexander's account of analysis and synthesis, and for further references to Greek theories of analysis and synthesis in Geometry, see \cite{Haas:2024}. Where Lobachevsky obtained this reference to Alexander is an interesting open question.} In today's Analytic Mathematics one similarly defines unknown numbers as if they were first principles  (\selectlanguage{Russian} как бы начала \selectlanguage{English}) from which one composes equations by adding to them certain data. Every computations aims at finding something unknown. This is why the rules of [all] computations constitute a single doctrine  (\selectlanguage{Russian} учение \selectlanguage{English}), namely \emph{Analytics}. Analytics can be divided into  \emph{Arithmetic}, \emph{Algebra} and \emph{Differential Calculus}.  In Arithmetic one begins with numerical examples; in Algebra one proceeds by using letters instead of numbers, which allows one to consider arbitrary numbers and preserve the conceptual continuity [with Arithmetic]. Yet, in Algebra one avoids infinitesimal methods and limits.\footnote{In Lobachevsky's word, (\selectlanguage{Russian} В Алгебре употребляют буквы, избегая, однако ж, способа бесконечно малых или границ \selectlanguage{English}). Russian (\selectlanguage{Russian} граница \selectlanguage{English}) in today's mathematical parlance stands for \emph{bound} or \emph{boundary} but the context suggests that Lobachevsky uses this word here in the sense of \emph{limit}.} Infinitesimal methods and limits require more mental efforts and thus constitute the last and the highest part of Analytics. Algebra is the same science that Newton called \emph{Universal Arithmetic}\footnote{See \cite{Newton:1707} and later editions.} in order to distinguish it from the numerical Arithmetic; the same discipline can be called the \emph{calculus of finites} to oppose it to the infinitesimal calculus, which indisputably involve new first principles.\footnote{Montucla similarly distinguishes between the \emph{ordinary} algebra that treats finite quantities and the \emph{infinitesimal} algebra \cite[p. 9]{Montucla:IV}.} [Various writers] present these new principles in different ways aiming at the [mathematical] rigour, which is an essential feature of every Mathematical doctrine. Lagrange achieved the rigour in his \emph{Calculus of Functions} [\cite{Lagrange:1797}]; it should be noticed, however, that the general theory of function does not reduce to the differential calculus.  It is much better to follow the example given by [Augustin--Louis] Cauchy  in his \emph{Algebraic Analysis} [\cite{Cauchy:1821}] and not deprive Algebra from the integrity that it may gain without using (what people call) the differential calculus. Number theory can be also included into Algebra but it should, nevertheless, be seen as a separate doctrine since it involves specific methods, which are not found elsewhere. 
 
I completed the [present] work  pursuing the aforementioned purpose. Its title [Algebra], which is nowadays commonly used, came to us from the Arabs along with the science itself. The Arabs [in their turn] borrowed this science in its infant state from the Greeks, namely from Diophantus' books.\footnote{Lobachevsky's lineage of Algebra from Diophantus appears naive but it is not wholly fantastic. Diophantus' \emph{Arithmetic} was first translated into Arabic later than al-Khw\={a}rizm\={i} accomplished his \emph{Algebra} ( by Qust\={a} b. L\={u}q\={a}) but the influence of Diophantus on al-Khwarezmi  cannot be simply ruled out \cite[p. 27-28]{Saliba:2007}. At the same time, Lobachevsky's presumption according to which the historical source of Algebra should be ultimately found in Greek Mathematics, is not justified.  The conceptualisation of Algebra by  al-Khw\={a}rizm\={i} was highly original and certainly could not be reduced to earlier contributions including the possible influence of Diophantus \cite{Rashed:1983}. A plausible source of this and the following Lobachevsky's historical notes is Jean-Etien Montucla's \emph{History of Mathematics} first published in 1758  \cite{Montucla:1758} and republished with certain revisions in an IV, i.e., in the period from September 23, 1795 to September 21, 1796, see \cite[p. 381--383]{Montucla:IV}. It should be noticed that pointing to Diophantus'  \emph{Arithmetic} as a possible source of Arab Algebra Montucla also stresses the originality of the Arab contribution.} The Eastern writers of that time added to the name ``Algebra'' word ``Mukabala'', meaning that solutions are grounded here on positions. \footnote{Name \emph{Algebra} derives from the title of a treatise by Muhammad ibn Musa al-Khwarizmi (c. 780--c. 850) written around 820 A.D.; The full Arabic title reads (in the standard transliteration to Latin alphabet) \emph{al-Kit\={a}b al-Mukhtasar f\={i} His\={a}b \textbf{al-Jabr} wal-\textbf{Muq\={a}balah}}. This treatise played a pivotal role in transferring algebraic concepts to Europe in the 13-15-th centuries, and thus a part of its title gave the name to the discipline.\cite{Rashed:2010}. Term \emph{al-Jabr} (that gave the name \emph{Algebra}) and term \emph{Muq\={a}balah} stand for syntactic operations involved in solutions of algebraic equations, which in the modern English are often rendered as \emph{Restoration} (or  \emph{Completion}) and \emph{Balancing}. Robert Chester who in 1145 produced the first Latin translation of  al-Khwarizmi's treatise rendered the title as \emph{Liber Algebrae et Almukabola}, i.e., transliterated it without providing full translation \cite{Karpinski:1915}. Later Luca Pacioli (aka Lucas de Burgo, 1445-1517), as quoted by Montucla, translated this title in Latin as \emph{restoratio et oppositio}   \cite[p. 382]{Montucla:IV}, which is the likely source of Lobachevsky's reference to ``position'' in this context. It makes perfect sense as long as \emph{position} is understood here in syntactic terms.} This distinguishing feature of Analytics appeared to European mathematicians so advantageous that they began to call it a \emph{strong art} (\selectlanguage{Russian} сильное искусство \selectlanguage{English}) ([in particular, it is called] \emph{Arte maggiore} by Pacioli and  \emph{Ars Magna} by Cardan).\footnote{Lobachevsky refers here to Pacioli's \emph{opus mangum} titled \emph{Summa de arithmetica, geometria, proporzioni di proporzionalita} first published in 1495, see modern English translation of a fragment in \cite{Pacioli:1994}. For some reason that I fail to explain Lobachevsky writes Pacioli's name as \emph{Patiolus}. Montucla writes it as \emph{Paccioli} \cite[p. 550]{Montucla:IV}. The second source referred to by Lobachevsky is \emph{Ars Magna} by Girolamo Cardano first published in 1545; see \cite{Cardano:1993} for a modern edition and translation into English. According to L.Ch. Karpinski, the name \emph{Ars Magna} and its variants derive from a remark in the \emph{Quadripartitum numerorum} by Johannes de Muris (c. 1290 -- c. 1355) who describes Arithmetic as a ``smaller'' art than that of Algebra. So some later authors including Pacioli and Cardano gave the name of \emph{Major Art} to Algebra and the name of   \emph{Minor Art} to Arithmetic \cite[p. 35]{Karpinski:1915}.} 

Algebra and Geometry had the same destiny. [\emph{Lobachevsky's footnote}: Hereafter we refer to Algebra without the help of Differential Calculus and to Geometry without Analytics.]  Rapid advances in the beginning [of the historical development of these sciences] were followed by very slow developments, which left these sciences at their [present] states, which are still far from being perfect. It likely happened because Mathematicians were primarily focused on higher parts of Analytics and neglected the foundations [of Algebra and Geometry]; they didn't want to cultivate this field, which they once already crossed and left behind. 

I hope that everyone agrees that in all Mathematical sciences the first principles are easily acquired but always have certain defects, which are very difficult to correct afterwards. The authors [of textbooks] who write for the beginners skip this point because they pursue a different aim and want to make their writings more accessible. At some point, however, one should return to the foundations and treat them rigorously. In my opinion, only Algebra brings to Mathematics the full generality and the conceptual rigour while Arithmetics is a preliminary step, which serves as a preparation and helps [students] to develop the [necessary computational] skill. This is why this Book begins with numbers, counting, and the four arithmetical operations. I judge it to be sufficient [for the purposes of the present Book] to treat positive and negative [magnitudes] as abstract [concepts] rather than borrow them from the nature of things directly; I treat these [magnitudes] as numbers artificially connected to their signs in order to introduce general rules for addition and subtraction.

Raising to power [as an arithmetical operation] gives [me] an occasion to consider functions and all dependences between numbers. The transition from the simple multiplication to exponents of any values serves [me] as the first example of this doctrine. The series for an exponent of the sum [of two magnitudes] known under the name of Newton's Binomial, was an important discovery in the Analytics. But util recently it lacked a fully general proof. It is based on the property of coefficients, which can be explained either by the method invented by myself or by Mr. Gauss' method (Ch. XI)\footnote{Hereafter Lobachevsky refers to chapters of his \emph{Algebra}  \cite[vol.4, p. 23-365]{Lobachevsky:1946-51}} The former method is preferable because it is a special case of expansion of roots of equation with three monomials (Ch. XV) found by Mr. [Johann Martin Christian] Bartels  (\emph{Vorlesungen \"uber Mathematische Analysis}, Dorpat 1833) and communicated to me [by him] while I was a student. The application of $\sqrt{-1}$ in Mathematics is so important that I judged it necessary to talk about imaginary powers more than this is usually done; I tried to provide [students] with clear concepts [of imaginary magnitudes] and lay dawn firm foundations for computations with help of  $\sqrt{-1}$. 

Solving equations was always the main subject of Algebra. It appeared that after the general expressions (called in this Book \emph{representative}) for equations of the first degree with any number of unknowns were found, there remained nothing else to be added to this subject. I managed, however, to replace such expressions by their proper values.\footnote{Lobachevsky's self-claimed contribution amounts, in modern terms, to a method of computing the determinant of a given square matrix using permutations of its elements, which is known in today's mathematical literature as  \emph{Leibniz formula for determinants}. The method using ``representatives''  (\selectlanguage{Russian} представительные выражения \selectlanguage{English}) is similar to that used by Cauchy in \cite{Cauchy:1815} but Lobachevsky's notation is quite different. } I expect that this discovery will not remain without use. I also included into [my] Algebra [a section on] the solution of linear equations in integers. Although treating this discipline in full generality is difficult and requires larger textbooks its foundations are necessary in Algebra and can be used in common problems like the problems of the Ecclesiastic Calculus (see Ch. X)\footnote{i.e., the problems of computing the dates of Orthodox Christian determined by multiple astronomical and other conditions, see  \cite[vol.4, p. 165-169]{Lobachevsky:1946-51}}

The increment and summing up of functions  (\selectlanguage{Russian} прирашение и суммирование функций \selectlanguage{English}) should also be taught as a part of Algebra because any calculus of finite magnitudes is a part of this discipline. The essential part of this content is found in Chapter XV; in addition to usual and well known examples the reader will find there a lot of new cases and further examples. This part of Mathematics deserves more attention than it has been given so far. \footnote{Chapter XV of Lobachevsky's \emph{Algebra} does not have counterparts in the contemporary Algebra textbooks. The ``increment and summing up of functions'' is an original attempt to develop finitary versions of Differential and Integral Calculi which could be useful for many practical purposes. This project is likely motivated by Lobachevsky's  scepticism about the introduction of infinitesimal and infinite quantities in mathematics; it squares with Lobachevsky's d'Alembertian epistemic strategy of closing the gap between Mathematics, Natural Sciences and Engineering. See \ref{sum7} for discussion.}

$\dots$

\section{Introduction and Chapter 1 of the \emph{New Foundations of Geometry with Complete Theory of Parallels} (1835)}\label{apH}  

\subsection{Introduction} \cite[vol.2, p. 158-167]{Lobachevsky:1946-51}\footnote{We include here only the concluding part of this Introduction in which Lobachevsky presents his philosophical ideas about Nature, perception, space, measurement and other general matters. In the begining of this Introduction Lobachevsky considers various historical attempts to prove Euclid's Fifth Postulate. An earlier full English translation of this \emph{Introduction} by G.B. Halsted has been published in 1897 \cite{Halsted:1897}; it is available at \url{http://philomatica.org/wp-content/uploads/2026/05/halsted1897-125-141.pdf}. The present translation from Russian is made by the author independently. There exists also a German translation of this \emph{Introduction} by Friedrich Engel published in 1898 \cite[p. 67--83]{Engel:1898} and available at \url{http://philomatica.org/wp-content/uploads/2026/05/engel1898-1-91-107.pdf}.}

In the Nature we properly grasp only motion without which no sensual perception is possible. All other concepts including geometrical concepts are produced in our minds artificially on the basis of properties of motion. It follows that space does not exist for us separately and by itself. So the assumption according to which some forces in Nature follow one geometry while some other forces follow another specific geometry, cannot produce a contradiction in one's mind. In order to clarify this thought let me suppose, following a widely shared belief, that gravitational forces weakens with their propagation along a sphere.\footnote{Lobachevsky refers here to the assumption according to which a gravitating point-like body $O$ exerts a gravitational force on a (simply connected) area of a sphere with center in $O$, which is proportional to the density of ``force lines'' propagating from $O$ onto the sphere. In case the space is Euclidean this assumption implies the inverse square law for gravitational (and any other central) forces.} In the \emph{Usual} geometry where the magnitude [ of the surface] of the sphere equals to $4\pi r^{2}$ (where $r$ is the half of the diameter [of the sphere]) the force should weaken according to the inverse square law. In the \emph{Imaginary} geometry I found that the surface [of the sphere] equals to 
$$\pi(e^{r} - e^{-r})^{2}$$. 
and it is possible that this [Imaginary] geometry governs molecular forces of various types, which are characterised by the value of [parameter] $e$\footnote{The value of $e$ in this formula depends on the unit of length used for measurements. In case the unit is the ``absolute'' unit equal to the absolute value of the \emph{radius of curvature} of the given Hyperbolic space then $e= \mathbf{e} = 2,71828\dots$ the base of natural algorithms. Generally, if the absolute value of the radius of curvature measured in arbitrary units equals to $R$ then we have $e = \mathbf{e}^{\frac{1}{R}}$.\label{hunit}} This is nothing but a hypothesis that needs to be confirmed with further evidence. It is, however, beyond a doubt, that motion, velocity, time, mass as well as distances and angles are produced by [physical] forces alone. Everything is related to forces, and since one doesn't know the essence of those forces, one is not in a position to claim that [mathematical] expressions that involve different types of magnitudes should involve only ratios of those magnitudes.\footnote{Hereafter Lobachevsky argues against Legendre's tentative prove of the Fifth Postulate using the ``homogeneity principle'', see footnote \footref{homogen} above. According to this principle angles of a given triangle cannot be possibly computed as functions of (lengths of) their sides.}  If we assume that [angle in a triangle] can depend on ratio [of its sides] why we cannot assume that it can depend on its sides directly? Some [physical] examples support this [latter] hypothesis. In particular, the magnitude of gravitational force is expressed by the mass divided by the square of distance [between gravitating bodies].\footnote{Lobachevsky points here to the fact that constant $G$ (called the \emph{Universal Gravitational Constant}) in Newton's Law of Gravity $F = G\frac{m_{1} \times m_{2}}{r^{2}}$ is not dimensionless: the standard dimensional analysis gives the physical dimension of $G$ equal to $M^{-1}L^{3}T^{-2}$ where $M$ stands for unit of mass, $L$ for unit of distance and $T$ for unit of time. Lobachevsky suggests that in the elementary geometry a similar situation may occur. Once again, it indeed occurs not only in the Hyperbolic geometry but also in the Spherical geometry.} 

If the distance [$r$] is equal to zero then this expression is meaningless. One should always begin with some actual [i.e., positive] distance whether it is large or small. Only in this case the force emerges. How the distance produces this force? How such heterogeneous objects [namely, the distance and the force] are related in the Nature? We will probably never learn it. But if forces depend on distances then [lengths of] lines may also depend on angles. The heterogeneity is similar in both these cases. The difference between the two cases is not conceptual but concerns only the fact that one [namely, the physical] dependence is learned from experiences while the other [namely, geometrical] dependence, given the shortage of [empirical] observations should be assumed intellectually either beyond the limits of the visible world or in the narrow realm of molecular forces.\footnote{Even if Lobachevsky's reference to a realm ``beyond the limits of the visible world'' naturally admits for a platonic metaphysical interpretation, the context suggests that he rather refers here to large astronomical scales inaccessible to the naked eye (as well as to the contemporary telescopes). Under this reading, the above Lobachevsky's argument amounts to saying that deviations from the Euclidean geometry can be possibly found in the physical world either at large astronomical scales or at microscopic ``molecular'' scales.}  

The hypothesis according to which the angles [of a given triangle] are determined by the ratios of the lengths  [of the three sides of this triangle] holds only in a special case, which is obtained when the sides [of the triangle] are taken to be infinitely small. Thus the Usual geometry always leads to correct conclusions. But these conclusions are not as general as those obtained in the generalised geometrical system, which I call the \emph{Imaginary Geometry}. The difference between equations [in the Usual and the Imaginary geometries] is due to the introduction of a new constant, which should be determined via [empirical] observations [i.e., measurements made in in the physical space].\footnote{Theoretically, this ``new constant'' can be identified with the (negative) curvature of Hyperbolic space. But since Lobachevsky does not apply the concept of curvature in his geometrical theories, in the given context the constant should be rather identified with $e$ in the above formula. See also footnote \footref{hunit}.\label{hunit2}} But since the observations do not detect any perceivable difference [between the two geometrical theories], the Usual geometry is sufficient [for all practical purposes] even if it is not strictly correct [with respect to the ``true'' geometry of the physical space, which remains out of the reach of the existing observational methods]. This means that such a system [i.e. the Usual (Euclidean) geometry] is either contingent in the Nature, or that all available [for human observations] distances are still infinitely small [in comparison with the true cosmological scales]\footnote{Lobachevsky wrote this words in the context of his contemporary rapidly developing observational astronomy, which progressively and dramatically extended the observable universe.}  
Generally, any proposition formulated in terms of distances, which is admissible in the Imaginary geometry, leads to rules of the Usual geometry when these distances are taken to be infinitely small.\footnote{Here we take the liberty to correct the original, which says just the opposite, namely, that assuming that in the equations of Imaginary (Hyperbolic) geometry the involved distances (lengths) are \emph{large}  one obtains the familiar equations of the Usual (Euclidean) geometry. This doesn't make mathematical sense, and moreover it overtly contradicts what Lobachevsky says in the preceding sentence and the following sentences. The mistake can be made by the note-taker who wrote down the text of \emph{New Foundations} while Lobachevsky dictated it, which somehow escaped the attention of later editors and translators including Engel in \cite[p. 77]{Engel:1898}. The original text reads \selectlanguage{Russian} ``[В]сякое положение, которое Воображаемая Геометрия допускает в элементах величины, будучи принято для линий в большом размере, должно необходимо приводить к правилам обыкновенной Геометрии, потому что с таким предположением удерживаются только первые степени тех чисел, которые представляют собой линии, а следовательно везде в уравнения войдут их содержания''  \selectlanguage{English}.  \cite[vol.2, p. 160]{Lobachevsky:1946-51}} This is because in this case one may ignore all powers of distances except the first power, so the equations will contain only ratios of the distances. Such are, for example the following propositions: [1] distances between two perpendiculars  [to the same straight line] are everywhere equal; [2] a perpendicular [$AB$ to given line $L$] describes a straight line by its end [$B$ when its base $A$ moves along line $L$]; [3] a circle transforms into a straight line when its diameter grows [unlimitedly]\footnote{Propositions [1]-[3] are logically equivalent to Euclid's Fifth Postulate}. 

Among all such known propositions the proposition that says that the ratio of lines [i.e., sides of a triangle] depends [only] on the angles [of the given triangle], has certain advantages. Indeed, the conceptual simplicity of this proposition is close to our first experiences. But this is all that can be said in its favour; any other judgement will be either false or unjustified. We shouldn't be confused by the idea that [in the Imaginary geometry] the immediate dependence of lines from angles brings along a constant, the value of which is as arbitrary as the choice of unit [of length]\footnote{See footnotes \footref{hunit2}, \footref{hunit}}. Against this argument [according to which to which admitting such an arbitrary constant is epistemologically unacceptable] we can respond that instead of considering ratios to a certain [arbitrarily chosen] line [taken for unit], we may rather consider ratios to a line, which is somehow determined in the Nature. I showed this in the Imaginary geometry by providing equations, where all lines stand in ratios with a single line, which could be found via [empirical] observations if the [precision of] observations would be sufficient for this purpose. 

Other similar propositions [equivalent to Euclid's Fifth Postulate] are either too artificial or arbitrary. Some attention deserves only proposition [3] concerning the transformation of a circle into a straight line. One can immediately see here a problem, which concerns a violation of continuity in the situation when a closed curve [namely, the given circle] should abruptly turn into an [open] infinite straight line thus loosing its essential property. The Imaginary geometry fills this gap more satisfactorily. In this geometry a growing circle where all its radii meet in its center finally [that is, in the limit] turns into a line [$H$] such that all straight lines which are normal to it approach each other unlimitedly albeit no longer intersect. Such a line  [$H$] having this property is no longer a straight line [like in the Euclidean geometry] but a curve that in my \emph{Foundations of Geometry} I call a \emph{limit circle}\footnote{See  \cite[vol.1, p. 195-196]{Lobachevsky:1946-51} and  \cite[vol.2, p. 289]{Lobachevsky:1946-51}. Lobachevsky also calls this line \emph{horocycle} (\selectlanguage{Russian} орицикл \selectlanguage{English},  \cite[vol.1, p. 106]{Lobachevsky:1946-51}). A two-dimensional counterpart of horocycle called a \emph{horosphere}, plays a very important role in Lobachevsky's presentation of his geometry: Lobachevsky observes that (in modern terms) the intrinsic geometry of horosphere (as a surface in the Hyperbolic 3-space) is Euclidean. This observation helps Lobachevsky to build a hyperbolic analogue of Euclidean Trigonometry, which becomes the formal analytic core of his novel geometrical theory. See \cite{Rodin:2015} and \cite{Braver:2011} for further details.}.

Finally, if the difficult problem of parallelism should be solved by [physical] experience then the experiment proposed by Legendre, namely, the 6-fold application of radius of a given circle to its circumference, should be judged insufficient. In my \emph{Foundations of Geometry} I established using astronomical observations that in the triangle two legs of which almost equal to the distance from the Earth to the Sun, the sum of [internal] angles cannot differ from the two right angles by more that 0,0003 seconds.\footnote{See \footref{astro1}.} This difference increases in the geometrical progression [i.e., exponentially] with the length of the legs, which shows that the Usual [Euclidean] geometry is sufficient for all actual measurements. One can come to the same conclusion using certain simple propositions appropriate in the foundations of a science but the complete theory requires to change completely the [usual] order of teaching; in particular, this concerns the Trigonometry. 

The imperfection of the [existing] theory of parallels concerns, among other things, the very definition of being parallel. But, contrary to Legendre, it doesn't depend either on the imperfection of the definition of straight line or on the defects of the primitive [geometrical] concepts, which I'm now going to highlight and try to fix. Usually the presentation of geometry begins with ascribing  three dimensions to solids, two dimensions to surfaces, one dimension to lines, and no dimension to a point. By giving to these dimensions the names of \emph{length}, \emph{width} and \emph{hight} and referring with these words to the three [spatial] coordinates, [the authors of geometry textbooks] harry to communicate immature concepts using words of common language, which don't have an exact scientific meaning. Indeed, how can one clearly think of measuring lengths without knowing what is a straight line? How can one speak of width and hight without saying anything beforehand about perpendiculars, planes, and the fact that that some perpendiculars belong to the same plane while some other perpendiculars  do not? Finally, if a point has no dimension, what makes it into a subject -mater of a judgement? Let us assume that everyone clearly represents to oneself a straight line without giving yet an account of this concept. Then there still remains the question of how to apply the straight line in order to attribute one dimension to a curve line and two dimensions to a curve surface.

Indeed, there is no need to require that length, width and hight be perpendicular to each other; it is sufficient to take for them lines going in different directions. But then new specific problems arise. How it is possible then, if one wants to avoid using concepts introduced only afterwards, to formulate the condition according to which the three sizes of a body belong to three straight lines in different planes? \footnote{Lobachevsky refers here to the condition that the three lines, which determine the three dimensions of a solid, should not be coplanar.} Further, different directions of the two parts of a broken line should not be confused with the double extension of the plane. Finally, one should wholly define what is meant by direction and by angle. In short, ``space'',  ``extension'', ``place'', ``body'', ``surface'', ``line'', ``point'', ``direction'', ``angle''  \textemdash\ are words with which the authors usually begin [their expositions of] geometry but which never convey a clear meaning. 

Meanwhile, these things can be also seen from a different perspective. It should be noticed that the conceptual obscurity is produced here by the abstractness, which is redundant in the actual measurements and  hence also unnecessary in the theory itself \label{abst}. Surfaces, lines and points as they are [usually] defined in Geometry exist only in our imagination while the [actual] measurement of surfaces and lines is done using bodies. This is why it is appropriate to talk about surfaces, lines and points only as so far as these concepts are used in the actual measurement. In this case we shall stick to concepts which are familiar to our imagination, and which can be verified in the nature straightforwardly without using other artificial and irrelevant concepts. In this case the science [of geometry] takes a different direction from the outset;  it follows this [new] direction until it becomes Analytics. Now I shall try to explain what this change amounts to.\footnote{See \ref{sum5}.}

In Mathematics one follows two methods, namely, Analysis and Synthesis. The distinctive feature of Analysis is the [presence of algebraic] equations, which serve as the first foundation for all judgements and lead to all conclusions. Synthesis or the method of constructions requires a [mental] representation, which is immediately associated with the basic concepts in our mind. The main advantage of Analysis is that it always provides for a direct path towards the intended goal. Synthesis is not a subject to general rules  but one should necessarily begin with this method in order to find the [appropriate] equations and thus to reach the point beyond which everything [in Geometry] becomes the science of numbers [i.e., Arithmetic]. For example, in Geometry one proves that [i] two perpendiculars [to the same straight line] do not intersect, that  [ii] triangles are equal in all their elements as soon as only some of those elements  are equal. It would be in vain to consider such cases, including the whole theory of parallels, analytically. This cannot be possible achieved. Similarly, Synthesis is indispensable in measuring plane [figures] bounded by straight lines, and in measuring solids bounded by planes. It goes without saying that in Synthesis one should apply [some elements of] Analysis; but it is beyond doubt that the foundations of Geometry and Mechanics cannot be built by analytic methods alone.  Geometry will always comprise properly geometrical [i.e., synthetic] features. One can narrow the circle of Synthesis but one cannot eliminate it from Geometry completely. In the desire to replace Synthesis by Analysis one should not harry to admit functions every time when one can assume a dependency [between magnitudes] without knowing what this dependency is and, moreover, how it can be expressed.\footnote{Lobachevsky argues here in favour of limiting the scope of Mathematics by ``reasonable'' functions, which are in some appropriate sense computable or, as he puts it  ``expressible''. Compare Lobachevski's misgivings about putative infinite mathematical objects in \ref{apINF}}.

These limitations of Analysis leave a room for Synthesis and help to understand its purpose. Only Synthesis can provide science with its basic concepts, which serve as foundations of all further judgements.  These further judgements are derived from the first data [i.e., the first principles provided by Synthesis]; thus they extend the horizon of our knowledge in all directions unlimitedly. The first data are, beyond any doubt, the concepts obtained in the Nature via our senses. The Mind can and should reduce their number to minimum in order they could serve as a firm foundation for the science. But it turns out that nobody follows the synthetic method as it is described here following the aforementioned rules but prefers to introduce Analysis prematurely relying on the imperfect development of the [relevant] concepts that constitute our natural mind, and naming those concepts without trying to explain them and provide them with exact definitions. Albeit this latter approach in teaching  [Geometry] is easy and simple, the rigorous truth has other advantages, which should not be ignored. My first attempt [in revising foundations of mathematics along these lines] has been made in Algebra \footnote{See  \emph{Algebra or Calculus of Finites}  (1834) \cite[vol.4, 23-27]{Lobachevsky:1946-51}, \ref{apG}}, and now I'm making a similar attempt in Geometry. 

A pure Analysis without an admixture of Synthesis cannot be appropriate in Geometry before every [geometrical] dependence is represented with equations, and geometrical magnitudes of all sorts are represented with [analytic] expressions. Every magnitude in Geometry should be understood in terms of measurement, which does not properly exist in the case of curve lines and surfaces. However small parts of a curve [line] one may take, they remain curve and thus cannot be measured with a straight line. A similar argument applies to the case of curve surface, any part of which, however small, never becomes flat. At the same time, we should not forget, that in the Nature there is no straight lines, no curve lines, no planes and no curve surfaces but there are only bodies; all the rest is the figment of our imagination and exists only in theory. Lagrange admits as a foundation the Archimedian proposition according to which two points on a curve can be always chosen sufficiently close, so the arc between them is larger than the chord but smaller than the sum of two tangents produced from those two points to the point of their intersection.\footnote{See \cite[p. 157]{Lagrange:1797}. The reference to Archimedes belongs to Lagrange. I did not immediately found this exact proposition in Archimedes' writings but in many occasions, in order to compute a length, an area or a volume, Archimedes considers its upper and lower approximations as in the given example.} This proposition is indeed necessary but it destroys the original idea of measuring curves with straight lines. The situation is similar in case of [curve] surfaces measured with planes. The computation of the length of a given curve or the computation of the area of a given curve surface is not, so to speak, a straightening of curvature. Such computations pursue a different goal, namely, finding a limit, which  the actual measurement approaches with the rise of its precision. This shows what justifies the Archimedian proposition. Along the pattern of curve lines one can also reason about the areas of given surfaces without claiming that its small parts are capable of rectification. 

For plane figures bounded by curve lines and for solids bounded by curve surfaces, strictly speaking, there is no measurement as far as they are supposed to be measured, correspondingly,  with a square and with a cube. Nevertheless, one may aim at finding a limit approached by the actual measurement. Then one needs to show that such a limit exists, represent the measuring procedure and explain how it may reach a desired precision. One cannot satisfy these requirements without making some special auxiliary assumptions, which are admitted as axioms: 1) that plane figures are equal when their parts are equal (even if in the two cases these parts are composed in different ways); 2) if one plane figure [$A$] wholly fits into another plane figure [$B$] but does not fill it completely then the former figure is smaller than the latter [$A<B$]; 3) that the magnitude [i.e., the area] of a given triangle vanishes when one of its sides vanishes. The last proposition is necessary in the [theory of] measurement of plane figures. Similar axioms are needed for the measurement of bodies. \footnote{As it has been shown by Max Dehn in 1901 polyhedra of equal volume may be not equicomposable, which shows that the theory of areas for plane figures sketched here by Lobachevsky does not, in fact, generalise to the 3D case \cite{Dehn:1901}.}

\subsection{Chapter 1} \label{apH1} \cite[vol.2, p. 168-170]{Lobachevsky:1946-51}
\textbf{First Concepts in Geometry} \footnote{German translation by F. Engel  \cite[p. 83--93]{Engel:1898} is available at \url{http://philomatica.org/wp-content/uploads/2026/05/engel1898-1-107-117.pdf}}

1) [Touch and Space]

The \emph{touch} is the distinctive feature of bodies, which gives them the name of \emph{geometrical} [bodies] when one focuses on this property ignoring all other properties [of bodies], both essential and contingent.\footnote{Compare \ref{touch1}.} 

In addition to bodies, possible subject-matters of [scientific] judgements include time, force, velocity of motion. The concept referred to by word ``touch'' does not belong to this list. We associate this concept only with bodies when we talk about their composition or about their division into parts. This [i.e., the concept of touch] is a primitive representation that we obtain in immediately in the Nature via our senses. It is not derived from some other concepts, and for this reason should not be a subject to further interpretation. 

When two bodies $A,B$ \emph{touch} each other they constitute one geometrical body $C$. In this case $A,B$ remain distinguishable constitutive parts of $C$; they do not lose themselves in the whole $C$.

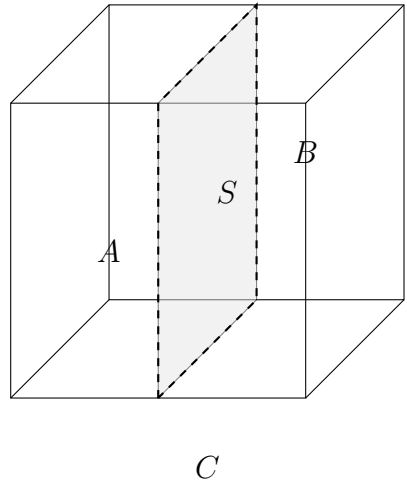
\begin{figure}[htp] 
\centering

\begin{tikzpicture}[scale=1.3]

\coordinate (P1) at (0,0);
\coordinate (P2) at (0,3);
\coordinate (P3) at (3,3);
\coordinate (P4) at (3,0);

\coordinate (P5) at (1,1);
\coordinate (P6) at (1,4);
\coordinate (P7) at (4,4);
\coordinate (P8) at (4,1);

\draw (P1)--(P2)--(P3)--(P4)--cycle;
\draw (P5)--(P6)--(P7)--(P8)--cycle;
\draw (P1)--(P5);
\draw (P2)--(P6);
\draw (P3)--(P7);
\draw (P4)--(P8);

\coordinate (M1) at ($(P1)!0.5!(P4)$);
\coordinate (M2) at ($(P2)!0.5!(P3)$);
\coordinate (M3) at ($(P5)!0.5!(P8)$);
\coordinate (M4) at ($(P6)!0.5!(P7)$);

\filldraw[fill=gray!20,opacity=0.5]
(M1)--(M2)--(M4)--(M3)--cycle;

\draw[dashed,thick]
(M1)--(M2)--(M4)--(M3)--cycle;

\node at (1.0,1.5) {$A$};
\node at (3.0,2.5) {$B$};

\node at (2.2,2.1) {$S$};

\node at (2,-0.7) {$C$};

\end{tikzpicture}

 \caption{Section $S$ of body $C$ into two bodies $A,B$}
\label{fig:section}
\end{figure}

Conversely, each body $C$ is divided into parts $A,B$ by arbitrary \emph{section} $S$ (Fig. \ref{fig:section}). Word ``section'' does not refer here to a new property of the given body [$C$] but only expresses the fact that $C$ is divided into two touching parts. Parts $A,B$ will be called \emph{sides} of section $S$ in body $C$.

All bodies in the Nature can be represented as parts of the same whole that we shall call the \emph{space}.\footnote{In his \emph{Foundations} of 1929 Lobachevsky qualifies the space as a body, see \ref{apF}. Here he doesn't make the same move. 
Lobachevsky's conception of space presented in the \emph{New Foundations} is basically Cartesian, which shares with the older Aristotelian relation conception the idea of the non-existence of a void. According to this conception the space is the totality of all spatial entities, i.e., of all bodies rather than a reservoir containing all bodies. This historical notion of Cartesian space stemming from works by Ren\'e Descartes should not be confused with Newtonian absolutist conception of space provided with Cartesian coordinates, which is more popular these days, see \cite{Slowik:2002},\cite{Zepeda:2014}.} 

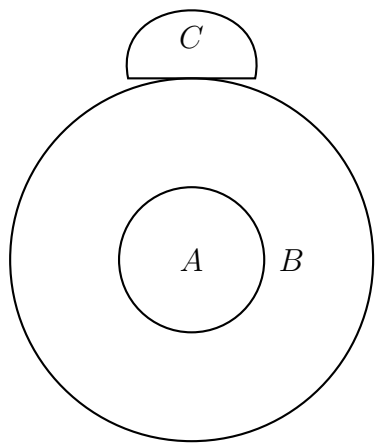
\begin{figure}[htp] 
\centering

\begin{tikzpicture}[scale=1.2]

\draw[thick] (0,0) circle (2);

\draw[thick] (0,0) circle (0.8);

\draw[thick]
(-0.7,2) 
.. controls (-0.9,3.0) and (0.9,3.0) .. 
(0.7,2)
-- cycle;

\node at (0,0) {$A$};
\node at (1.1,0) {$B$};
\node at (0,2.45) {$C$};

\end{tikzpicture}

 \caption{$B$ is the ambient space aka place of $A$: $A$ and $B$ cannot be jointly touched by another body $C$}
\label{fig:place}
\end{figure}

 2) [Place and measurement]
 Every given body $A$ can be touched by another body $B$ in such a way, that no third body $C$ can possibly touch both $A$ and $B$ (Fig.\ref{fig:place}). In this case, if one [additionally] assume that $B$ contains every third body $C$, $B$ is called the \emph{ambient space} of body $A$ \footnote{i.e., the ambient space of $A$ is the total space $U$ without $A$, in set-theoretic symbols $U\setminus A$.} If this further assumption is dropped, and a third body $C$ is allowed to be found outside of $B$ [as shown at Fig.\ref{fig:place}], $B$ is called a \emph{place} of body $A$ \emph{This is a geometrically articulated version of  Aristotle's  definition of \emph{place} found in his\emph{Physics}: ``Place is the primary motionless limit of that which surrounds [a body]'' (see. \cite[Book IV, chapter 4, 212a20–21]{Aristotle:2008}.} When $A$ touches $B$ in such a way that any body [$C$] which touches $A$ will also touch $B$ then we say that $A$ \emph{fills} its place $B$. All other bodies, which [may] also fill place $B$ without undergoing a change [during the replacement of one such body by another] will be geometrically \emph{identical} (\selectlanguage{Russian} одинаковы \selectlanguage{English}) in all respects to each other.\footnote{It is clear that by ``geometrical identity in all respects'' Lobachevsky means here congruence. Russian adjective \selectlanguage{Russian} одинаковый \selectlanguage{English}} in other contexts may mean ``indistinguishable'' or ``equal''.

Two bodies [$A,B$] are \emph{equal} when parts of one [say, part of $A$] can be rearranged in a new order, so the resulting body [$B'$] fills the place of the other body [$B$].
\footnote{Lobachevsky's definition of equality of bodies amounts to their \emph{equicompositionaly}, i.e. a possibility to compose two non-congruent given bodies $A,B$ from the same finite sets of parts arranging these parts in the two cases in different ways. Euclid develops a similar notion of equality for polygons \cite[Books 1-2]{Euclid:2008} (the difference is that he also counts as equal \emph{equicomplementable} polygons; this, however, is redundant from a logical point of view since every equicomposable Euclidean polygon is also provably equicomposable. Euclid's equality of polygons $A,B$ can be interpreted in modern terms as equality of their areas $S(A) = S(B)$ (since two polygons have the same area if and only if they are equicomposable). But this latter propositions does not generalise to the case of solids. Euclid's concept of equality of solids developed in \cite[Book 12]{Euclid:2008} can be understood as the equality of volumes but its foundations are different: in order to define the equality of solids Euclid applies his geometrical theory of proportions that covers the case of incommensurable magnitudes (\cite[Book 5]{Euclid:2008}). In 1901 Max Dehn showed that this feature of Euclid's \emph{Elements} cannot be fixed because no counterpart of  the theory of Books 1-2 can exist for volumes: he proved that polyhedra of equal volume are not necessary equicomposable \cite{Dehn:1901}. Thus Lobachevsky's definition of equality for polyhedra probably does not reach its intended purpose. One should take into account, however, Lobachevsky's view on the concept of non-commensurable geometrical magnitudes as an artificial mathematical figment having neither practical nor theoretical importance \ref{apDIFF}. So Lobachevsky might have here in his mind some approximating measuring procedure for solids and rightly believe that it be sufficient for grounding his intended theory of volumes.}

A measurements of given geometrical body [$A$] with another body [$M$], taken as a \emph{measure}, is accomplished when both bodies [$A,M$] are divided into equal parts [$P$]; then the \emph{magnitude} [of $A$] is expressed as ratio [$\frac{m}{n}$] where [$m$] is the number of parts [$P$] in [$A$] and [$n$] is the number of parts [$P$] in [$M$]. Thus the possibility to find the magnitude of bodies assumes the possibility of composing every body via the repetition of one [distinguished body], which is [repeatedly] applied to itself and in this way fulfils the measured body as well as the measure itself.\footnote{This distinguished body is thus equal to $\frac{1}{n}$-th part of measure $M$.} If this cannot be done strictly then one begins with fixing a certain type (\selectlanguage{Russian} вид \selectlanguage{English}) of measure, which fills the space unlimitedly, and divides it into congruent parts of arbitrary size [i.e., into some arbitrarily fixed units of measurements]. Then one can form any given body by composing [those parts of the measure, i.e., units] progressively and thus achieve the degree of equality [between the given measured body and the body composed of the units], where are senses cease to perceive imperfections.\footnote{By the ``imperfections''  (\selectlanguage{Russian} недостатки \selectlanguage{English}) Lobachevsky means here the difference between the real size of the measured body and its measured size. According to this account a measure is a rational number rather than real number.} In this way the errors of measurement will not exceed deviations occurring  in the Nature itself \label{nature}. The Nature is the principal goal of sciences and the source of our fundamental concepts. Yet the rigour in these concepts is due the imperfection of our senses [rather than the Nature itself]. Indeed, in the Nature there is no continuous composition, which implies that, contrary to how it appears for the whole at a larger scale,  the regularity of formation [of bodies from its constituents] is not preserved up to the smallest particles.\footnote{Lobachevsky believes that at the microscopic scales the Nature is made of  atoms (or molecules), see \ref{apH}. Given his relationalist account of space, this implies that the apparent continuity of physical space does not extend to the microscopic level. At this  respect Lobachevsky's view of physical space is similar to today's scientific view.}. On the other hand, when we produce in our imagination incommensurable bodies like a ball and a cube, we should agree about criteria of being larger and being smaller\footnote{i.e., criteria of comparison, which may allow one to claim, for example, that a given cube is \emph{larger} than a ball inscribed into this cube.}, then enclose the measured [body] into certain bounds\footnote{i.e., to provide for its volume an upper and a lower bounds}, which can be made arbitrarily close [to each other]. In this way we may treat the measurement of the mental [geometrical] objects with an arbitrary precision.\footnote{In the remaining part of Chapter 1 of the \emph{New Foundations} Lobachevsky presents a more detailed version of his theory of dimension, which has been earlier sketched in the \emph{Foundations}, see \ref{apF}}.            

$\dots$
\bibliographystyle{plain}
\bibliography{loba2}
\end{document}